\RequirePackage{fix-cm}
\documentclass[]{article}
\usepackage{arxiv}
\usepackage[english]{babel}
\usepackage{lipsum}
\usepackage{graphicx} 
\usepackage{subcaption}
\usepackage{algorithm}
\usepackage{algpseudocode}
\usepackage{placeins}
\usepackage[table]{xcolor}
\usepackage{bm}
\definecolor{mygray}{gray}{0.9} % Adjust the shade of gray as needed
\usepackage{booktabs} 
\usepackage{csquotes}
\usepackage{multirow}
\usepackage{listings}
\usepackage{xcolor}
\definecolor{abstract_background}{RGB}{235,235,235}%{251,231,201}
\usepackage{float}
\usepackage{color, colortbl}
\definecolor{ao(english)}{rgb}{0.0, 0.5, 0.0}
\usepackage{xr-hyper}
\usepackage{hyperref}
\hypersetup{colorlinks,linkcolor={ao(english)},citecolor={blue},urlcolor={red}} 
\usepackage{stackrel,amssymb}
\usepackage{amsfonts,amsmath,amstext,amsbsy,amsthm}
\theoremstyle{definition}

\newtheorem{definition}{Definition}[section]
\newtheorem{assumption}{Assumption}[section]
\newtheorem{theorem}{Theorem}[section]

\usepackage{cleveref}
\usepackage{amsmath}
\definecolor{LightCyan}{rgb}{0.88,1,1}
\definecolor{LightRed}{rgb}{1,0.7,0.7}

\title{Nonlinear Discrete-Time Observers with Physics-Informed Neural Networks}

\author{
\textbf{Hector Vargas Alvarez\textcolor{teal}{$^{1}$}, Gianluca Fabiani\textcolor{teal}{$^{1,2}$}, Ioannis G. Kevrekidis\textcolor{teal}{$^{2,5,6}$}, Nikolaos Kazantzis\textcolor{teal}{$^{3,}$}\thanks{Corresponding author, email: \texttt{nikolas@wpi.edu}}\hspace{0.1cm},}\\ \textbf{Constantinos Siettos\textcolor{teal}{$^{4,}$}\thanks{Corresponding author, email: \texttt{constantinos.siettos@unina.it}}\hspace{0.1cm}} \\
{}\\
\textcolor{teal}{$^{(1)}$} Modelling Engineering Risk and Complexity, \emph{Scuola Superiore Meridionale}, Naples 80138, Italy \hspace{1cm}\\
\textcolor{teal}{$^{(2)}$}Department of Chemical and Biomolecular Engineering, \emph{Johns Hopkins University}, Baltimore, MD 21218, USA\\
\textcolor{teal}{$^{(3)}$} Department of Chemical Engineering, \emph{Worcester Polytechnic Institute}, Worcester, MA 01609, USA\\
\textcolor{teal}{$^{(4)}$}Dipartimento di Matematica e Applicazioni ‘‘Renato Caccioppoli", \emph{Universit\`a degli Studi di Napoli}\\ \hspace{0.38cm}\emph{Federico II}, Naples 80126, Italy\\
\textcolor{teal}{$^{(5)}$}Department of Applied Mathematics and Statistics, \emph{Johns Hopkins University}, Baltimore, MD 21218, USA\\
\textcolor{teal}{$^{(6)}$}Department of Urology, School of Medicine, \emph{Johns Hopkins University}, Baltimore, MD 21218, USA\\}
\usepackage{amsopn}

\begin{document}

\maketitle

% REQUIRED
\begin{abstract}
\colorbox{abstract_background}{\begin{minipage}{1\linewidth}
We use Physics-Informed Neural Networks (PINNs) to solve the discrete-time nonlinear observer state estimation problem. Integrated within a single-step exact observer linearization framework, the proposed PINN approach aims at learning a nonlinear state transformation map by solving a system of inhomogeneous functional equations. The performance of the proposed PINN approach is assessed via two illustrative case studies for which the observer linearizing transformation map can be derived analytically.
We also perform an uncertainty quantification analysis for the proposed PINN scheme and we compare it with conventional power-series numerical implementations, which rely on the computation of a power series solution.
%It is demonstrated that the proposed PINN approach outperforms, in terms of numerical approximation accuracy, traditional/conventional numerical implementation approaches relying on the computation of a power series solution to the associated system of functional equations. 
%A physics-informed machine learning (PINN) approach is discussed in this investigation to address the observation problem associated %with discrete-time nonlinear dynamical systems. In this case, the PINN approach aims to learn a nonlinear transformation that %converts the nonlinear dynamics into linear dynamics driven by an injection term. Due to the fact the nonlinear transformation law %exhibits steep gradients due to singularities, we adopt a greedy-wise training strategy to alleviate the method's convergence. We %assess the performance of the proposed PINN approach via two different nonlinear discrete system systems whose nonlinear %transformation can be derived analytically.
%
%We evidence that the proposed PINN outperforms, in terms of numerical approximation accuracy, the established numerical %implementation, which involves the construction of a power series solution related to a system of functional equations. The great %importance of continuation techniques in the training procedure of PINN is especially remarked and suggested.
\end{minipage}
}
\end{abstract}

\keywords{Physics-Informed Neural Networks \and Nonlinear Discrete Time Observers \and Nonlinear Discrete Time Systems \and Uncertainty quantification}

% REQUIRED
%\textbf{\emph{AMS subject classification.}} 68T07, 68T20, 65L80, 65M22, 37N30, 65P99
%\newpage
\section{Introduction\label{sec:intro}}

In modern feedback control systems theory and practice, reliable access to the dynamically evolving system states is needed at both the implementation stage of advanced control algorithms and for process/system condition and performance monitoring purposes \cite{Isidori1995, Chen2013,Sepulchre2011,gelb1974applied,Wu2020}. Traditionally, an explicit use of an available dynamic model complemented by sensor measurements, involving measurable physical and chemical variables of the system of interest, represented a first option to respond to the above need. However, in practice, key critical state variables are often not available for direct on-line measurement, due to inherent physical as well as practically insurmountable technical and economic limitations associated with the current state of sensor technology as it is invariably deployed in cases of considerable system complexity \cite{Wu2020,Isidori1995,Chen2013, Sepulchre2011}. In light of the above remarks, a better, scientifically sound and practically insightful option is the design of a state estimator (an observer). This is itself an appropriately structured dynamical system itself that utilizes all information provided by a system model as well as available sensor measurements to accurately reconstruct the dynamic profiles of all other unmeasurable state variables \cite{Wu2020,Isidori1995,Chen2013,gelb1974applied}. It should be pointed out, that a state observer is endowed with the requisite design degrees of freedom to arbitrarily assign the rate of convergence of the state estimate it generates to the actual state or, equivalently, the dynamic modes of the induced estimation error dynamics.  Furthermore, since digital and Artificial Intelligence (AI) technology are increasingly integrated into advanced modern control and system performance monitoring strategies, the design of state estimators/observers with mathematical representations realized in the discrete-time domain (known also as ``software'' or ``soft'' sensors in state estimation) appears an amply motivated, theoretically and practically attractive line of contemporary research activity. Even though the well-known Kalman filter and its deterministic analogue, the Luenberger observer, represent comprehensive solutions to the linear discrete-time state estimator/observer design problem \cite{gelb1974applied, Chen2013}, the nonlinear discrete-time one is considerably more challenging and has attracted appreciable attention in the pertinent body of literature. For nonlinear systems, the conventional approaches have focused on the linearization of the dynamic equations around: i) either a reference equilibrium point and the subsequent implementation of a standard linear Luenberger observer, or ii) a reference trajectory followed by a standard Kalman filter design method (“extended Kalman filter”) \cite{Wu2020,Isidori1995, Chen2013, gelb1974applied}. However, both traditional approaches exhibit local validity within a small region in state space around the steady state, and could therefore lead to unacceptable observer performance since dominant system nonlinearities are conveniently ignored. Within the geometric system-theoretic framework, systematic discrete-time nonlinear observer design methodologies have been introduced \cite{ChungGrizzle1990, LeeNam1991, LinByrnes1995} in an effort to rigorously develop discrete-time analogues of the pioneering work by Krener and Isidori (1983) \cite{KrenerIsidori1983}, previously developed in the continuous-time nonlinear observer design problem. In all the aforementioned approaches, the authors seek \textit{for a nonlinear state transformation to linearize the system up to an additional output injection term}, followed by a second step whereby a linear Luenberger observer design method is applied to the transformed system (known as exact observer linearization methods). However, such nonlinear observer design methods, while theoretically rigorous and insightful, rely on  a set of rather restrictive conditions that can hardly be met by every physical or engineering system \cite{ChungGrizzle1990, LeeNam1991, LinByrnes1995}. Other interesting approaches, presented in \cite{Ciccarela1993, Boutayeb2000, Lilge1998, Califano2009}, propose nonlinear discrete-time observer design methods based on either restrictive global Lipshitz conditions, on past values of certain output variables and/or become applicable to special classes of nonlinear systems with mild nonlinearities.  Within the exact observer linearization framework, Kazantzis and Kravaris (2001) \cite{KazantzisKravaris2001} introduced a new nonlinear discrete-time observer design method through a single-step problem reformulation that no longer requires the redundant linearization of the system’s output map, thus effectively overcoming the restrictive conditions associated with the aforementioned approaches. In particular, it was shown that for linearly observable, real analytic nonlinear discrete-time systems whose spectrum of the linear part lies wholly in the Poincar\'e domain, the single-step nonlinear discrete-time observer design problem admits a solution that is based on a fairly general set of assumptions. In subsequent work by Xiao et al. (2003) \cite{Xiao2003} the above Poincar\'e requirement on the spectrum of the system's linear part was relaxed and replaced by an even milder assumption of the Siegel-type. Finally, recent important contributions further extended this approach to cases where relaxed regularity assumptions imposed on the system dynamic equations guarantee the existence of such an observer \cite{Brivadis2019,Venkat2022}.
%In real applications, all state variables are rarely accessible for direct online assessment in practice. Most of the time, a %trustworthy estimation of the unmeasurable state variables is highly necessary, especially when they are employed to create model-%based output feedback controllers \cite{Christophides,Isidori1995} or for process monitoring. Thus, a state observer is typically %used for this precise assignment in order to precisely reconstruct the unmeasurable state variables. Both the well-known Kalman %filter \cite{Eykhovt,gelb1974applied} and its deterministic analogue realized by Luenberger's observer design theory \cite{Luenberger1963} %provide an exhaustive solution to the problem for continuous or discrete-time linear systems. The nonlinear observer design problem &is substantially more difficult in the field of nonlinear systems and has attracted a lot of interest in the literature. There have &been numerous attempts to build continuous-time nonlinear observer design techniques. Krener and Isidori \cite{Krener1999} offered a &novel method for the construction of a theory of nonlinear observers.
%The authors' groundbreaking work involved linearizing the original system up to an additional output injection term using a nonlinear %state transformation. The development of a systematic nonlinear observer design technique for discrete-time systems, which is also %based on a fairly comprehensive set of requirements, was proposed by Kazantzis and Kravaris \cite{NikoKat} based on the concepts of %Luemberger.

Recent theoretical and technological advances have reignited the control community's interest in the development of new machine learning-based schemes  by revisiting techniques developed back in the 90s \cite{nguyen1990neural,hunt1992neural,miller1995neural,chen1995adaptive,taprantzis1997fuzzy,lee1997application, siettos1999advanced,Gonzalez-Garcia1998identification, kim1999neural,leu1999observer,alexandridis2002modelling,siettos2002truncated,siettos2002semiglobal} and introducing new ones \cite{Wu2019,wu2019machine,Wu2020,Wu2022,lombardi2021using,lombardi2021dynamic,alvarez2023discrete,patsatzis2023data} to name just a few. Regarding the design of nonlinear observers, a recurrent neural network (RNN) structure integrated within the extended Luenberger observer design framework has been recently proposed in \cite{alhajeri2021machine} in an effort to develop data-driven state estimators used for Lyapunov-based Model Predictive Control (MPC) purposes. In Janny et al. (2021) \cite{Janny2021}, a machine learning approach for the continuous-time nonlinear observer design problem that only depends on a few measured system trajectories has been suggested; in this approach, deep networks and gate recurrent units (GRU) are deftly employed. On a similar line of inquiry, for nonlinear systems represented through neural ODEs,  a procedure to develop state observers with embedded machine learning capabilities has been introduced in Miao and Gatsis (2023) \cite{miao2023learning} being structurally reminiscent of Luenberger ones. In Abdollahi (2006) \cite{Abdoll}, an adaptive observer design method for a broad class of multiple input-multiple output (MIMO) nonlinear continuous-time systems modeled through appropriately parameterized neural network structures has been proposed. Furthermore, another observer design method for the rotor resistance of an indirect vector controlled induction motor drive is presented in Karanayi (2006) \cite{Karanayil2006}, that utilizes  a fuzzy logic-based  observer, leading to further enhancement of the performance characteristics of the artificial neural network (ANN) structure. The results of the first rigorous research study aiming at integrating neural networks into the continuous-time single-step exact observer linearization framework introduced in Kazantzis and Kravaris (1998) \cite{kazantzis1998} were presented in Miao and Gatsis (2023) \cite{miao2023learning}. In particular, the linearizing transformation map was identified through a supervised learning method, which may occasionally suffer by limitations in the exploration of the full state space during the training stage. Recently, in another research study whose primary aim is to overcome these limitations, an unsupervised deep learning neural network method has been introduced in Peralez et al. (2022) \cite{Peralez2022}, and integrated into the discrete-time single-step exact observer linearization framework introduced in Kazantzis and Kravaris (2001) \cite{KazantzisKravaris2001}. It should be pointed out that deep learning neural network structures employed for learning nonlinear maps have been evaluated in the pertinent literature as well \cite{lusch2018deep}. Within a conceptually similar nonlinear discrete-time observer design context, in Peralez et al. (2022) \cite{Peralez2022}, a new method is proposed that leads to an enhancement of the numerical approximation of the map and its extension in the transient phase with an ANN architecture.

Motivated by previous work on the feedback linearization of discrete time maps using PINNs in Alvarez et al. (2023) \cite{alvarez2023discrete},  PINNs are used in the present research study to solve the nonlinear observer design problem in the discrete-time domain within the single-step exact observer linearization framework originally introduced in Kazantzis and Kravaris (2001) \cite{KazantzisKravaris2001}. In particular, we use PINNs to learn the state transformation map that linearizes the observer dynamic equations up to an output injection map. The PINN-based discovery of the above map, that in the original single-step problem formulation has to satisfy a system of first-order inhomogeneous functional equations, enables the design of a discrete-time observer with linearizable error dynamics and assignable convergence rates under a set of rather mild assumptions. 
%Compared to previous studies based on ANNs and attendant ML methods, where solutions to neural ODE systems are sought that are directly dependent on the parameters of the underlying network, the present research work differentiates itself by placing its focus on training the Neural Network (NN) in a physics-informed manner. Indeed, by following such an approach a considerable improvement of the network's convergence and numerical approximation precision properties is realized in comparison to the Taylor series expansion method suggested in \cite{KazantzisKravaris2001} for solving the associated set of functional equations and deriving the linearizing transformation map. In this case, 
%the coefficients of the Taylor series representation of the transformation map can be calculated recursively
%by solving a system of linear algebraic equations and utilizing a symbolic software program such as MAPLE \cite{KazantzisKravaris2001}. However, even for medium-scale system dimensionality, such a power-series expansion approach could become computationally impractical and thus unable to attain the requisite numerical approximation precision levels {\em in the entire domain}, particularly in state space regimes with gradients that are extremely steep and occasionally exhibit singularities. 
In particular, for the training process, a step/greedy-wise training method is employed in Alvarez et al. (2023) \cite{alvarez2023discrete} to enhance the PINN scheme's ability to  numerically approximate nonlinear maps in the presence of steep gradients. %Within such a context,  one begins by learning the nonlinear transformation {\em on a subset of the entire domain} where the transformation is sought, gradually increasing its size. This task is known as "warm-restarting" the training procedure, which uses the results derived from the previous algorithmic step as initial assumptions for the ANN's unknown weights. Notice that the proposed PINN scheme is quite simple to implement, and indeed, a "homemade" pertinent code was created on Matlab 2022. For illustration and training purposes the Levenberg-Marquard optimization algorithm was wrapped around it using the  \texttt{nonlinsq} function. A Python version of the TensorFlow library's Keras API was developed for comparison purposes as well. In order to illustrate and evaluate the aforementioned computational procedure two benchmark case studies involving nonlinear discrete-time dynamic system models were considered. Their associated observer linearizing transformation maps can be derived analytically and therefore they were deemed suitable for assessing the performance and effectiveness of the proposed continuation/greedy-wise PINN method. In particular, both benchmark models considered in this work exhibit a singularity at the domain boundary, thus making it difficult to approximate the respective maps numerically close to that point, especially through a power-series expansion method. For illustration purposes, we also compared the numerical approximation accuracy of the proposed PINN greedy scheme against the standard power-series expansion method, as well as Matlab and Python's TensorFlow-based implementations with automatic differentiation that were used to learn the transformation in the entire domain. Furthermore, we considered two different scenarios, namely: (a) one where  we assumed that the equations of the model are explicitly known, and (b) one where we assumed that only a black-box simulator is available, i.e., pertinent equations are not available explicitly in closed form.\par

The paper is organized as follows: in Section \ref{sec:framework}, we present a brief overview and the necessary mathematical preliminaries associated with the single-step exact observer linearization framework in the discrete-time domain, followed by a description of the structure of the proposed PINN scheme.  In Section \ref{sec:benchmarks}, the two benchmark case studies are introduced. The numerical results derived from the proposed PINN-methods are presented in Section \ref{sec:numerical_results}, along with a discussion on a comparative performance assessment against the classical power series expansion method. Finally, concluding remarks are offered in Section \ref{sec:conclusion}.

\section{Preliminaries and Methodological Framework}
\label{sec:framework}
We consider nonlinear discrete-time autonomous systems admitting the following state-space representation:
\begin{gather}
x(t + 1) = \Phi(x(t)) \nonumber \\
y(t) = h(x(t)),
\label{eq:DiscreteSystem}
\end{gather}
where  $t=0,1,2,...$ is the discrete-time index,  $x(t) \in \mathbb{R}^n $ is the state variables vector and $y(t)  \in \mathbb{R}$ represents the system's measured output
variable. It is also assumed that $\Phi(x)$ is a real analytic vector function defined on $\mathbb{R}^n$ and $h(x)$ is a real analytic scalar function on $\mathbb{R}^n$. Furthermore, without loss of generality, let  the origin $x = 0$ be an equilibrium point of (\ref{eq:DiscreteSystem}): $\Phi(0) = 0$ with $h(0) = 0$. The following
assumptions are made:
\begin{assumption}
\label{ass:assumption1}
   The Jacobian matrix $F$ of $\Phi(x)$ evaluated at $ x = 0: F = (\frac{\partial \Phi}{\partial x})(0)$ has eigenvalues $k_i,\quad i\in \mathbb{N}$ that all lie in the Poincar\'e domain. 
\end{assumption}

\begin{assumption}
\label{ass:assumption2}
    Denoting by $H$ the
$1 \times n$ matrix: $H=\displaystyle{\frac{\partial h}{\partial x}}(0)$, 
it is assumed that the following $n \times n$ matrix $O$:
\begin{equation}
O=\begin{bmatrix}
H\\
HF\\
\vdots\\
HF^{n-1}
\end{bmatrix}
\end{equation}

has rank $n$. This assumption states that system in Eq. \eqref{eq:DiscreteSystem} 
is locally observable around the origin $x=0$ \cite{Chen2013, Kazantzis2001}.
\end{assumption}

A nonlinear discrete-time state observer driven by the available measurements of the output variable $y(t)$ and being capable of reconstructing the vector of state variables of system in Eq. \eqref{eq:DiscreteSystem} conforms to the structure delineated in the definition below \cite{Isidori1995, KazantzisKravaris2001}:
\begin{definition}
\label{def:def1}
    A dynamical system 
\begin{gather}
 z(t + 1) = \Psi(z(t),y(t))   
 \label{eq:GeneralFormLinEquation}
\end{gather}
with $z \in \mathbb{R}^n$; $y\in \mathbb{R}$; $\Psi:\mathbb{R}^n \times \mathbb{R} \rightarrow \in \mathbb{R}^n$, is called a \emph{full-order observer} for Eq. \eqref{eq:DiscreteSystem}, if there exists a locally invertible
(around the origin) map $T(x)$ with $T: \mathbb{R}^n \rightarrow \mathbb{R}^n$ and $T(0) = 0$, such that if $z(0) = T(x(0))$ with $x(0)$ near
zero, then $z(t) = T(x(t)), \quad \forall t>0$. In particular, when $T$ is the identity map,  Eq. \eqref{eq:GeneralFormLinEquation} is called the \emph{identity observer} for Eq. \eqref{eq:DiscreteSystem}.
\end{definition}

Definition \ref{def:def1} implies that if $\Psi$ and $T$ are related through the following system of functional equations:
\begin{gather}
  T(\Phi(x))= \Psi (T(x),h(x)) \nonumber \\
T(0) = 0,
\label{eq:ConditionFunequations}
\end{gather}
then the $y$-driven dynamical system \eqref{eq:GeneralFormLinEquation} represents an observer for system \eqref{eq:DiscreteSystem}, whenever the map $T(x)$ is invertible. It is important to note that the initial condition $T(0)=0$ that goes along with the system of  functional equations \eqref{eq:ConditionFunequations} reflects the fact that equilibrium properties are preserved under the proposed state transformation. Please also notice that in the special case of an identity observer, condition \eqref{eq:ConditionFunequations} reduces to
\begin{gather}
\Phi(x) = \Psi(x, h(x)).
\label{eq:identityobsever}
\end{gather}
Condition \eqref{eq:identityobsever} expresses the well-known mathematical requirement that the standard definition of a nonlinear observer entails \cite{Krener1999,LinByrnes1995, Sepulchre2011}.
Please notice that the vector function $\Psi(z,y)$ in the observer dynamic equations \eqref{eq:GeneralFormLinEquation} can be arbitrarily selected, as long as the system of functional equations \eqref{eq:ConditionFunequations}
admits an invertible solution $T(x)$. Furthermore, any stability requirement can be imposed at the observer design stage, where $\Psi$ will be
\emph{a priori} selected to ensure asymptotic stability of the observer \eqref{eq:GeneralFormLinEquation}. For simplicity reasons, it would be more practical to request linear observer dynamics in the
transformed states $z$ up to an output injection term by choosing:
\begin{equation}
\Psi(z,y)=Az+b(y),
\end{equation}
where $A$ is a constant matrix and $b$ a vector function defined on $\mathbb{R}$ with appropriate dimensionalities. Please notice that the 
asymptotic stability of the observer in Eq. \eqref{eq:GeneralFormLinEquation} can now be enforced by appropriately selecting the eigenvalues of the matrix $A$. Under such a requirement, the
map $T(x)$ of Definition \ref{def:def1} ought to satisfy the following system of
first-order non-homogeneous functional equations:
\begin{gather}
  T(\Phi(x))=AT(x)+b(h(x)) \nonumber \\
T(0)=0,
\label{eq:LFEs}
\end{gather}
with $T: \mathbb{R}^{n} \longrightarrow \mathbb{R}^{n}$ being the solution of \eqref{eq:LFEs}.  Within the above single-step exact observer linearization method, the following Theorem ensures the existence of a nonlinear discrete-time observer in the sense of Definition \ref{def:def1} \cite{KazantzisKravaris2001}:
\begin{theorem}
\label{thm:thm1}
Suppose that Assumptions \ref{ass:assumption1} and \ref{ass:assumption2} hold true. Let $B=\displaystyle{\frac{\partial b}{\partial y}}(0)$ and assume that the pair of matrices
$\{A,B\}$ is chosen to be controllable and the eigenvalues
$k_{i}, (i=1,...,n)$ of matrix $F=\displaystyle{\frac{\partial \Phi}{\partial x}(0)}$ are not related to the
eigenvalues $\lambda_{i}, (i=1,...,n)$ of matrix $A$ through any equation of the form:
\begin{equation}
\prod_{i=1}^{n}k_{i}^{m_{i}}=\lambda_{j}
\label{eq:cond1}
\end{equation}
$(j=1,...,n)$, where all the $m_{i}$'s are non-negative integers that satisfy the condition:
\begin{equation}
\sum_{i=1}^{n}m_{i}>0
\label{eq:cond2}
\end{equation}
Then, there exist a unique analytic and locally invertible solution $T(x)$ to the system of functional equations \eqref{eq:LFEs} as well as a nonlinear map $z=T(x)$ that make the dynamical system:
\begin{equation}
z(t+1)=Az(t)+b(y(t))
\label{TransfObserv}
\end{equation}
an observer 
for system \eqref{eq:DiscreteSystem} in the
sense of Definition \ref{def:def1}.
\end{theorem}

Notice that the proposed nonlinear discrete-time state observer expressed in the original state variables attains the following form:
\begin{equation}
\hat{x}(t+1)=T^{-1} \big[ T(\Phi(\hat{x}(t)))+b(y(t))-b(h(\hat{x}(t))) \big]
\label{OrigCoordObserv}
\end{equation}
where $\hat{x} \in \mathbb{R}^{n}$ is the state estimate and $T(x)$, $T :\mathbb{R}^{n} \longrightarrow \mathbb{R}^{n}$ is the solution to the associated system of
 functional equations \eqref{eq:LFEs}. The structure of the above observer is mathematically realized through the following two terms: 

(i) $\Phi(\hat{x})$ is a replica of the system dynamics \eqref{eq:DiscreteSystem}, and

(ii) $b(y(t))-b(h(\hat{x}(t)))$ (which drives the observer dynamics in Eq. \eqref{OrigCoordObserv}) is a feedback term that accounts for the mismatch between the actual sensor measurement $y(t)$ and its estimate $h(\hat{x}(t))$.

Furthermore, the proposed state observer in Eq. \eqref{TransfObserv}
induces linear estimation error dynamics with assignable dynamic modes in the transformed variables $z=T(x)$:
\begin{gather}
e_{z}(t+1)=z(t+1)-\hat{z}(t+1)= T(x(t+1))-T(\hat{x}(t+1))=\nonumber \\
=T(\Phi(x(t)))-T(\Phi(\hat{x}(t)))-b(y(t))+b(h(\hat{x}(t)))=\nonumber \\
=AT(x(t))+b(h(x(t)))-AT(\hat{x}(t))-b(h(\hat{x}(t)))-b(y(t))+b(h(\hat{x}(t)))=\nonumber \\
=AT(x(t))-AT(\hat{x}(t))=A(z(t)-\hat{z}(t))=Ae_{z}(t)
\end{gather} 
Therefore, if the matrix $A$ is chosen with an eigenspectrum ensuring stability, i.e. a set of eigenvalues located inside the unit disc on the complex domain, the magnitude of these eigenvalues will also directly regulate the rate of decay of the estimation error to zero, or equivalently the convergence rate of the state estimate to the actual state. Notice that the local invertibility of the transformation map $T(x)$ would imply that the state estimate $\hat{x}$ asymptotically approaches the actual state $x$.

Finally, it should be pointed out, that the above exact observer linearization approach is fundamentally different from the ones presented
in other studies (see e.g., in \cite{ChungGrizzle1990, LeeNam1991, LinByrnes1995}) where as a first step, a nonlinear coordinate transformation that transforms the original system  \eqref{eq:DiscreteSystem}
into a linear one (up to an 
output injection term) with
a linear output map is sought:
\begin{gather}
  z(t+1)=Az(t)+b(y(t)) \nonumber \\
y(t)=cz(t)
\label{LinTransf}
\end{gather}
The design of the observer in \cite{ChungGrizzle1990, LeeNam1991, LinByrnes1995} is then completed in a second step where a standard linear Luenberger observer is used 
for the transformed system \eqref{LinTransf}:
\begin{equation}
\hat{z}(t+1)=A\hat{z}(t)-L(y-c\hat{z}(t))+b(y(t)),
\end{equation}
with $L$ being the Luenberger observer gain and $\hat{z}$ the state estimate in the transformed coordinates.The
estimate $\hat{x}$ of the original state vector $x$ is then recovered, through the
inverse transformation employed in the first step of the design procedure. Please notice that the 
requirement of a linear output map in the 
transformed system \eqref{LinTransf} introduced in the above two-step approach, represents the main mathematical reason
for the emergence of quite restrictive conditions in \cite{ChungGrizzle1990, LeeNam1991, LinByrnes1995}. Within the 
proposed exact observer linearization framework however, the state observer is a dynamical system that is driven by the measured output variable and designed through an appropriate
coordinate transformation in a single step without the redundant requirement of linearity of the output map imposed on the transformed linear system \eqref{LinTransf}. 

A practical implementation of the proposed nonlinear discrete-time observer design method requires the development of a solution scheme for the nonlinear functional equations \eqref{eq:LFEs}. 
As mentioned earlier, the functions 
$\Phi(x)$, $h(x)$, as well as the transformation map $T(x)$
are all locally analytic. Therefore, as suggested in \cite{KazantzisKravaris2001}, a solution method could be based on a multivariate Taylor series expansion of $\Phi(x)$, $h(x)$ and $T(x)$, followed by a procedure that equates same-order Taylor coefficients of both sides of functional equations \eqref{eq:LFEs}. As a result, recursive
algebraic formulas can be generated that are linear with respect to the Taylor coefficients of the unknown solution. Indeed, one can calculate the $N$-th order Taylor coefficients of the $T(x)$ map given the Taylor coefficient values up to the order $N-1$ already calculated in previous recursive steps. The above linear recursive formulas admit a compact mathematical representation if tensorial notation is used (for more details, see \cite{KazantzisKravaris2001,siettos2006,alvarez2023discrete}).

The linearization of Eq. \eqref{eq:DiscreteSystem} around the equilibrium $(0,0)$ gives:
% \begin{equation}
%     x(t + 1) = \Phi(x(t))
% \end{equation}
\begin{equation}
dx(t+1)=\frac{\partial \Phi}{\partial
x}(0) dx(t);
\label{eq:linearizedsystem}
\end{equation}
on the other hand, since $\Phi(x(t))=x(t+1)$, the linearization of the transformed system \eqref{eq:LFEs} around the equilibrium reads:
\begin{equation}
%\textcolor{red}{T(\Phi(x))=AT(x)+b(h(x)) \nonumber \\
    \frac{\partial T}{\partial x}(0)dx(t+1)=\biggl(A \frac{\partial T}{\partial x}(0)+b\frac{\partial T}{\partial x}({h}(0))\frac{\partial y}{\partial x}(0)\biggl) dx(t).
    \label{eq:lin_transformed_around_equilibrium}
\end{equation}
Multiplying both sides of Eq.(\ref{eq:linearizedsystem}) by $\frac{\partial T}{\partial x}(0)$, the following is obtained:
\begin{equation}
  \frac{\partial T}{\partial x}(0)dx(t+1)=\bigg(\frac{\partial T}{\partial x}(0)\frac{\partial \Phi}{\partial x}(0)\bigg)dx(t).
   \label{eq:linearizedsystem2}
\end{equation}
Consequently, from Eqs.\eqref{eq:lin_transformed_around_equilibrium},\eqref{eq:linearizedsystem2}, it can be deduced that the nonlinear transformation centered on the equilibrium point $(x_0)$ must meet the following (phase) condition:
\begin{equation}
\frac{\partial T}{\partial x}(0) \frac{\partial \Phi}{\partial x}(0)= A\frac{\partial T}{\partial x}(0) + b\frac{\partial T}{\partial x}(h(0))\frac{\partial y}{\partial x}(0).
\label{eq:pinning}
\end{equation}
% \begin{equation}
% \frac{\partial T}{\partial x}(0) \frac{\partial \Phi}{\partial x}(0) - A\frac{\partial T}{\partial x}(0) = b\frac{\partial T}{\partial x}(0)\frac{\partial y}{\partial x}(0).

% \end{equation}
Notice that the elements of the Jacobian matrix $\frac{\partial T}{\partial x}(0)$ can be computed, by solving the system of equations \eqref{eq:pinning}, if $\Phi$ and $h$ are explicitly available. \par
Here, exploiting the properties of universal approximation \cite{chen1995universal}, we use a PINN to learn the transformation $T(x)$. The proposed PINN has two hidden layers and $N_1$, $N_2$ neurons for the first and second hidden layers respectively, along with a linear output layer. This can be written as:
% \begin{equation}
%     \hat{T}_j(x)=\sum_{i=1}^{N_2} W^{o}_{ij} \phi^{(2)}_i\biggl(\sum_{s=1}^{N_1}W^{(2)}_{si}\phi_s^{(1)}\biggl(\sum_{k=1}^n W_{ks}^{(1)}x_k+\beta^{(1)}_s\biggr)+\beta^{(2)}_i\biggr) +\beta^{(o)}_j,
% \end{equation}
% or equivalently, in matrix form:
\begin{equation}
   \hat{T}(x) = {W^{(0)}}^T \Phi_2({W^{(2)}}^T \Phi_1({W^{(1)}}^T x+\beta^{(1)})+\beta^{(2)})+\beta^{(0)}.
   \label{eq:FNN}
\end{equation}
$W^{(0)}$ is an $N_2 \times n$ matrix that encompasses the weights $W^{(0)}_{ji}$ linking the second hidden layer to the linear output layer. The multivariate vector-valued functions $\Phi_1:\mathbb{R}^n \rightarrow \mathbb{R}^{N_1}$ and $\Phi_2:\mathbb{R}^{N_1} \rightarrow \mathbb{R}^{N_2}$ are associated with the activation functions $\phi^{(1)}_s$ and $\phi^{(2)}_i$ in the first and second hidden layers, respectively.
$W^{(1)}$ is a $n \times N_1$ matrix that encompasses the weights connecting the input to the first hidden layer, while $W^{(2)}$ is a $N_1 \times N_2$ matrix whose elements are the weights linking the first hidden layer to the second hidden layer. The biases of the nodes in the first and second layers are denoted by $\beta^{(1)}\in \mathbb{R}^{N_1}$ and $\beta^{(2)}\in \mathbb{R}^{N_2}$, respectively; $\beta^{(0)} \in \mathbb{R}$ encapsulates the bias(es) of the output node(s).\par
To learn $T(x)$ with PINNs, we considered a certain domain $D\subset\mathbb{R}^n$ around the equilibrium point $(0,0)$ discretized in a grid of $M$ points $x_i$ with $i=1,\dots,M$. Thus, finding  $T(x)$ reduces to the task of minimizing the loss function:
\begin{gather}
\mathcal{L}(P)=\sum_{i=1}^M \sum _{j=1}^n {r^{(1)}_{ij}}^2(x_{ij},\hat{T}(x_{ij};P))+
\sum _{j=1}^n {r^{(2)}_j}^2( \hat{T}_j(0;P))+
\sum_{j=1}^n \sum _{k=1}^n {r^{(3)}_{jk}}^2\biggl(\frac{\partial\hat{T}_{j}}{\partial x_k}(0;P)\biggr),
\label{eq:loss1}
\end{gather}
with respect to the unknown parameters $P =({W}^{(0)},{W}^{(2)}, {W}^{(1)}, {\beta}^{(0)},{\beta}^{(2)},{\beta}^{(1)})$ of the FNN  given by (\ref{eq:FNN}). In the above:
\begin{gather}
r^{(1)}_{ij}(x_i, \hat{T}(x_i))= \hat{T}_j\bigg(\Phi(x_i) \bigg)-\biggl( {\alpha_j} \hat{T}(x_i) + b_i(h(x_i))\biggl), \qquad i=1,\dots,M, \quad j=1,2,\dots,n,
\label{eq:NFE_discretized}
\end{gather}
where ${\alpha_j}$ is the $j$-th row of the matrix $A$ and $\hat{T}_j$ is $j$-th output component of $\hat{T}$, and:
\begin{gather}
r^{(2)}_j(\hat{T}_{j}(0))= \hat{T}_{j}(0), \quad j=1,2,\dots,n,
\nonumber \\ r^{(3)}_{jk}(x_{k}, \hat{T}_{j}(0))= \frac{\partial \hat{T}_j}{\partial x_k}(0)-\frac{\partial {T}_j}{\partial x_k}(0), \quad j,k=1,2,\dots,n,
\end{gather}
where $\frac{\partial T_j}{\partial x_k}(0)$ is the $(j,k)$-th element of the Jacobian matrix of $T(x)$ computed at the equilibrium and obtained by solving the system of equations in (\ref{eq:pinning}). It should be noted that in the illustrative examples considered below, we assign to each of the three terms equal weights in the above loss function expression.
% It should be also pointed out, that an optimization method could be employed to solve the least squares problem by using (at the very least) first-order derivatives if the objective function in \eqref{eq:loss1} is sufficiently smooth. 

To train the network, we used a greedy zero-th order continuation approach presented in \cite{alvarez2023discrete}, in order to deal with steep gradients resembling singularities in the domain $D$ of interest. In particular, the PINN structure is exploited to solve the system of Eq.\eqref{eq:ConditionFunequations} progressively using a \emph{nested family} of $n$ subdomains:
\begin{equation}
    D_1 \subset D_2 \subset \dots \subset D_k \subset \dots \subset D_n=D.
\end{equation}
The subdomain sizes are chosen to provide a level of sufficient accuracy.

\subsection{Inverse of the Nonlinear Transformation}

As stated in Theorem \ref{thm:thm1}, the nonlinear transformation $T: \mathbb{R}^n \rightarrow \mathbb{R}^n$, is locally invertible and the inverse function, denoted $T^{-1}(z): \mathbb{R}^n \rightarrow \mathbb{R}^n$, has the property that $T^{-1}(T(x)) = x$ for all $x$ in the neighborhood of the equilibrium point $x_0$. %, i.e., $\forall x \in \mathbb{R}^n$ and $T(T^{-1}(\hat{z})) = \hat{z}$ for all $\hat{z}$ in the range of $T$.
In order to reconstruct the trajectories of the system in the phase space, after finding the transformation map $T(x)$ through the PINN approach, a second step is taken aiming at finding the inverse operator $T^{-1}$. In the next paragraph, we show how to compute the inverse operator using  Newton iterations.
Within the proposed framework, simulations are conducted for both the discrete system \eqref{eq:DiscreteSystem} and the corresponding state observer \eqref{TransfObserv} dynamics since the observer dynamic equations are driven by the measurable states (''sensor measurements'') of the discrete system over a specific time horizon. Within  such a context, to determine the state estimate in the original coordinates $x(t)$ through the inverse transformation map, given the simulated current value of the observer state vector in the transformed coordinates $z(t)$ at each time step,  a standard Newton's method is employed by numerically solving w.r.t. $x(t)$ the following equations:
%we consider the transformation z = T(x), where z represents the state observer in the latent space, and T(x) is the solution to a functional equation system, which is approximated by the PINN approach. Our goal is to find the inverse of T, and we formulate this as the following problem.}
 % given a fixed value of $z$, the inverse function \(T^{-1}(z)\), i.e. the output $x$ corresponding to $z$, can be found indirectly using Newton's method by solving w.r.t. $x$ the equation
\begin{equation}
    G(x)=T(x(t)) - z(t) = 0,  
    \label{G(x)}
\end{equation}
where $T(x)$ is the approximation of the transformation map calcualted through the proposed PINN scheme.
The objective is to compute the state vector $x$  which satisfies equations \eqref{G(x)} at each time step of the simulation of the system's \eqref{eq:DiscreteSystem} and observer's \eqref{TransfObserv} discrete-time dynamic equations. The iterative update formula for Newton's method in the multivariate case is given by:
\[
x_{n+1} = x_n - \left(\frac{\partial G}{\partial x}\right)^{-1} G(x_n)
\]
Here, \(x_n\) is the current estimate of the root, \(\frac{\partial G}{\partial x}\) is the Jacobian matrix of \(G\) with respect to \(x\), and \(\left(\frac{\partial G}{\partial x}\right)^{-1}\) is its inverse matrix. %This iterative process continues until a convergence criterion is met, e.g. $||G(x)||<tol$ for a given tolerance value.
Clearly, the solution \(x\) obtained through this process provides a numerical approximation of the inverse function \(T^{-1}(z)\), allowing the reconstruction of the dynamic trajectories/profiles of the system's \eqref{eq:DiscreteSystem} ummeasurable states. %It signifies the value of the input \(x\) corresponding to the output \(z\) under the function \(T\).
For our computations, we have set the relative and absolute tolerance to be 1E$-$06.%, and our condition for the Newton Rhsapson alogorithm is to iterate until both conditions are met i.e. $Error>1E$-$06$ and $iter\leq 100$.
% Moreover, we evaluate the transformation $T(\bm{x})$, which in this case is the output of the PINN in a uniform distributed grid ranging from $[x_0,x_L]\times[x_0,x_L]$ where we select the value $x_L$ according to the region in which the transformation $T$ where found. \textbf{Algorithm 1} describes the procedure for applying Newton's method to compute the inverse of $T(x)$:
% \begin{algorithm}
% \caption{Inverse of PINN using Newton's Method}
% \begin{algorithmic}[1]
%     \State Initialize network parameters and initial guess $X$
%     \State Set iteration counter $iter \gets 0$
%     \State Set maximum iterations $max\_iter \gets 100$
%     \State Set convergence threshold $tolerance \gets 1E$-$06$
    
%     \While{$\text{Error} > \text{Tolerance}$ and $iter < \text{max\_iter}$}
%         \State $iter \gets iter + 1$
        
%         \State Evaluate the function $F(X)$
        
%         \State Compute the Jacobian matrix $J_X = \frac{\partial F}{\partial X}$
        
%         \State Solve for the update $\Delta X$ using $\Delta X = -\text{pinv}(J_X) \cdot F(X)$
        
%         \State Update the guess $X \gets X + \Delta X$
        
%         \State Compute the error $\text{Error} \gets \|\Delta X\|$
%     \EndWhile
    
%     \State \textbf{Output:} Inverse of PINN parameters $X$
% \end{algorithmic}
% \end{algorithm}

%%%%%%%%%%%%%HVA Here two Benchmark problems divided again in two different sections each one
\section{Illustrative Benchmark Problems}
\label{sec:benchmarks}
The performance of the proposed method is evaluated in the following two benchmark case studies \cite{KazantzisKravaris2001}.
\subsection{Benchmark Problem 1}
The following nonlinear discrete-time system is considered with a state space representation:
\begin{gather}
x_1(t+1)=exp\bigg(0.2\frac{x_2(t)}{1+x_2(t)}\bigg)\sqrt{(1+x_1(t)+x_2(t))}-1-0.4x_2(t)-0.5\ln(1+x_1(t)+x_2(t))\nonumber \\
x_2(t+1)=0.5\ln(1+x_1(t)+x_2(t))+0.4x_2(t) \nonumber \\
y(t)=x_2(t),
\label{eq:benchmark1}
\end{gather}
where the measured output variable $y$ coincides with the second state variable $x_{2}$.
Notice that $(x_1, x_2) = (0, 0)$ is an equilibrium point of \eqref{eq:benchmark1} and $1 + x_1 + x_2 > 0,\quad \forall x_1, x_2 \in \mathbb{R}$ (for well-defined right-hand side functions in the above system of difference equations). In addition, the Jacobian matrix $F$ of \eqref{eq:benchmark1} evaluated at the equilibrium point is:
\begin{gather}
    F=\begin{bmatrix}
0.0\quad-0.2 \\ 
0.5\quad0.9
\end{bmatrix}
\end{gather}
The eigenvalues of $F$ are: $\lambda_1 =  0.1298$ and $\lambda_2 = 0.7702$. Consider now the following choice for matrix $A$:
\begin{gather}
    A=\begin{bmatrix}
0.5\quad0.3 \\ 
0.5\quad0.4
\end{bmatrix},
\end{gather}
with eigenvalues: $\mu_1 =0.8405$ and $\mu_2 = 0.0515$. Moreover, the following output
injection map was chosen:
\begin{gather}
b(y)=b(x_2)=\begin{bmatrix}
0.2\big(\frac{x_2}{1+x_2} \big)-0.3x_2 \\ 
0.0
\end{bmatrix}
\end{gather}
such that the pair $(A,B)$ is controllable. Under the above choice for $A$, all other assumptions of Theorem \ref{thm:thm1} are satisfied. Consequently, the associated system of functional equations \eqref{eq:LFEs} takes the following form:
\begin{gather}
T_1\big( u_1,u_2) = 0.5T_1 + 0.3T_2 +0.2\big(\frac{x_2}{1+x_2} \big)-0.3x_2 
% T_1\biggl(exp\bigg(0.2\frac{x_2}{1+x_2}\bigg)\sqrt{(1+x_1+x_2}-1-0.4x_2-0.5\ln(1+x_1+x_2),\nonumber\\ 0.5\ln(1+x_1+x_2+0.4x_2\biggl)=0.5T_1 + 0.3T_2 +0.2\big(\frac{x_2}{1+x_2} \big)-0.3x_2 
\nonumber \\
T_2\big( u_1,u_2) = 0.5T_1 + 0.4T_2
% T_2\biggl(exp\bigg(0.2\frac{x_2}{1+x_2}\bigg)\sqrt{(1+x_1+x_2}-1-0.4x_2-0.5\ln(1+x_1+x_2),\nonumber\\0.5\ln(1+x_1+x_2+0.4x_2\biggl)=0.5T_1 + 0.4T_2 
\nonumber \\
T_1(0,0) = 0 \nonumber\\
T_2(0,0) = 0,
\label{FEsys}
\end{gather}
Where $u_1=exp\big(0.2\frac{x_2}{1+x_2}\big)\sqrt{(1+x_1+x_2}-1-0.4x_2-0.5\ln(1+x_1+x_2)$ and $u_2=0.5\ln(1+x_1+x_2+0.4x_2)$; in addition \eqref{eq:LFEs} admits a unique real analytic and locally (around the equilibrium point of interest) invertible solution. In particular, the solution can be analytically derived and expressed in closed-form as shown below:
\begin{gather}
T_1(x_1,x_2)=\ln(1+x_1+x_2), \quad T_2(x_1,x_2)=x_2. 
\label{Nonlinear_Transformations1}
\end{gather}
Since:
\begin{gather}
\frac{\partial T}{\partial x}(0,0)=\begin{bmatrix}
1\quad1 \\ 
0\quad1
\end{bmatrix},   
\end{gather}
and $det[T]\neq 0$, the above solution $T(x)$ is indeed locally invertible around the equilibrium point
$(x_1, x_2) = (0, 0)$. According to Theorem \ref{thm:thm1}, the proposed discrete-time observer expressed in the transformed state variables $z$ is given by:
\begin{gather}
\hat{z}_1(t+1)=0.5\hat{z}_1(t)+0.3\hat{z}_2(t)+0.2\frac{y(t)}{1+y(t)}-0.3y(t) \nonumber\\
\hat{z}_2(t+1)=0.5\hat{z}_1(t)+0.4\hat{z}_2(t) 
\end{gather}
Therefore, through the inverse transformation map, the state estimates expressed in the original state variables are given by the following equations:
\begin{gather}
\Bar{x}_1(t)=exp(\hat{z}_1(t))-\hat{z}_2(t)-1 \nonumber \\
\Bar{x}_2(t)=\hat{z}_2(t)
\end{gather}
%%%%%%%%%%%%%%%%%%%%%% 
Please notice that the nonlinear transformation map in the first benchmark example exhibits singularities when $x_1+x_2=-1$ in the system's state space. The primary objective is to learn the transformation map in the domain $[x_L,0] \times [x_L,0] =[-0.495,0]\times[-0.495,0]$.\par
\begin{figure}[htbp]
 \centering
 \begin{subfigure}[h]{0.45\textwidth}
     \centering
\includegraphics[width=\textwidth]{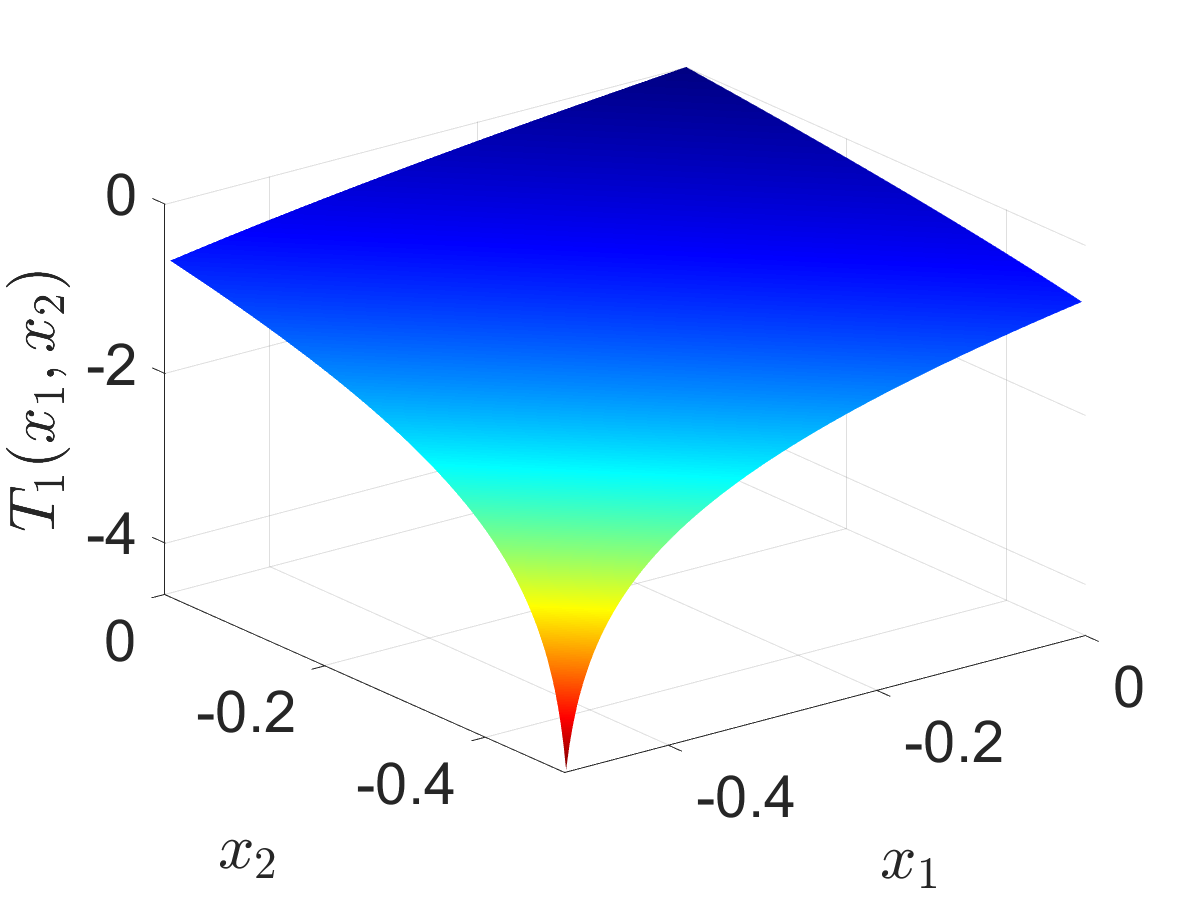}
     \caption{}
 \end{subfigure}
 \hfill
 \begin{subfigure}[h]{0.45\textwidth}
     \centering
     \includegraphics[width=\textwidth]{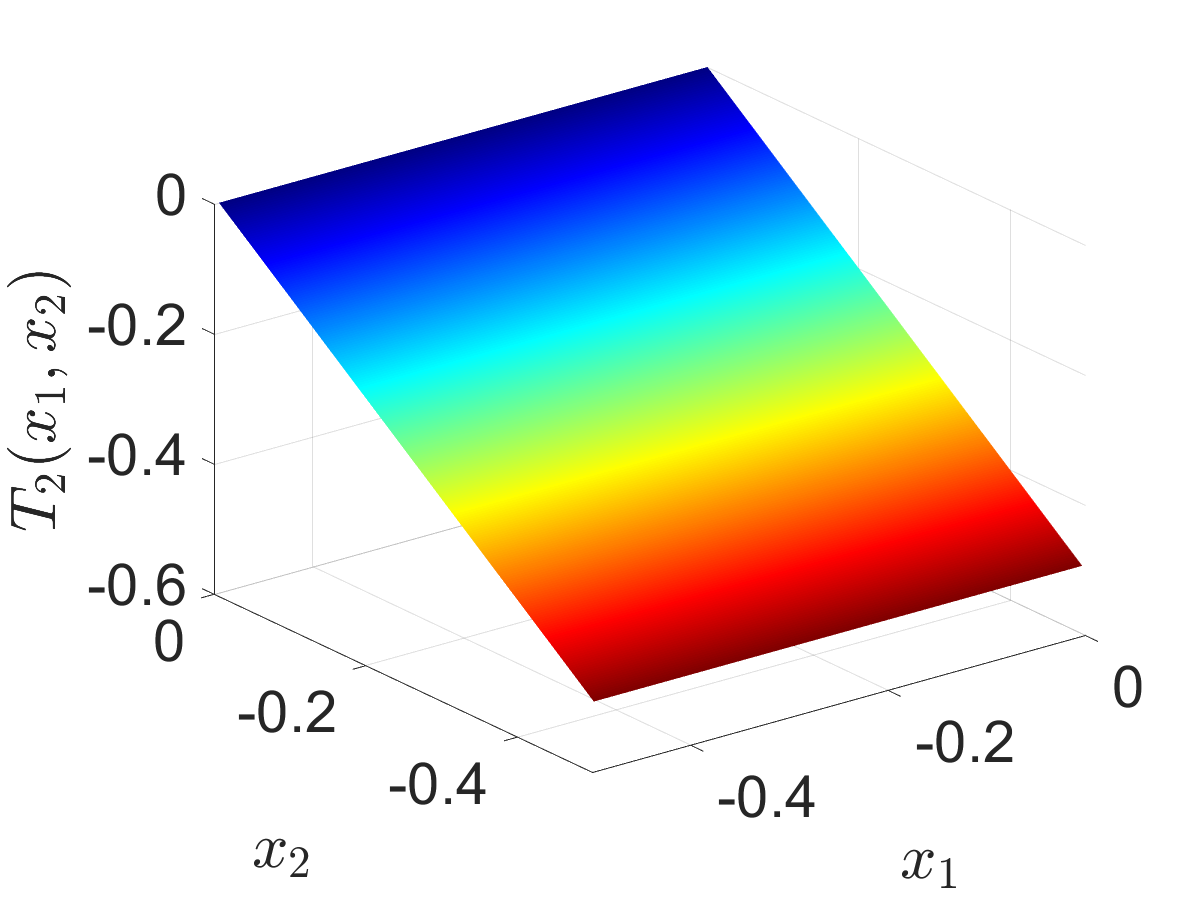}
     \caption{}
 \end{subfigure}
    \caption{Analytical solution of the NFEs (\ref{FEsys}) in $[-0.495,0]\times [-0.495,0] $. (a) $T_1(x_1,x_2)=\ln(1+x_1+x_2)$. A steep-gradient at $(-0.495,-0.495)$ is due to the presence of a singular point at $(x_1,x_2)=(-0.5,-0.5)$. (b) $T_2(x_1,x_2)=x_2$.}.
  \label{fig:sX}
\end{figure}

 
Figure (\ref{fig:sX}) depicts the analytical solutions $T_1(x_1,x_2), T_2(x_1,x_2)$, of the associated system of functional equations \eqref{FEsys}. 
We note that it encompasses the solution domain which has been deliberately chosen to be $x_1,x_2\in [-0.495,0]$ since $T_1(x_1,x_2)$ exhibits a singular point at $(x_1,x_2)=(-0.5,-0.5)$. Therefore, the first reasonable step would be to comparatively evaluate the numerical approximation accuracy of the proposed PINN scheme in this region against the existing analytical solution and assess its impact on the performance profile of the resulting discrete-time observer.
\subsection{Benchmark Problem 2}
Consider also the following nonlinear discrete-time dynamical system \cite{Xiao2003}:
\begin{gather}
 x_1(t+1)=\frac{0.5\frac{x_1(t)}{1+x_1(t)}-0.9x_2(t)}{1-0.5\frac{x_1(t)}{1+x_1(t)}+0.9x_2(t)} \nonumber \\
 x_2(t+1)=x_2(t) \nonumber \\
 y(t)=x_1(t)
 \label{eq:Benchmark2}
\end{gather}
where $x_2$ represents an unidentified constant parameter or disturbance term which needs to be estimated and $x_1$ the measured state variable. Please notice that $(x_1, x_2) = (0, 0)$ is an equilibrium point of
the above system and the Jacobian matrix $F$ evaluated at this point is:
\begin{gather}
    F=\begin{bmatrix}
0.5\quad-0.9\\ 
0.0\quad1
\end{bmatrix}.
\end{gather}
The eigenvalues of $F$ are: $\lambda_1 =  0.5$ and $\lambda_2 = 1$, and therefore, the original nonlinear discrete-time observer design method presented in \cite{KazantzisKravaris2001} can not be applied since the second eigenvalue is not located within the Poincar\'e domain; instead it lies on the unit circle. However, under the following choice for the matrix $A$
\begin{gather}
    A=\begin{bmatrix}
0.0\quad0.0 \\ 
0.0\quad0.1
\end{bmatrix}  
\end{gather}
with eigenvalues: $\mu_1 = 0$ and $\mu_2 = 0.1$, and an output injection map:
\begin{gather}
b(y)=b(x_{1})=\begin{bmatrix}
0.5\big(\frac{x_1}{1+x_1} \big) \\ 
\big(\frac{x_1}{1+x_1} \big)
\end{bmatrix}
\end{gather}
all assumptions in the method's extension presented in \cite{Xiao2003} (where the Poincar\'e domain assumption was relaxed and replaced by a more general one of the Siegel-type) are met, and therefore, the proposed nonlinear discrete-time observer exists and can be similarly designed. Indeed, the associated system of functional equations takes the following form:
\begin{gather}
    T_1\biggl(\frac{0.5\frac{x_1}{1+x_1}-0.9x_2}{1-0.5\frac{x_1}{1+x_1}+0.9x_2} ,x_2\biggl)=0.5\frac{x_1}{1+x_1} \nonumber \\
     T_2\biggl(\frac{0.5\frac{x_1}{1+x_1}-0.9x_2}{1-0.5\frac{x_1}{1+x_1}+0.9x_2} ,x_2\biggl)=0.1T_2+\frac{x_1}{1+x_1} \nonumber \\
       T_1(0,0) = 0 \nonumber\\
T_2(0,0) = 0,
\label{FEs2}
\end{gather}
The above system of functional equations (\ref{FEs2}) admits a unique analytically derived closed-form solution shown below:
\begin{gather}
    T_1(x_1,x_2)=\frac{x_1}{1+x_1}+0.9x_2 \nonumber\\
    T_2(x_1,x_2)=\frac{5}{2}\biggl(\frac{x_1}{1+x_1} + x_2\biggl).
    \label{Nonlinear_Transformations2}
\end{gather}
Since:
\begin{gather}
\frac{\partial T}{\partial x}(0,0)=\begin{bmatrix}
1\quad0.9 \\ 
2.5\quad2.5
\end{bmatrix},   
\end{gather}
and  $det[T]\neq 0$,  the above solution $T(x)$ is locally invertible around the origin. The proposed discrete-time observer is then given by the following dynamic equations:
\begin{gather}
\hat{z}_1(t+1)=0.5\frac{y(t)}{1+y(t)} \nonumber \\
\hat{z}_2(t+1)=0.1\hat{z}_2(t)+ \frac{y(t)}{1+y(t)}.
\label{observerEx2}
\end{gather}
Therefore, using the inverse transformation map, the state estimates expressed in the original state variables are given by the following equations:
\begin{gather}
\Bar{x}_1(t)=\frac{10\hat{z}_1(t)-3.6\hat{z}_2(t)}{1-10\hat{z}_1(t)+3.6\hat{z}_2(t)} \nonumber \\
\Bar{x}_2(t)=4\hat{z}_2(t)-10\hat{z}_1(t).
\label{AinverseEx2}
\end{gather}
In this second benchmark problem, the nonlinear transformation map exhibits a singularity when $x_1=-1$. Here, we aim to learn this map in the domain $[x_L,0] \times [x_L,0] =[-0.91,0]\times[-0.91,0]$.\par
\begin{figure}[htbp]
 \centering
 \begin{subfigure}[h]{0.45\textwidth}
     \centering
     \includegraphics[width=\textwidth]{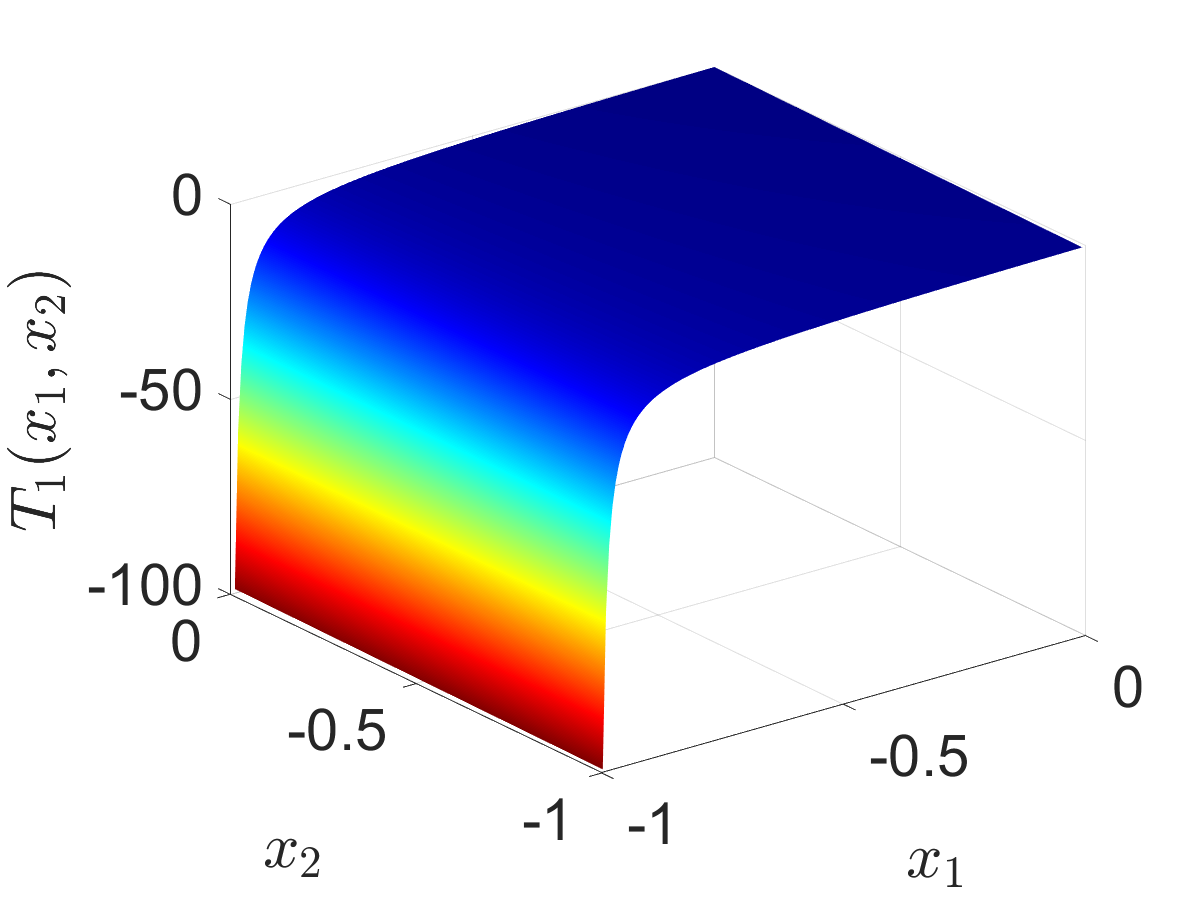}
     \caption{}
 \end{subfigure}
 \hfill
 \begin{subfigure}[h]{0.45\textwidth}
     \centering
     \includegraphics[width=\textwidth]{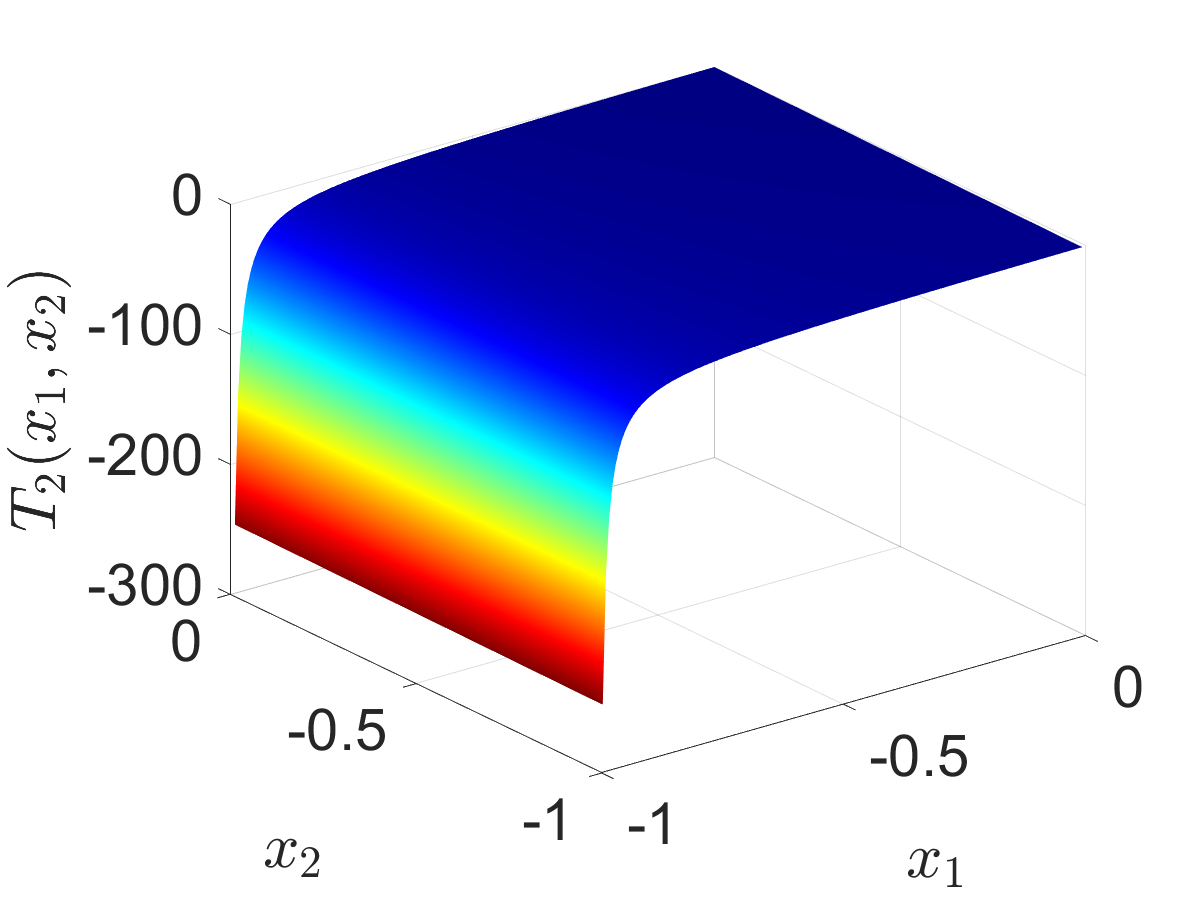}
     \caption{}
 \end{subfigure}
    \caption{Analytical solution of the functional equations (\ref{FEsys}). (a) $ T_1(x_1,x_2)=\frac{x_1(t)}{1+x_1(t)}+0.9x_2$. (b) $ T_2(x_1,x_2)=\frac{5}{2}\biggl(\frac{x_1(t)}{1+x_1(t)} + x_2\biggl)$. A steep-gradient at $(-0.91,-0.91)$ is due to the presence of a singular point at $x_1=-1$.}
  \label{fig:Analytical_transformation_ex2}
\end{figure}

 
Figure (\ref{fig:Analytical_transformation_ex2}) depicts the analytical solutions $T_1(x_1,x_2), T_2(x_1,x_2)$, of the associated system of functional equations \eqref{FEs2} . 
Please notice that it includes the solution domain, which in this case is $(x_1,x_2),\ \in$ $[-0.91,0]$ because $T_1(x_1,x_2)$ displays a singular point at $(x_1,x_2)=(-1,-1)$. As in the previous benchmark example, the focus is placed on the evaluation of the impact of the proposed PINN scheme on the performance profile of the resulting nonlinear observer by comparing the identified transformation map to the analytical solution and hence, demonstrating the enhanced numerical approximation accuracy in this domain attained by PINN.
\newpage

\section{Numerical results}\label{sec:numerical_results}
We utilized the Matlab deep-learning toolbox, where we wrote a custom code to implement the Physics-Informed Neural Network (PINN) strategy. In this process, we employed the Levenberg-Marquardt (LM) algorithm, to optimize the learnable, parameters. The scheme utilizes finite differences for computing the required derivatives. However, it's worth noting that these derivatives can also be computed analytically, numerically or through Automatic Differentiation (AD), as demonstrated in \cite{alvarez2023discrete}. The simulation results obtained when the PINN was trained in the entire domain without the greedy approach are also presented in order to evaluate the efficacy of the proposed greedy strategy.
For training, we used $15 \times 15$ equispaced distributed points in the interval $[a, b]$:$,\ a=-0.495$,$,\ b=0$, for the model (\ref{eq:benchmark1}), and $a=-0.91$, $b=0$, for the model (\ref{eq:Benchmark2}).
For the test set, we employed Chebyshev-Lobatto nodes, specifically selecting the roots of the (second kind) Chebyshev-Lobatto polynomial of degree $k$, i.e.  %In the interval $[a, b]$ 
%with $a,b\in\mathbb{R}$, the grid was established using these Chebyshev collocation points:
\begin{gather}
x_n=\frac{1}{2}(a+b)-\frac{1}{2}(a-b)cos\left ( \frac{2n-1}{2n}\pi \right )\quad n=1,\cdots,k,
\label{Cheby}
\end{gather}

\paragraph{Power-series solution.} 
In accordance with the theoretical framework outlined in reference \cite{KazantzisKravaris2001}, we also used the power series solution method delineated in Section \ref{sec:framework}, to find the nonlinear transformation. In particular, we used sixth-order multivariate Taylor polynomial approximations. The method produces a system of linear algebraic equations involving the unknown Taylor coefficients of the nonlinear transformation $T(x_1, x_2)$. The solution to this system can be computed using MATLAB's symbolic toolbox, as explained in Section \ref{sec:framework}. We opt for a sixth-order multivariate Taylor polynomial expansion for two reasons: firstly, the number of coefficients aligns well with the number of unknowns in our PINN, facilitating a fair comparison between the two schemes. Secondly, solving the system for an expansion of this order is already computationally challenging, highlighting the underlying issues of this traditional methodological approach usually employed to solve this type of problems for high order series expansions.
 
\paragraph{PINN-based solution.} 
We have used feedfordward network (FNNs) with two hidden layers, with each layer containing five neurons. We have used sigmoid activation functions $\sigma(x)$.
In the training process, we computed the derivatives for the optimization process using finite differences for the LM algorithm, we have set the maximum function evaluations to $500,000$
while the maximum iteration parameter was set to $50,000$. the finite difference step was $1e^{-05}$ and the step size tolerance was set to $1e^{-12}$.
%It is worth mentioning that despite experimenting with various configurations, including the exploration of additional neurons and hidden layers, it was observed that the convergence of the scheme did not exhibit any improvement. 
%\textcolor{blue}{Costas: Maybe I missed it somewhere, but you dont say how many neurons you use for each layer}.
This realization prompted the exploration of alternative strategies to enhance the training process. Subsequently, the adoption of a greedy-wise procedure emerged as a pivotal step in refining the PINN's convergence behavior. By incorporating this approach, the greedy-wise procedure effectively addressed challenges related to convergence of the PINN close to the singular points.

\paragraph{Greedy-wise procedure.}
We initiated training within the region $[-0.1, 0] \times [-0.1, 0]$. Subsequently, we systematically expanded the grid size in the domain, with the step size in the grid varying depending on the specific benchmark problem at hand. In the first benchmark problem, we employed a grid step size of $-0.1$ in both directions, until reaching a domain size of $[-0.4, 0] \times [-0.4, 0]$, then the step size is decreased to a value of $-0.01$ until reaching $[-0.49, 0] \times [-0.49, 0]$. Finally, the step size is decreased again to a value of $-0.001$ until reaching the domain of interest, i.e. $[-0.495, 0] \times [-0.495, 0]$, while for the second benchmark problem,  we initiated training within the region $[-0.1, 0] \times [-0.1, 0]$ and we adopted a step of $-0.1$ until reaching a grid of  $[-0.8, 0] \times [-0.8, 0]$. After that we selected a step size of $-0.01$ until reaching the domain of interest $[-0.91, 0] \times [-0.91, 0]$.
In both instances, following the completion of training within each subdomain, we utilized the learnable parameters specifically, the weights and biases obtained from the preceding subdomain as the initial guess for the optimization algorithm in the subsequent subdomain.
%%%%%%%%----------------------------------
\subsection{Benchmark Problem 1}
\begin{figure}[htbp]
 \centering
 \begin{subfigure}[h]{0.45\textwidth}
     \centering
     \includegraphics[width=7.5cm]{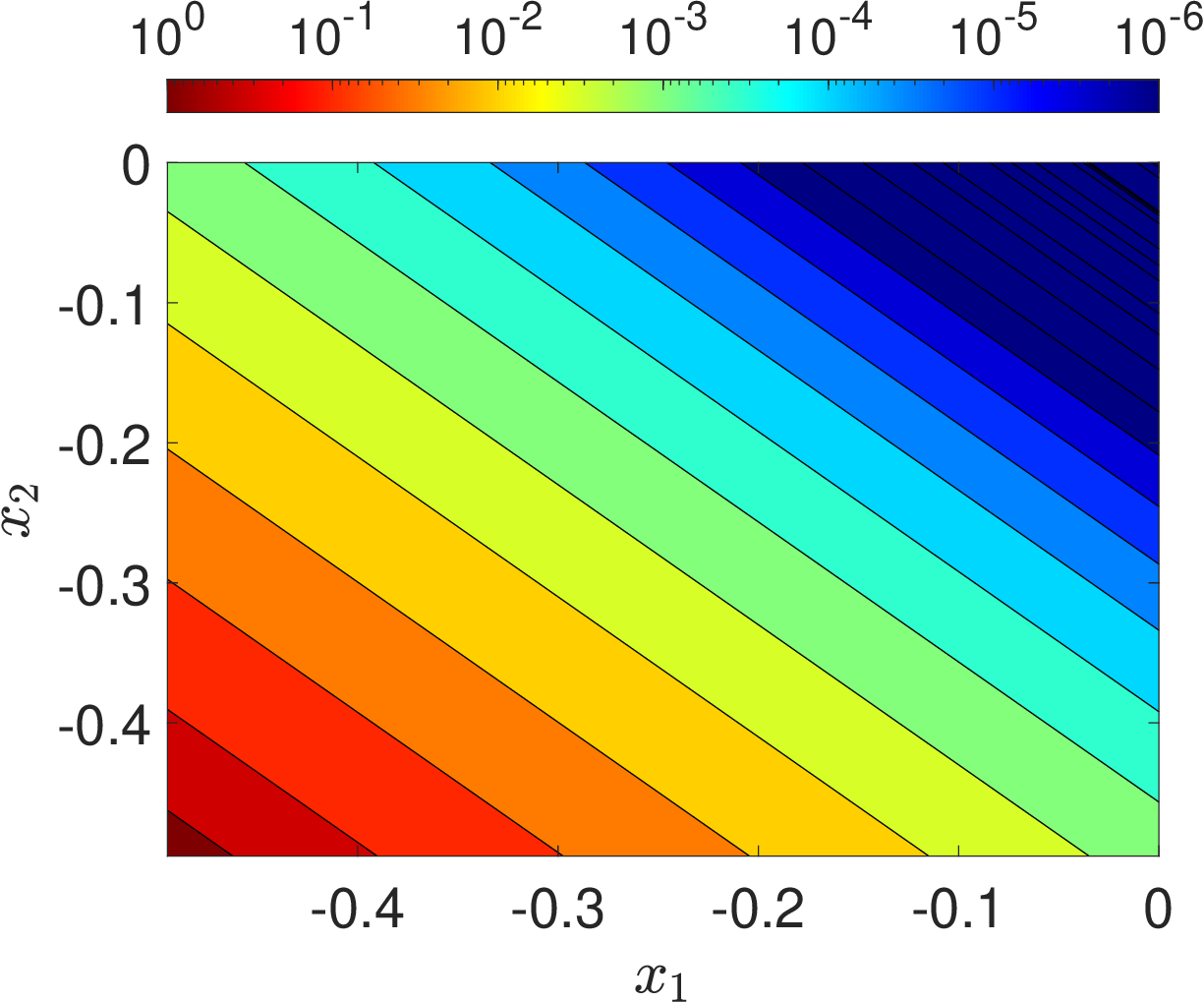}
     \caption{}
 \end{subfigure}
 \hfill
 \begin{subfigure}[h]{0.5\textwidth}
     \centering
     \includegraphics[width=7.5cm]{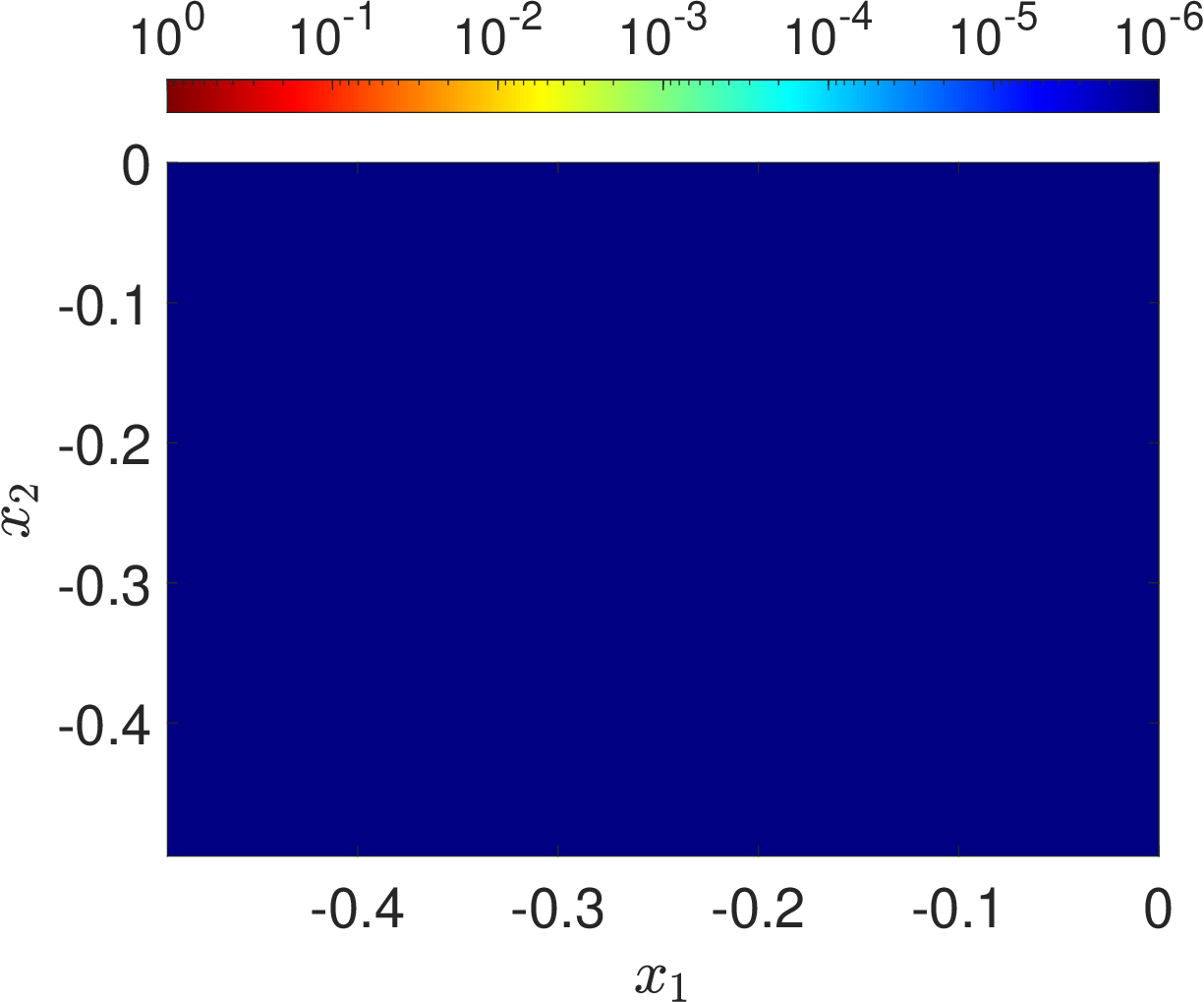}
     \caption{}
 \end{subfigure}
 \begin{subfigure}[h]{0.45\textwidth}
     \centering
     \includegraphics[width=7.2cm]{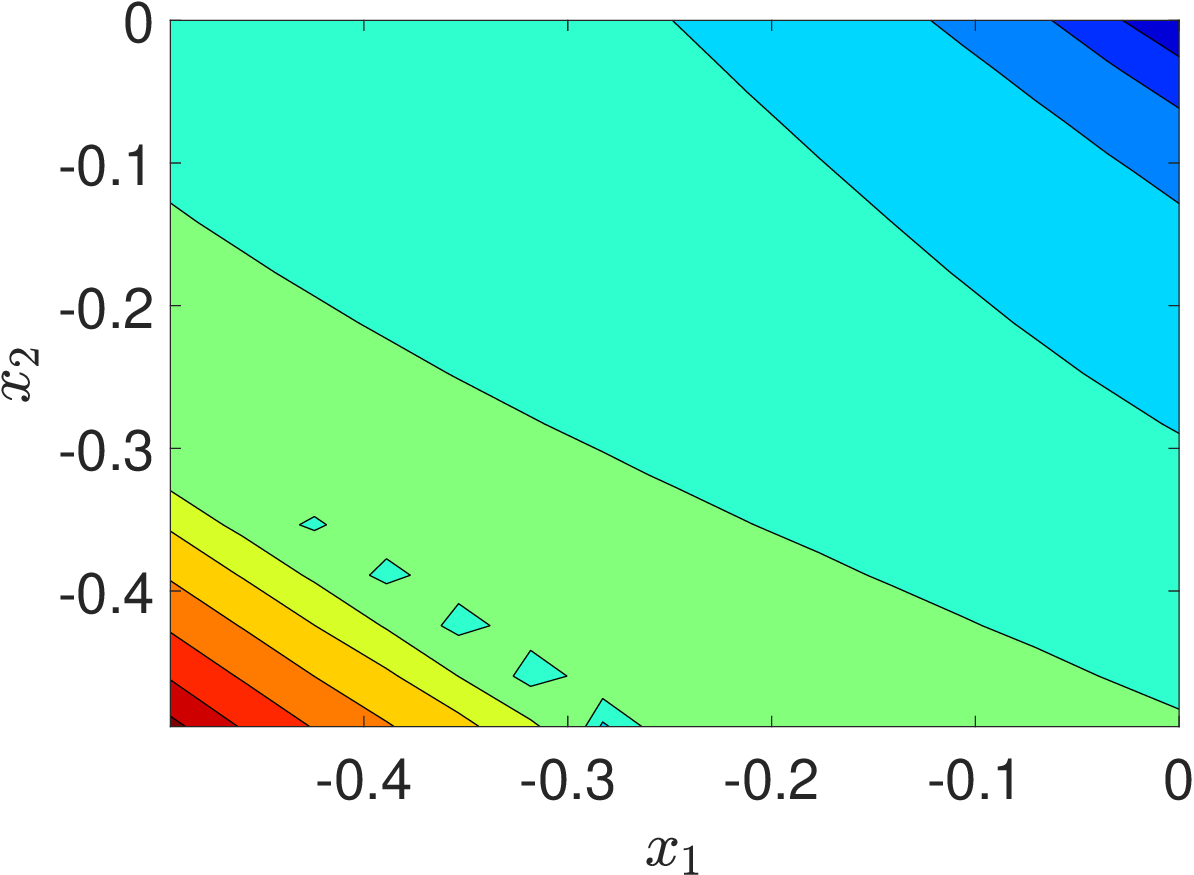}
     \caption{}
 \end{subfigure}
 \hfill
 \begin{subfigure}[h]{0.5\textwidth}
     \centering
     \includegraphics[width=7.2cm]{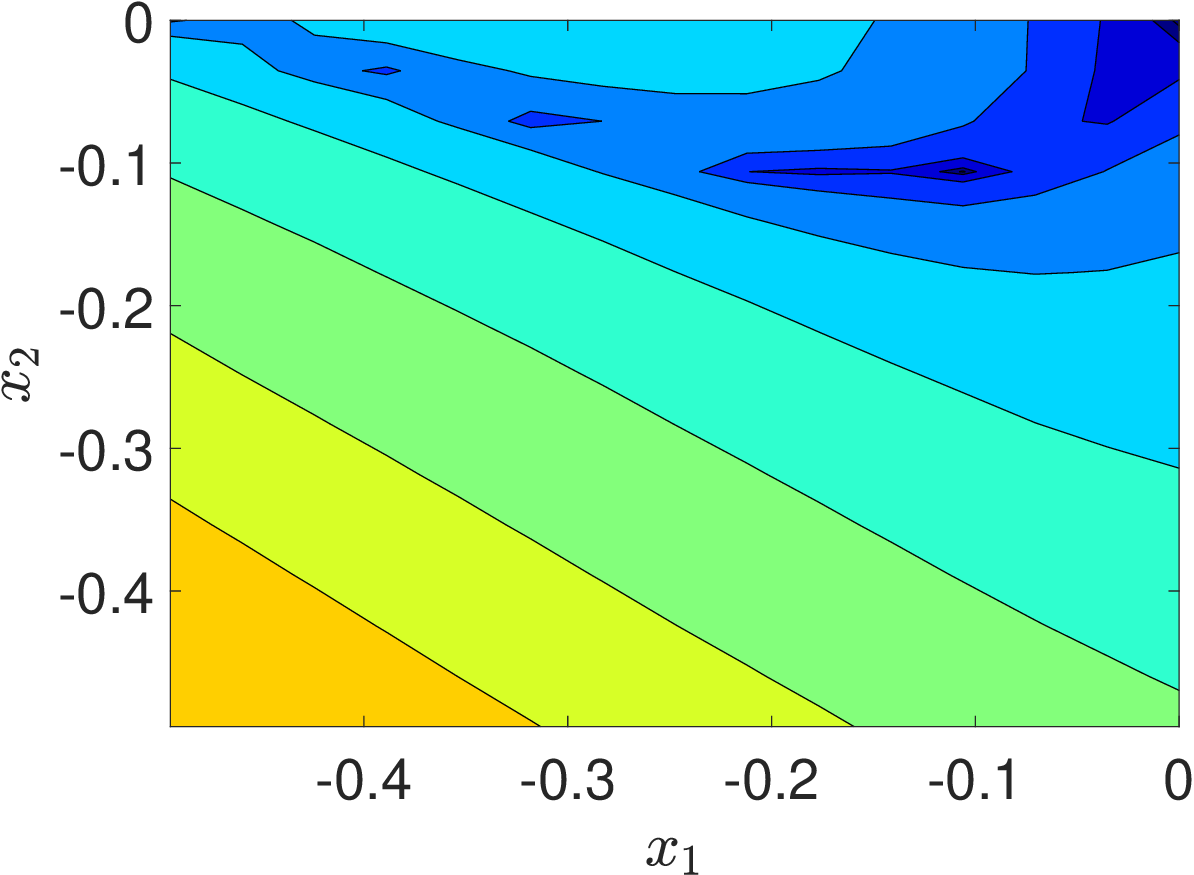}
     \caption{}
 \end{subfigure}
   \centering
 \begin{subfigure}[h]{0.45\textwidth}
     \centering
     \includegraphics[width=7.2cm]{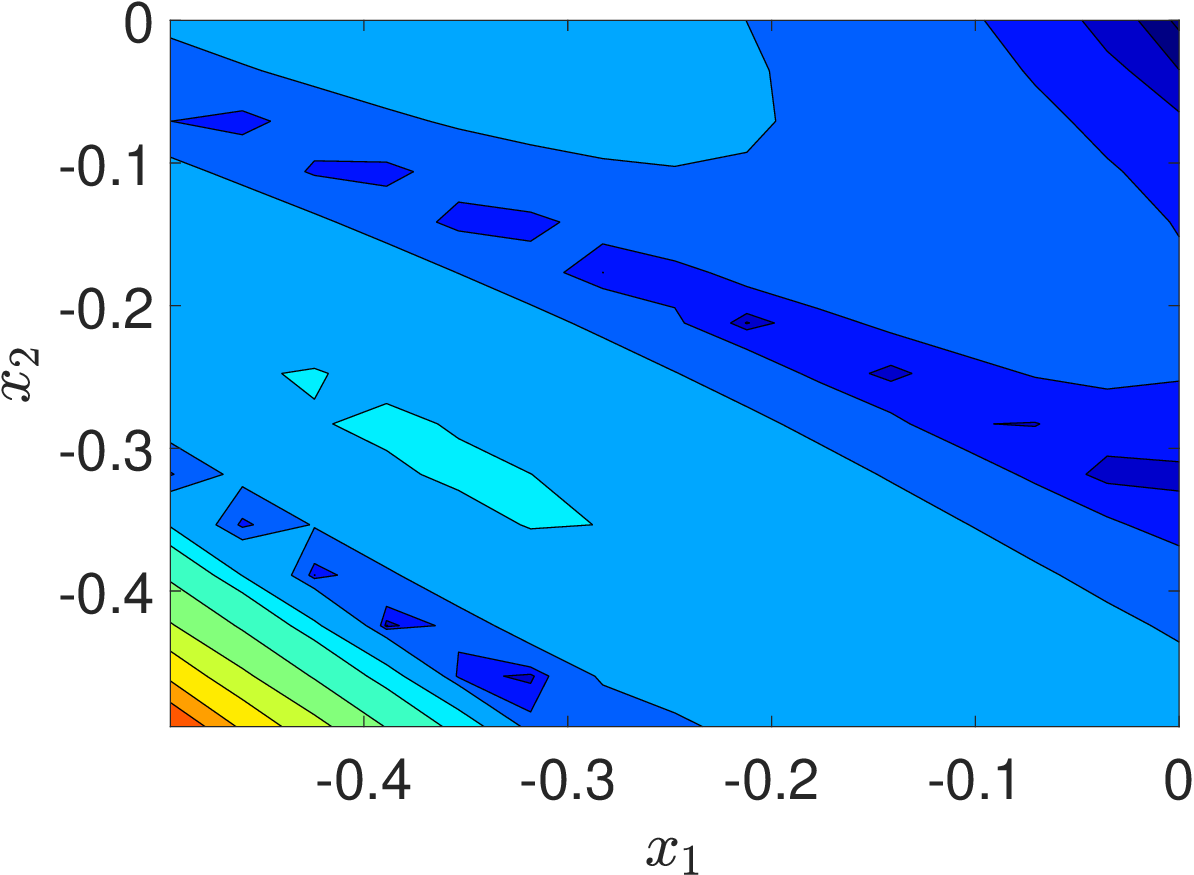}
     \caption{}
 \end{subfigure}
 \hfill
 \begin{subfigure}[h]{0.5\textwidth}
     \centering
     \includegraphics[width=7.2cm]{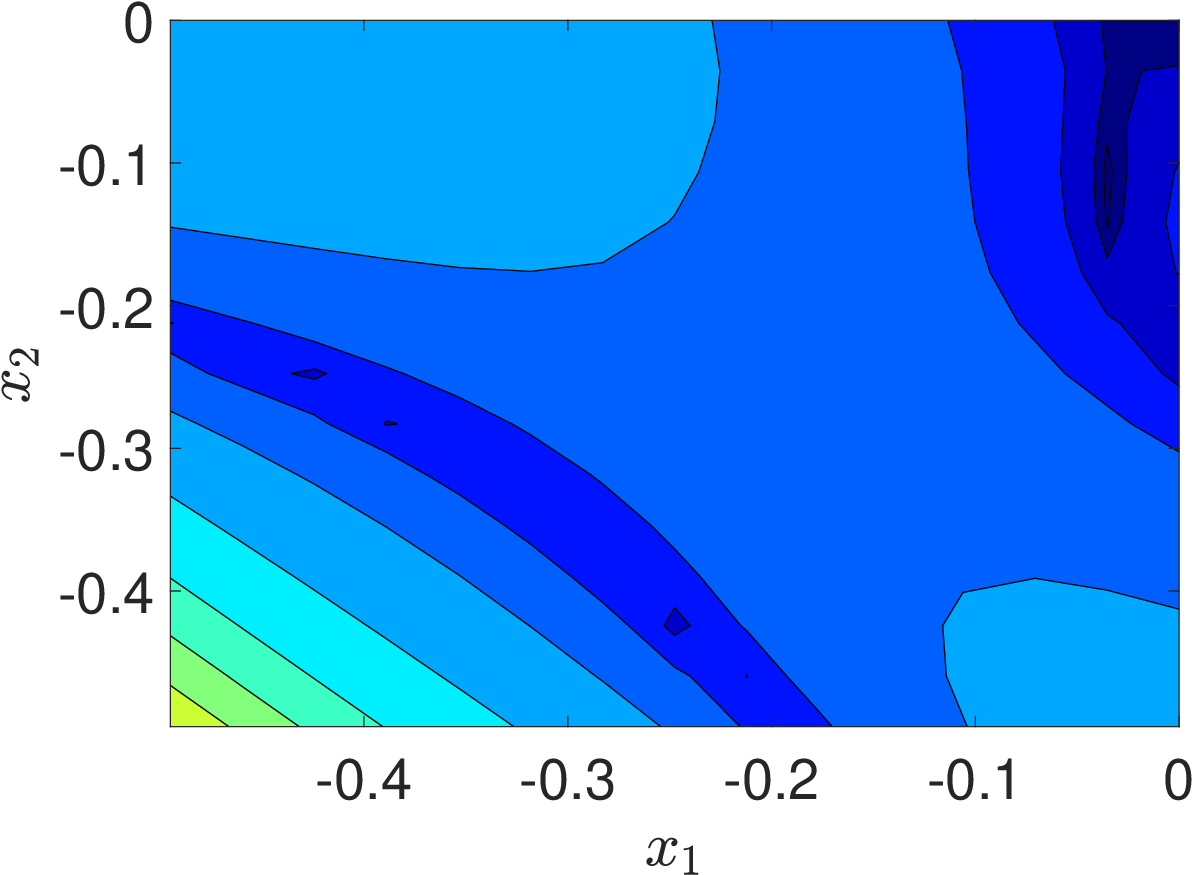}
     \caption{}
 \end{subfigure}
\caption{Benchmark problem 1. Training sets (grids of $15 \times 15$ equispaced distributed points). Numerical approximation accuracy (difference between the computed and analytical solution) of $T_1(x_1,x_2)$ (left column) and $T_2(x_1,x_2)$ (right column) using the various schemes. (a),(b) $6th$ order power-series expansion of $T_1(x_1,x_2)$ and $T_2(x_1,x_2)$ and the right-hand side of the model (\ref{eq:benchmark1}) in $[-0.495,0]\times[-0.495,0]$. (c),(d) PINN trained in the entire domain $[-0.495,0]\times[-0.495,0]$. (e),(f) PINN trained via the greedy-wise approach. }
\label{fig:Modelavailable_TrS}
\end{figure}

%%%%%%%%%%---------------------------------
Figure (\ref{fig:Modelavailable_TrS}) depicts the numerical approximation accuracy, measured as the difference between computed and analytical solutions for the training set, which is built using a grid of $15\times15$ nodes uniformly distributed. The results using the 6th-order power series expansion of both the nonlinear transformation and the right-hand side of the discrete model are displayed in subplots 3(a) and 3(b), respectively. As expected, the power series expansion yields zero error for the $T_2(x_1,x_2)$ component and demonstrates good approximation accuracy for the $T_1(x_1,x_2)$ component, but primarily in the vicinity of the linearization-relevant point at $[0,0]$. Beyond this point, the numerical approximation deteriorates significantly, reaching the order of $10^0$ at the grid edge, particularly close to the singular point. This situation could potentially undermine the precision with which the observer estimates the unmeasured state.\par
\begin{figure}[htbp]
 \centering
 \begin{subfigure}[h]{0.45\textwidth}
     \centering
     \includegraphics[width=7.5cm]{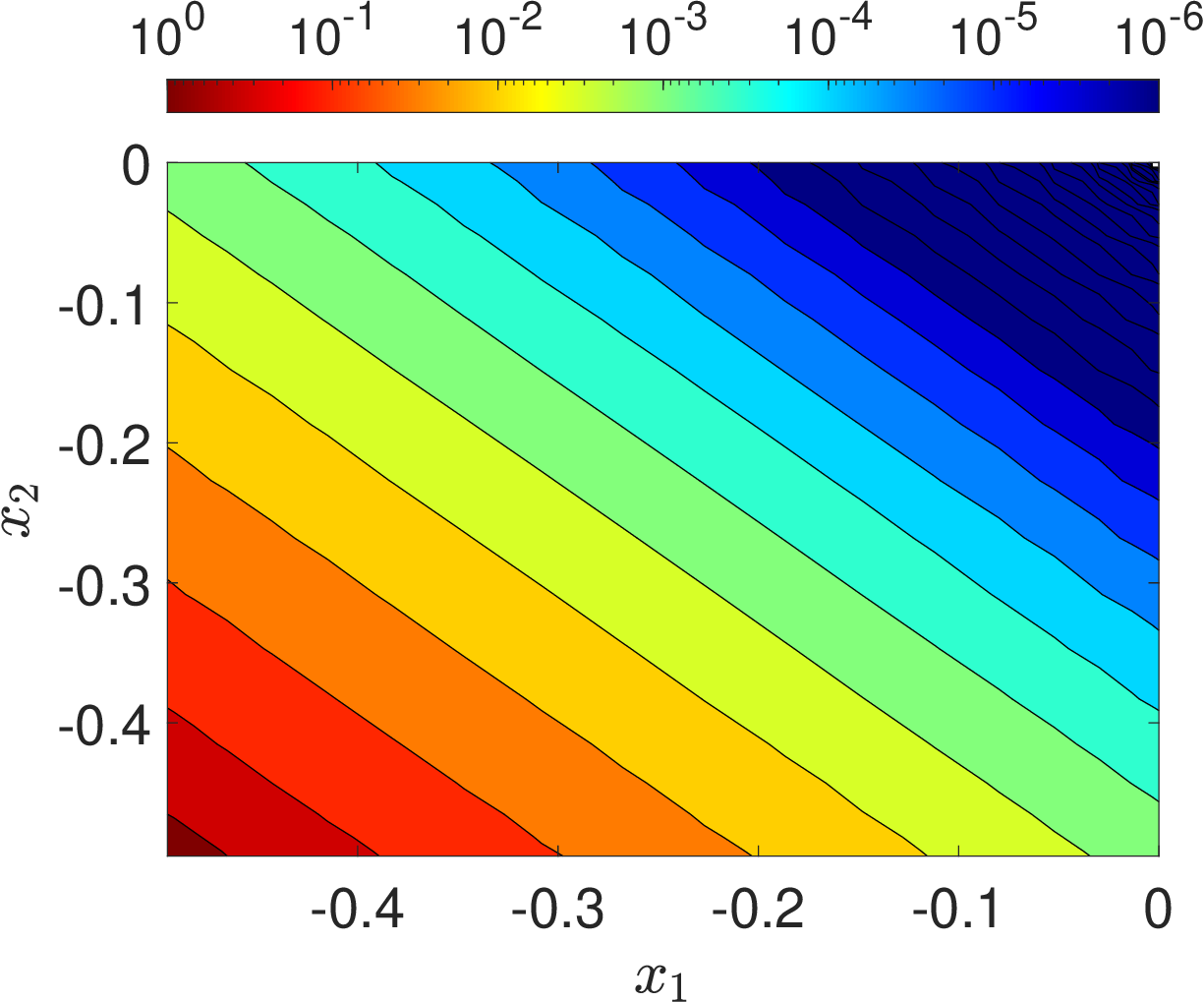}
     \caption{}
 \end{subfigure}
 \hfill
 \begin{subfigure}[h]{0.5\textwidth}
     \centering
     \includegraphics[width=7.5cm]{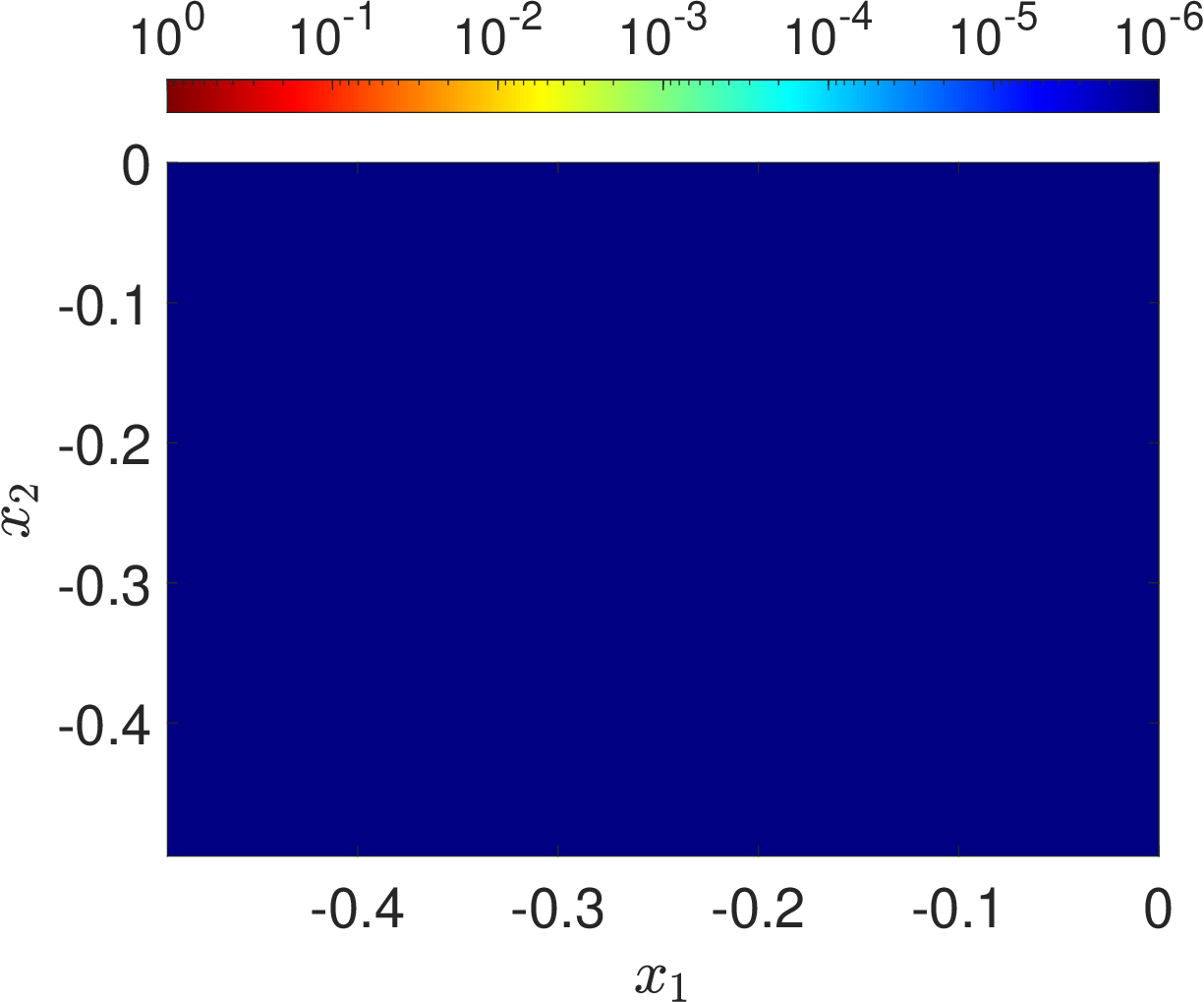}
     \caption{}
 \end{subfigure}
 \begin{subfigure}[h]{0.45\textwidth}
     \centering
     \includegraphics[width=7.2cm]{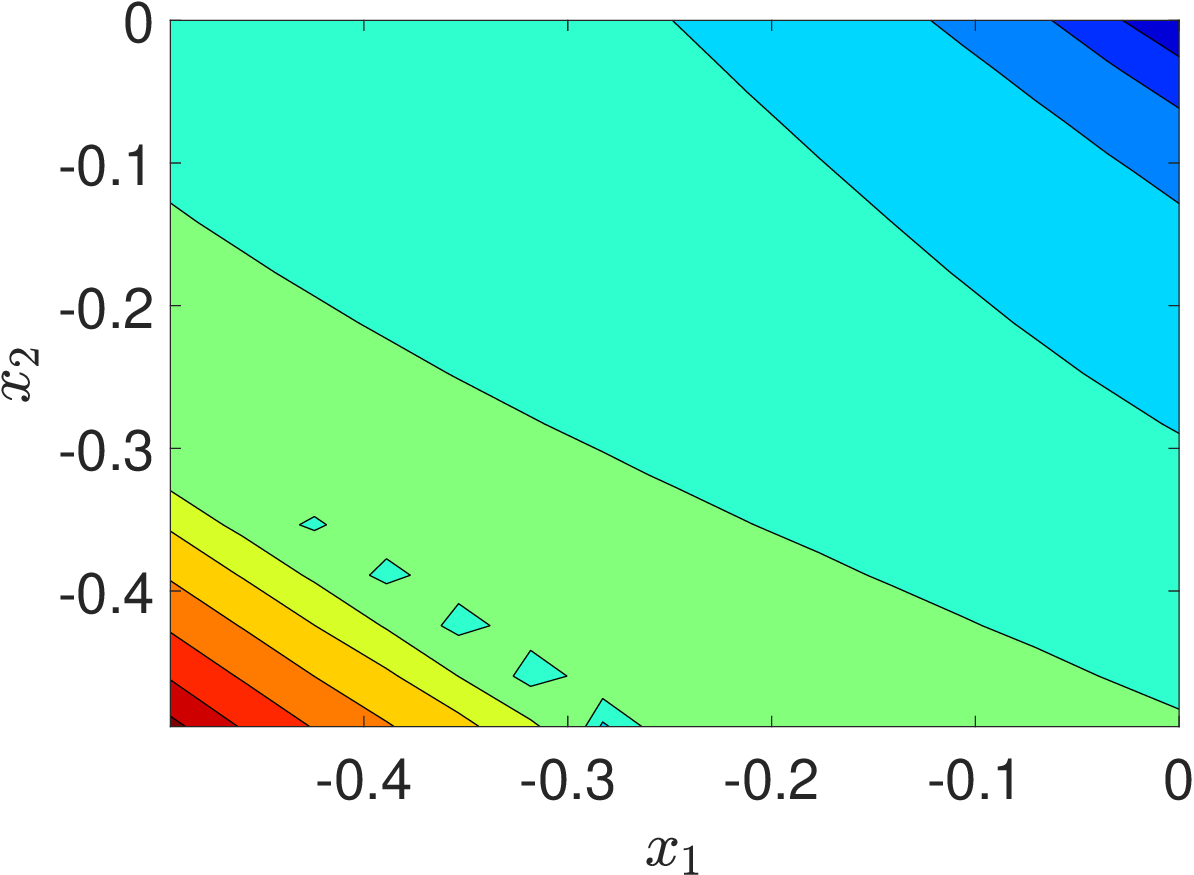}
     \caption{}
 \end{subfigure}
 \hfill
 \begin{subfigure}[h]{0.5\textwidth}
     \centering
     \includegraphics[width=7.2cm]{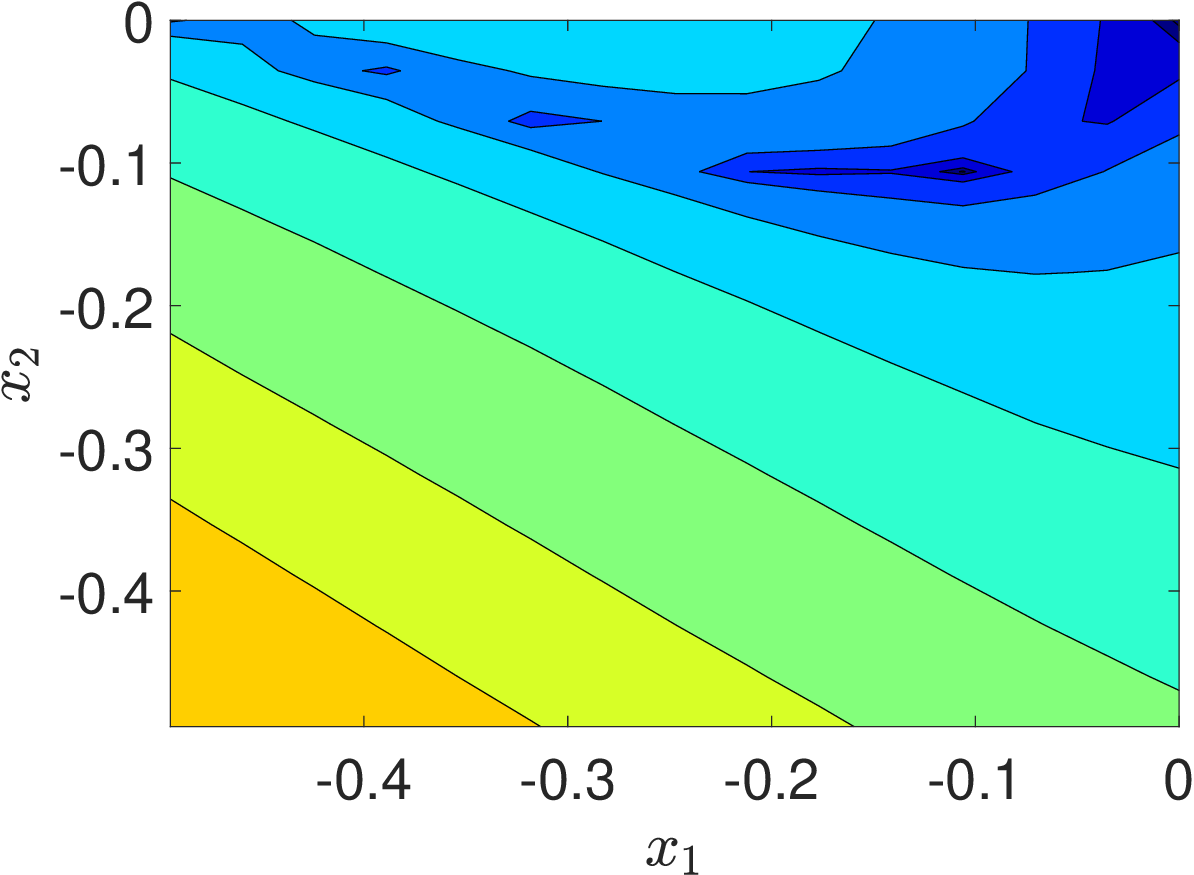}
     \caption{}
 \end{subfigure}
 \centering
 \begin{subfigure}[h]{0.45\textwidth}
     \centering
     \includegraphics[width=7.2cm]{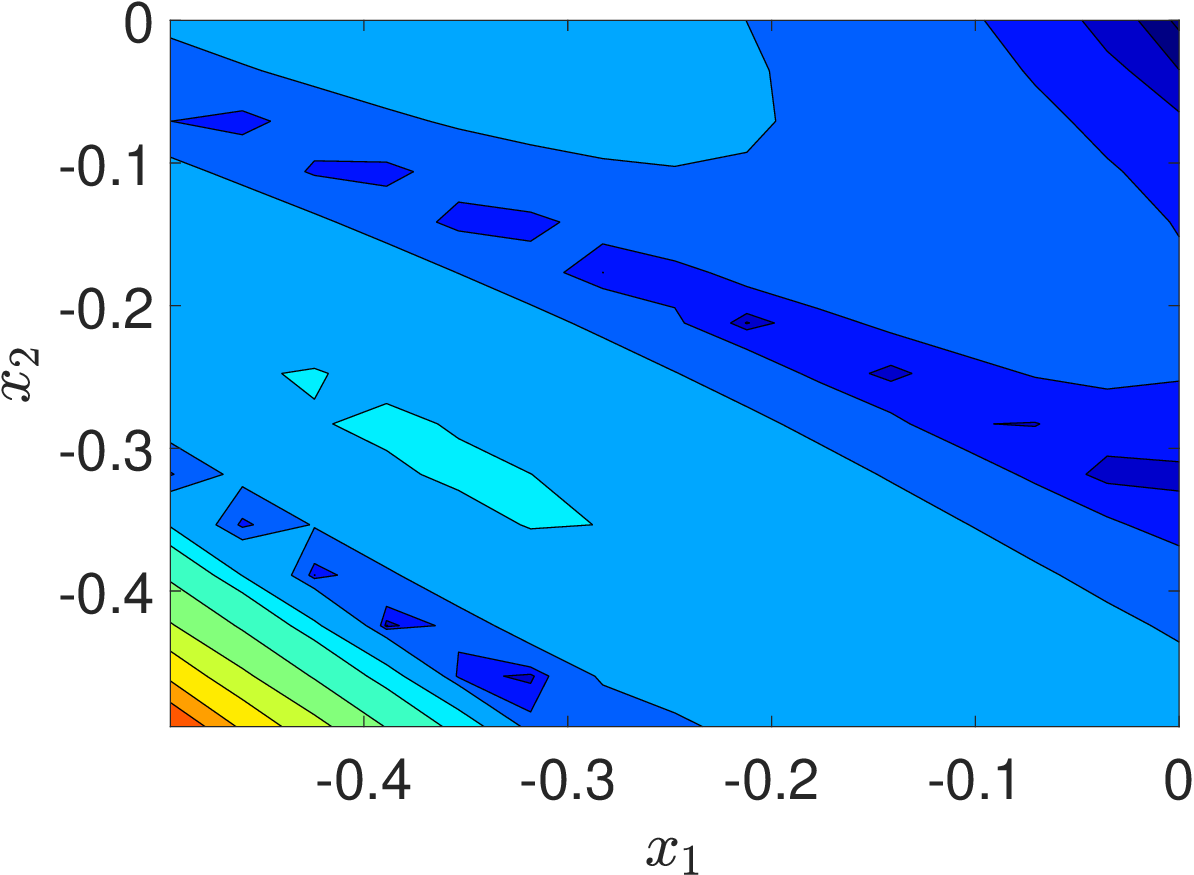}
     \caption{}
 \end{subfigure}
 \hfill
 \begin{subfigure}[h]{0.5\textwidth}
     \centering
     \includegraphics[width=7.2cm]{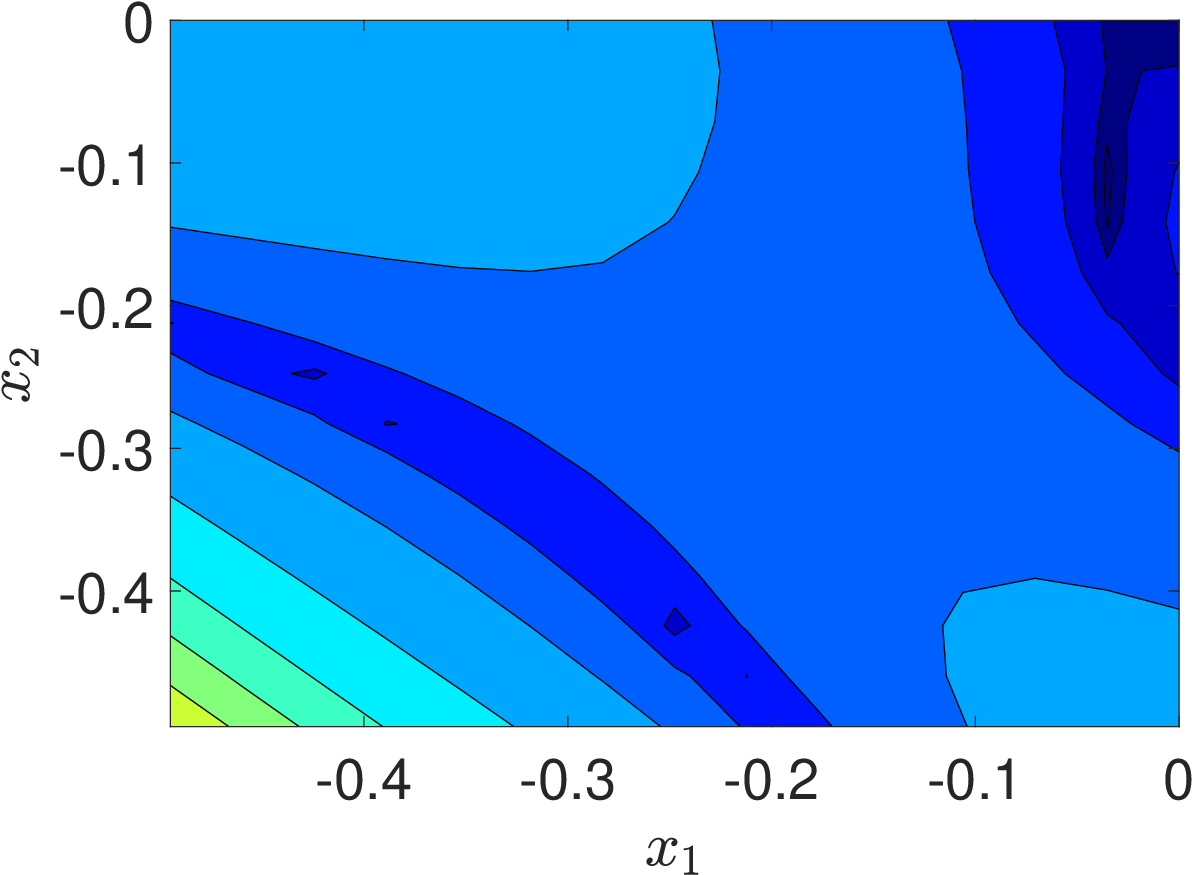}
     \caption{}
 \end{subfigure}
\caption{Benchmark problem 1. Test sets (grids of $20 \times 20$ Chebyshev-distributed points). Numerical approximation accuracy (difference between the computed and analytical solution) of $T_1(x_1,x_2)$ (left column) and $T_2(x_1,x_2)$ (right column) using the various schemes. (a),(b) $6th$ order power-series expansion of $T_1(x_1,x_2)$ and $T_2(x_1,x_2)$ and the right-hand side of the model (\ref{eq:benchmark1}) in $[-0.495,0]\times[-0.495,0]$. (c),(d) PINN trained in the entire domain $[-0.495,0]\times[-0.495,0]$. (e),(f) PINN trained via the greedy-wise approach.}
\label{fig:Modelavailable_TS}
\end{figure}

Figures (\ref{fig:Modelavailable_TrS})(c),(d) depict the approximation accuracy of the PINN, when trained in the entire domain (without the greedy approach). As seen, the approximation accuracy of the PINN remains rather poor, of  the order of $1$ in the vicinity of the singularity. Figures (\ref{fig:Modelavailable_TrS})(e),(f) illustrate the approximation error using PINN with the greedy-wise training approach. As shown, the performance of the PINN is significantly better compared to the one involving training across the entire domain. Furthermore, Figure (\ref{fig:Modelavailable_TS}) showcases the performance of the schemes for the test set, which consist of a grid of  $20\times20$ Chebyshev-Lobatto distributed nodes. Please notice that these results are similar to those observed in the training set.\par
Table (\ref{tab:1}) presents a comparison of $L_1$, $L_2$, $L_{\infty}$ error norms between the analytical and numerical solutions on the test set. Specifically, it summarizes the results obtained with the 6th order power series expansions of $T_1(x_1, x_2)$ and $T_2(x_1, x_2)$, and also with PINN solution approximations both trained in the entire domain at once using the greedy approach. The reported errors include the median, 5th percentile, and 95th percentile over 100 runs.
\\

Table (\ref{tab:2}), on the other hand, reports the results of uncertainty quantification carried out to investigate the convergence of the PINN scheme in approximating the nonlinear transformation maps $T_1(x_1, x_2)$ and $T_2(x_1,x_2)$ on a test set involving the use of a grid of $20 \times 20$ Chebyshev-Lobatto distributed points. The investigation encompassed 100, 200, 300, and 400 independent runs, providing a thorough evaluation of the result solution variability. %The analysis primarily aimed to capture central tendencies, achieved through median calculations, and scrutinized outcome dispersion via 5th and 95th percentile confidence intervals. This meticulous approach facilitated a comprehensive understanding of the reliability and stability of the PINN model across numerous simulations, offering insights into the robustness of the considered methodology.\\

\begin{figure}[htbp]
 \centering
    \begin{subfigure}{0.45\textwidth}
        \includegraphics[width=\linewidth]{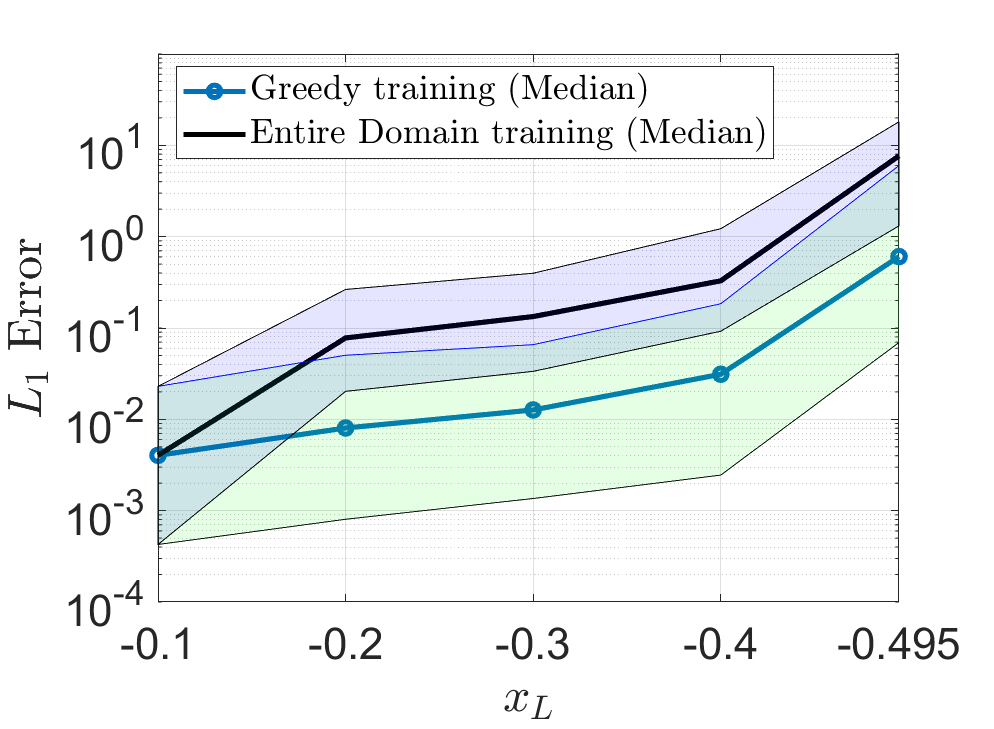}
        \caption{Greedy vs Single training, $L_1$ norm }
    \end{subfigure}
 \hfill
\begin{subfigure}{0.45\textwidth}
        \includegraphics[width=\linewidth]{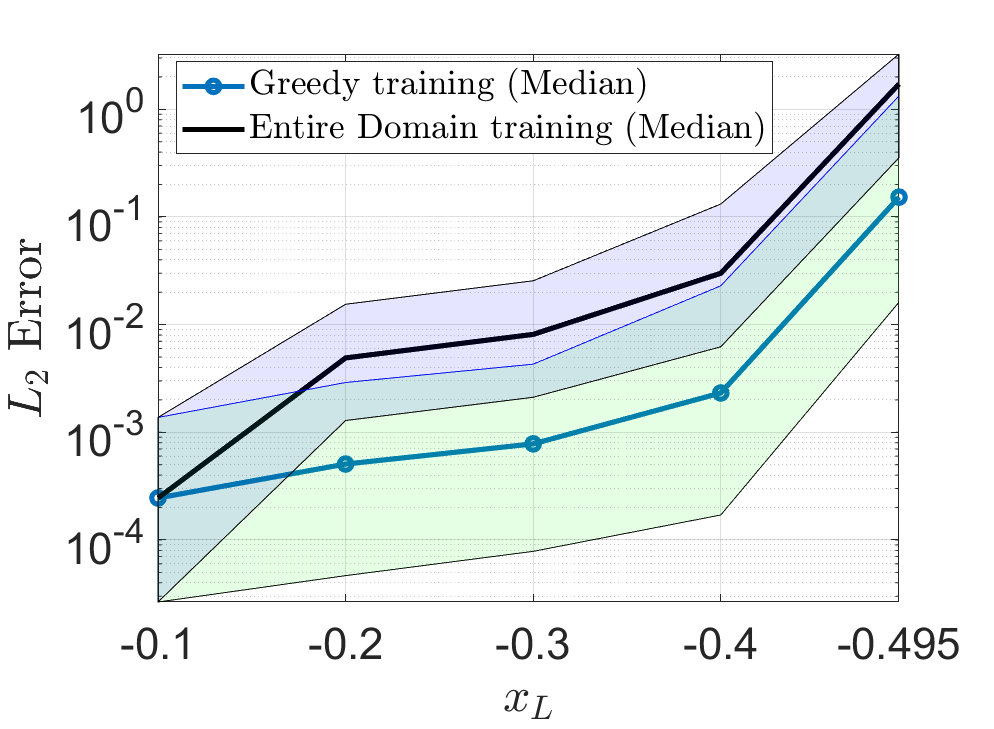}
        \caption{Greedy vs Single training, $L_2$ norm}
\end{subfigure}
 \hfill
\begin{subfigure}{0.45\textwidth}
        \includegraphics[width=\linewidth]{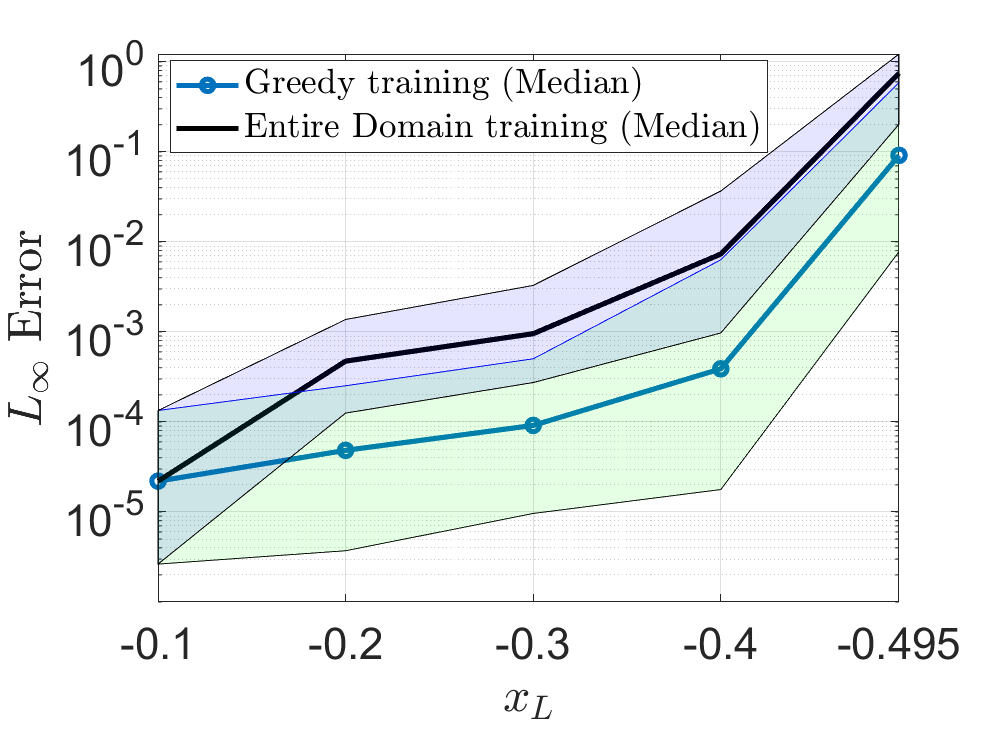}
        \caption{Greedy vs Single training, $L_\infty$ norm}
 \end{subfigure}
 \hfill
 \begin{subfigure}{0.45\textwidth}
     \includegraphics[width=\linewidth]{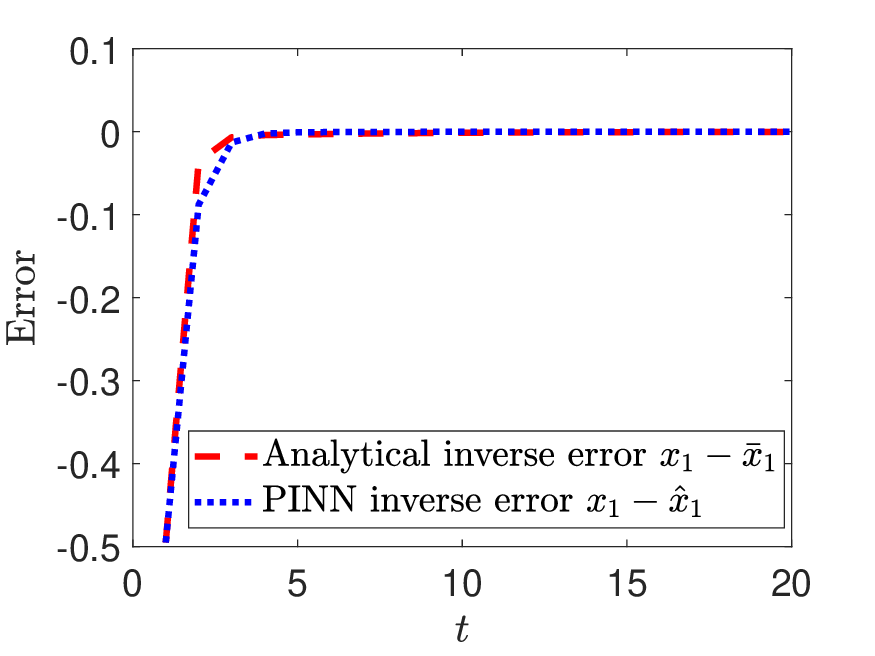}
     \caption{Analytical vs PINN-based inverse error}
 \end{subfigure}
    \caption{Benchmark Problem 1: (a)-(c) Panels involve a comprehensive comparison between the performance of Physics-Informed Neural Networks (PINN) trained using two distinct approaches: one trained on the entire domain at once (Black line) and another trained using a greedy approach (blue line). The evaluation is conducted on a test set comprising $20\times20$ Chebyshev nodes of the second type. The reported performance metrics include the $L_1$, $L_2$, and $L_\infty$ error norms, presented in terms of the median and the 5th to 95th percentiles. The evaluation is performed over 100 independent runs for the transformation $T_1(x_1, x_2)$. Panel (d) depicts a comparison between observation errors. The red line represents the difference between the analytical inverse $\Bar{x}_1$ of the transformation $T_1(x_1, x_2)$ and the real state $x_1$, whereas the blue line depicts the difference between the true state $x_1$ and its approximation obtained through Newton's method, denoted as $\hat{x}_1$.} 
    \label{fig:L-errors_BP1}
\end{figure}

Furthermore, Figure (\ref{fig:L-errors_BP1}) depicts the error characteristics between the analytical solution and the outcome generated by the PINN for $T_1(x_1, x_2)$. In this benchmark problem, the transformation map is quite important as it encompasses the singular point and relates to the state assumed to be unmeasurable. In Figure (\ref{fig:L-errors_BP1})(a), we juxtapose the solution approximations obtained through training across the entire domain and in a greedy-wise manner. The $L_1$ norm for both PINN approximations is presented; the blue line signifies the median of the $L_1$ error evaluated at specific points in the domain. Additionally, the green band represents a confidence interval spanning from the 5th to the 95th percentile, reflecting a PINN trained in a greedy-wise fashion. On the other hand, the black line depicts the median of the $L_1$ error when the network is trained across the entire domain at once, and the green band denotes a confidence interval ranging from the 5th to the 95th percentile.

Figure (\ref{fig:L-errors_BP1})(b), depicts the $L_2$ norm; the blue line depicts the median of the $L_2$ error with a confidence interval ranging from the 5th to the 95th percentile represented by the green band; the black line shows again the median of the $L_2$ error when the network is trained across the entire domain at once, considering the blue band as a confidence interval in the same sense as before. Figure (\ref{fig:L-errors_BP1})(c) displays results for the $L_\infty$ norm.
It can be inferred that the greedy-wise training significantly outperforms the one obtained when training the network across the entire domain in a single step. The numerical approximation from the single training exhibits a poorer performance, particularly in regions proximate to the singular point, as is also evident in Figure (\ref{fig:Modelavailable_TS}).\par

\begin{table}[ht!]
    \centering
    \rowcolors{1}{mygray}{white}
    \caption{Benchmark problem 1. Test set using $20 \times 20$ Chebyshev-Lobatto distributed points $L_1$, $L_2$, $L_{\infty}$ error norms (median, 5th- 95th percentiles) between the analytical and computed solutions of $T_1(x_1,x_2)$ and $T_2(x_1,x_2)$ of the various  schemes trained in the entire domain and with the greedy approach, over 100 runs.}
    \begin{tabular}{||c|c| c c c ||}
    \toprule
    \toprule
   & Error  & Power-Series & PINN Entire domain & PINN Greedy \\
    $T(x_1,x_2)$& norm &$6th$ order &Median  &Median \\
    &&-&(5th,95th) &(5th,95th) 
    \\
    \midrule
    \midrule
    & $\lVert \cdot \rVert_1$ &6.89E$+$01 &1.186E$+$01,(1.048E$+$01,1.372E$+$01) &6.051E$-$01,(6.94E$-$02,6.00E$+$00) \\
    $T_1(x_1,x_2)$ & $\lVert \cdot \rVert_{2}$ &4.55E$+$00  & 1.062E$+$01,(9.94E$+$00,1.243E$+$01) & 1.528E$-$01,(1.62E$-$02,1.33E$+$00)\\
    & $\lVert \cdot \rVert_{\infty}$ &2.88E$+$00  & 8.28E$+$00,(7.00E$+$00,9.63E$+$00)&9.15E$-$02,(7.80E$-$03,5.87E$-$01)\\
    \midrule
    & $\lVert \cdot \rVert_1$ & 0 &2.44E$+$00,(7.40E$-$01,4.00E$+$00)&3.595E$-$01,(4.39E$-$02,2.90E$+$00)\\
   $T_2(x_1,x_2)$ & $\lVert \cdot \rVert_{2}$ & 0 &2.34E$+$00,(6.68E$-$01,3.75E$+$00)&6.31E$-$02,(1.02E$-$02,4.988E$-$01)
   \\& $\lVert \cdot \rVert_\infty$ & 0 & 1.65E$+$00,(4.55E$-$01,2.68E$+$00)&2.94E$-$02,(3.90E$-$03,1.85E$-$01)\\ 
    \bottomrule
    \bottomrule
    \end{tabular}
    \label{tab:1}
\end{table}

\begin{table}[ht!]
    \centering
    \scriptsize
    \rowcolors{1}{mygray}{white}
    \caption{Benchmark problem 1. Uncertainty quantification systematically conducted on a test set employing a grid of $20 \times 20$ Chebyshev-Lobatto distributed points to investigate the convergence behavior of the Physics-Informed Neural Networks (PINN). The study involves performing 100, 200, 300, and 400 independent runs to thoroughly assess the variability in the results. The analysis focused on capturing the central tendency by calculating the median and examining the dispersion of outcomes through the 5th and 95th percentile confidence intervals.}
    \begin{tabular}{||c|c| c c c c ||}
    \toprule
    \toprule
   & Error  & 100 runs & 200 runs & 300 runs & 400 runs \\
    $T(x_1,x_2)$& norm &Median&Median &Median  &Median \\
    &&(5th,95th)&(5th,95th) &(5th,95th) &(5th,95th)
    \\
    \midrule
    \midrule
    & $\lVert \cdot \rVert_1$ &6.051E$-$01,&6.34E$-$01,&6.27E$-$01,&6.12E$-$01,\\
    &&(6.94E$-$02,6.00E$+$00)&(6.20E$-$02,5.60E$+$00)&(6.188E$-$02,5.45E$+$00)&(6.12E$-$02,5.32E$+$00)
     \\
    $T_1(x_1,x_2)$ & $\lVert \cdot \rVert_{2}$ &1.528E$-$01,&1.77E$-$01,&1.72E$-$01,&1.705E$-$01,\\
    &&(1.62E$-$02,1.33E$+$00)&(1.95E$-$02,1.23E$+$00)&(1.88E$-$02,1.21E$+$00)&(1.76E$-$02,1.18E$+$00)
    \\
    & $\lVert \cdot \rVert_{\infty}$ &9.15E$-$02,&8.87E$-$02,&8.75E$-$02,&8.66E$-$02,\\
    &&(7.80E$-$03,5.87E$-$01)&(6.58E$-$03,5.12E$-$01)&(6.55E$-$03,5.10E$-$01)&(6.502E$-$03,5.108E$-$01)
    \\
    \midrule
    & $\lVert \cdot \rVert_1$ &3.595E$-$01,&3.90E$-$01,&3.84E$-$01,&3.78E$-$01,\\
    &&(4.39E$-$02,2.90E$+$00)&(4.98E$-$02,3.25E$+$00)&(4.85E$-$02,3.21E$+$00)&(4.78E$-$02,3.16E$+$00)
    \\
   $T_2(x_1,x_2)$ & $\lVert \cdot \rVert_{2}$ &6.31E$-$02,&6.21E$-$02,&6.19E$-$02,&6.15E$-$02,
   \\
   &&(1.02E$-$02,4.988E$-$01)&(1.80E$-$02,6.01E$-$01)&(1.75E$-$02,6.00E$-$01)&(1.69E$-$02,5.98E$-$01)
   \\
   & $\lVert \cdot \rVert_\infty$ &2.94E$-$02,&2.81E$-$02,&2.75E$-$02,&2.701E$-$02,\\ 
   &&(3.90E$-$03,1.85E$-$01)&(3.28E$-$03,2.10E$-$01)&(3.22E$-$03,2.03E$-$01)&(3.19E$-$03,1,98E$-$01)
   \\
    \bottomrule
    \bottomrule
    \end{tabular}
    \label{tab:2}
\end{table}

Finally, Figure (\ref{fig:L-errors_BP1})(d) depicts the difference between the analytical inverse map $T_1^{-1}(x_1,x_2)$ \eqref{AinverseEx2} of the observed variable \eqref{observerEx2} and the numerical approximation of the inverse $\hat{T}_1^{-1}(x_1,x_2)$ as computed by Newton's method (presented in Section \ref{sec:framework}). As shown, the difference between the true state and the inverse transformation value generated by the PINN converges to zero. For this particular simulation, we have initialized the $z$-states as follows: $z_1(1)=0$ and $z_2(1)=0$, and we have used as initial guess for Newton's method the following values for the $x$-states:  $x_0 = [0.1,0.1]$. Moreover, we have initialized the analytical inverse map expression using: $\hat{x}_1(1) = 0$ and $\hat{x}_2(1) = 0$; the actual states were: $x_1=-0.495$ and $x_2=0.35$. Regarding the Newton's method, we have set the relative and absolute tolerances at: $tol<1E$-$06$ (the scheme converged in $5$ iterations).
% \begin{figure}[ht!]
%     \centering
%     \includegraphics[width=0.7\linewidth]{Error_Ex1.eps}
%     \caption{Comparison between observation errors. The red line represents the difference between the analytical inverse $\Bar{x}_1$ of the transformation $T_1(x_1, x_2)$ and the real state $x_1$, whereas the blue line depicts the difference between the true state $x_1$ and its approximation obtained through Newton's method, denoted as $\hat{x}_1$.}    \label{fig:Erroydynamics1}
% \end{figure}

%%%%%%%%%%

%%%%%%%%%%% Beggining of Benchmark problem 2
\newpage
\subsection{Benchmark Problem 2}
\begin{figure}[htbp]
 \centering
 \begin{subfigure}[h]{0.45\textwidth}
     \centering
     \includegraphics[width=7.5cm]{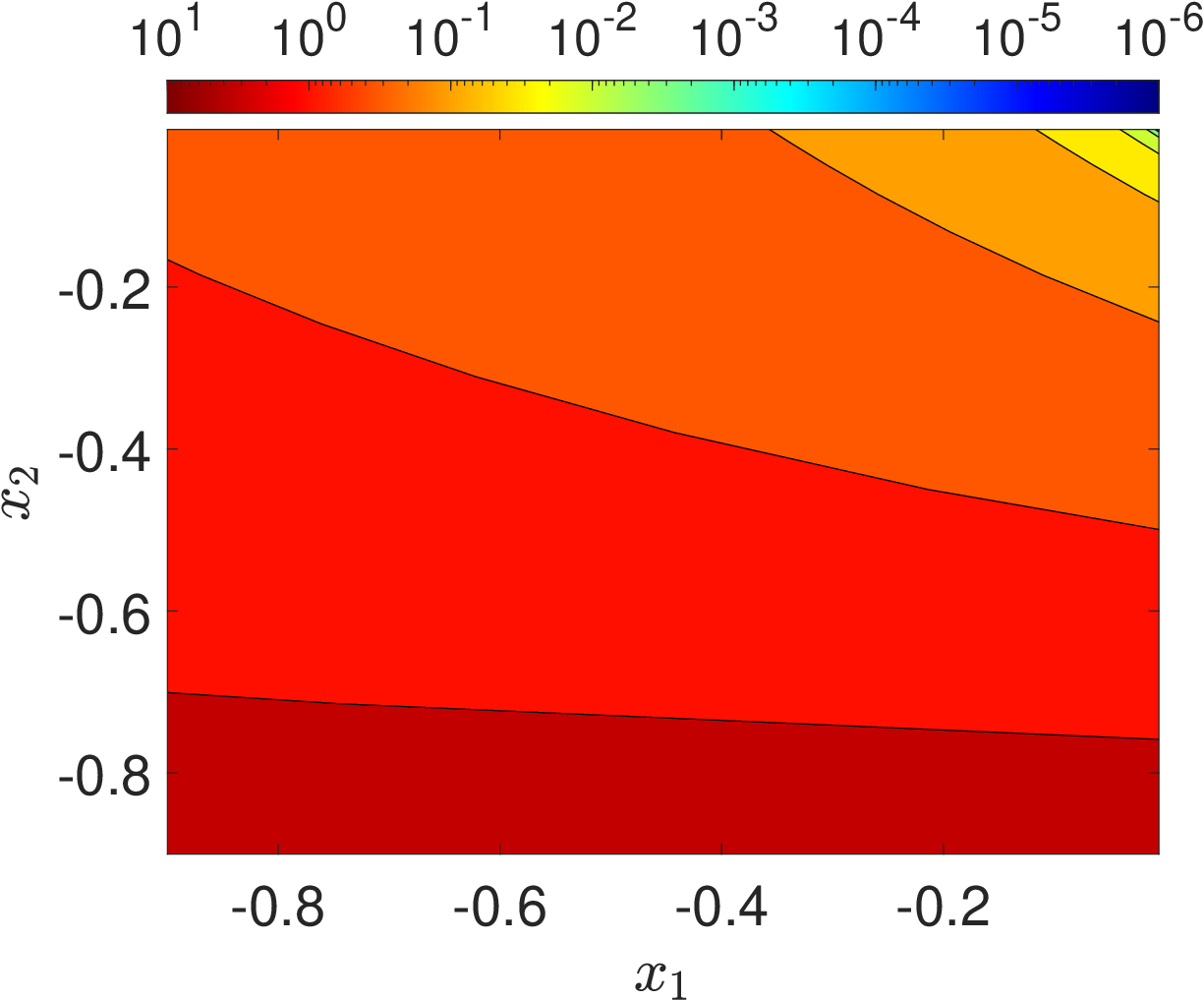}
     \caption{}
 \end{subfigure}
 \hfill
 \begin{subfigure}[h]{0.5\textwidth}
     \centering
     \includegraphics[width=7.5cm]{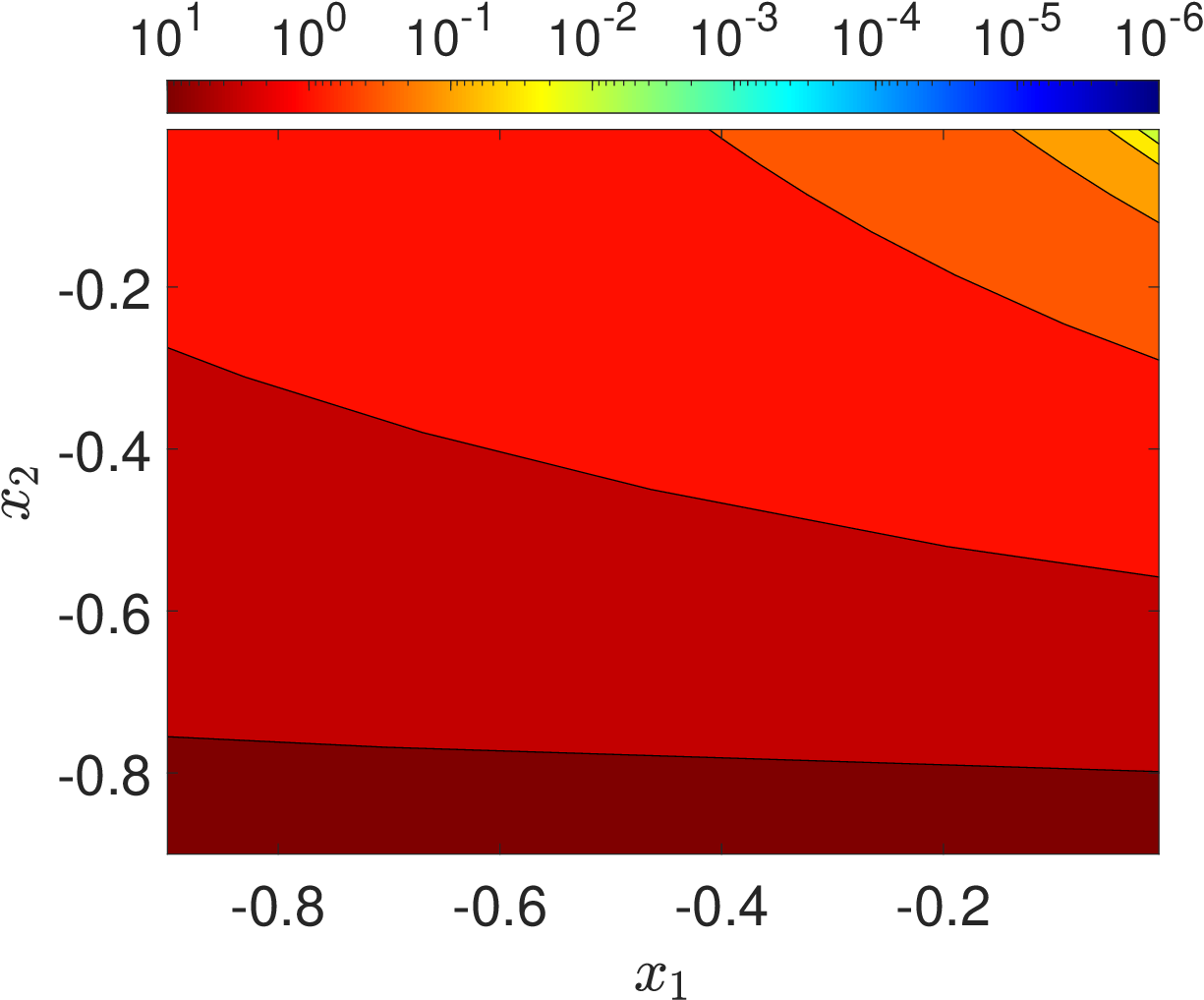}
     \caption{}
 \end{subfigure}
 \begin{subfigure}[h]{0.45\textwidth}
     \centering
     \includegraphics[width=7.2cm]{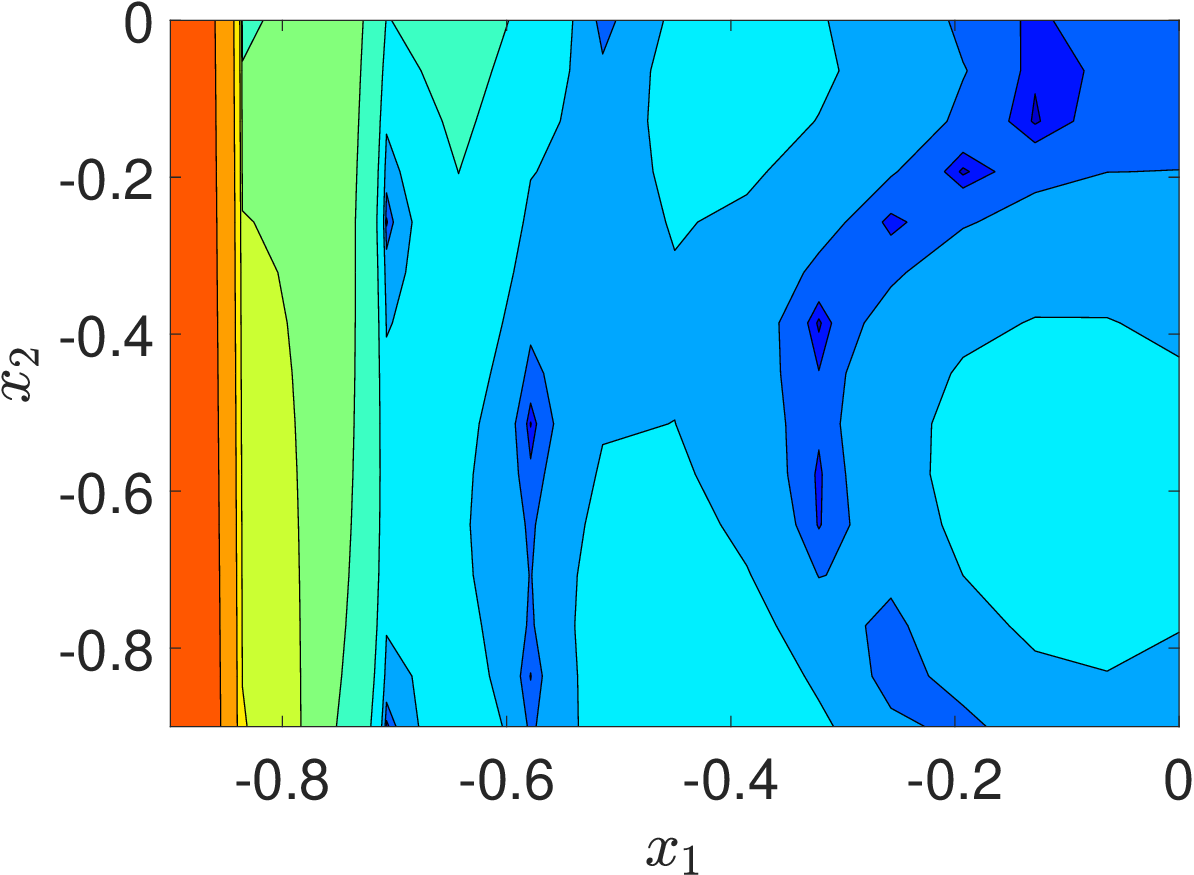}
     \caption{}
 \end{subfigure}
 \hfill
 \begin{subfigure}[h]{0.5\textwidth}
     \centering
     \includegraphics[width=7.2cm]{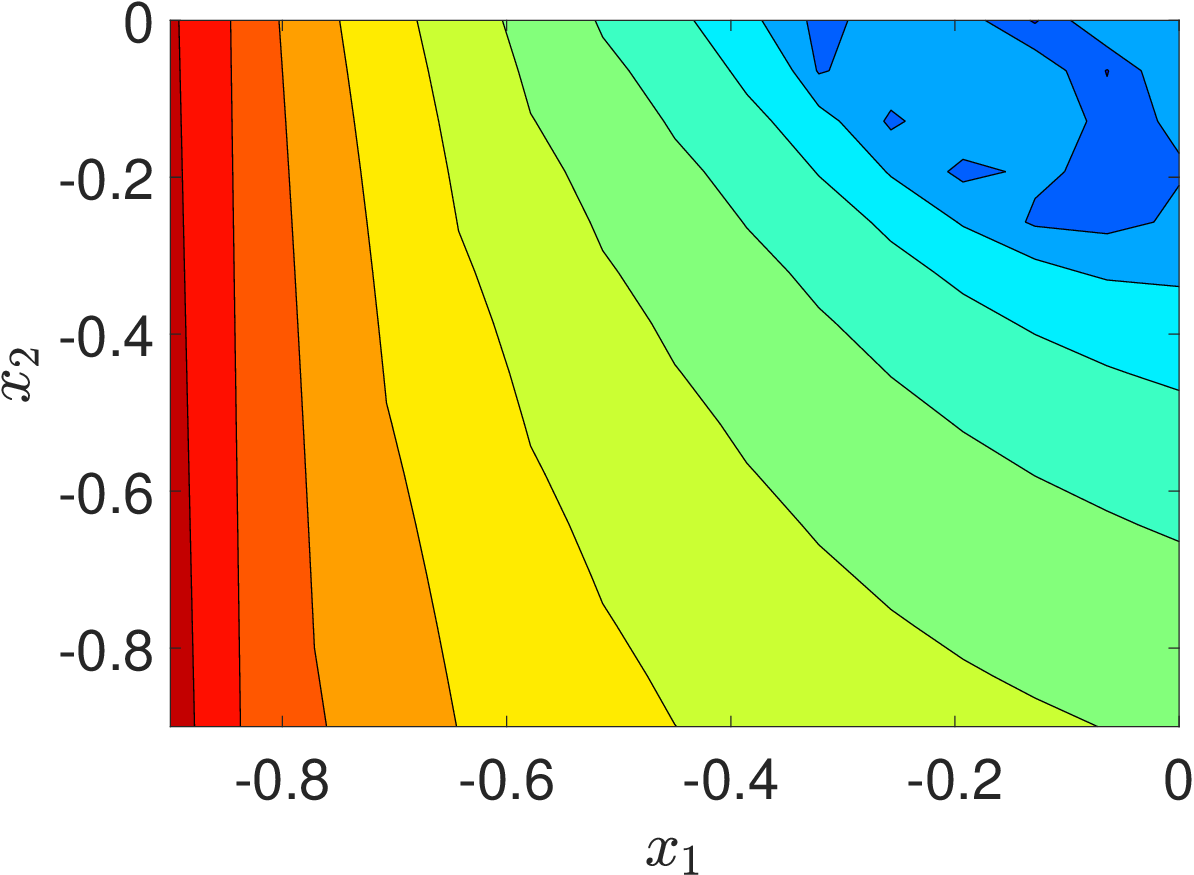}
     \caption{}
 \end{subfigure}
   \centering
 \begin{subfigure}[h]{0.45\textwidth}
     \centering
     \includegraphics[width=7.2cm]{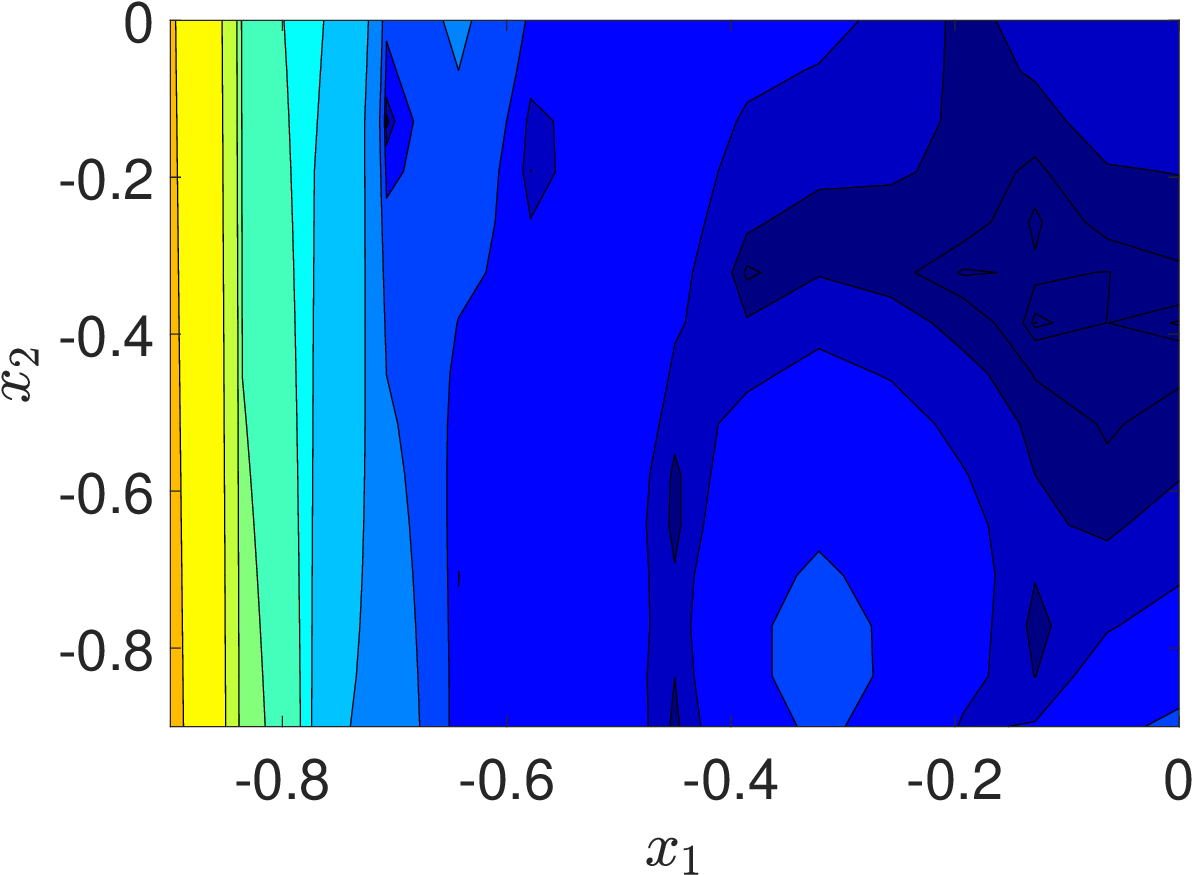}
     \caption{}
 \end{subfigure}
 \hfill
 \begin{subfigure}[h]{0.5\textwidth}
     \centering
     \includegraphics[width=7.2cm]{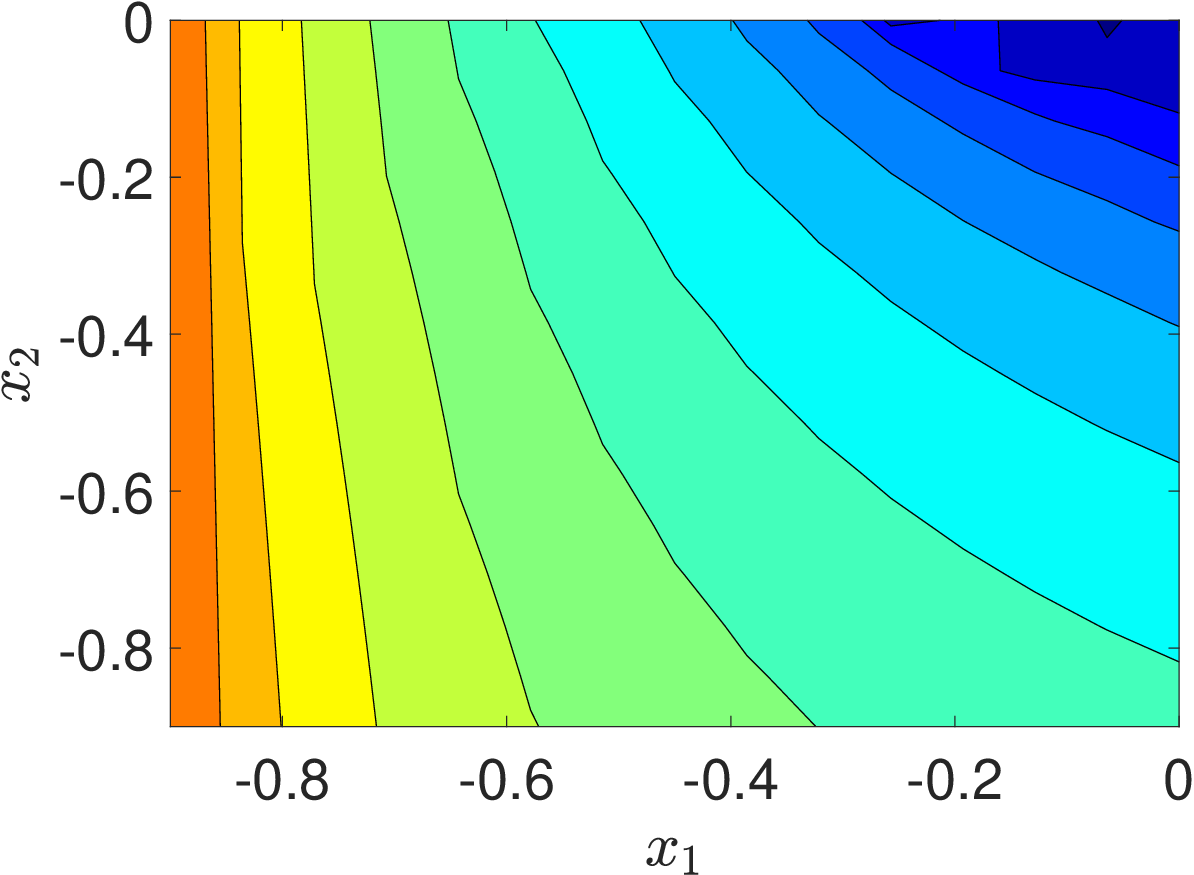}
     \caption{}
 \end{subfigure}
\caption{Training sets (grids of $15 \times 15$ equispaced distributed points). Numerical approximation accuracy (difference between the computed and analytical solution) of $T_1(x_1,x_2)$ (left column) and $T_2(x_1,x_2)$ (right column) using the various schemes. (a),(b) $6th$ order power-series expansion of $T_1(x_1,x_2)$ and $T_2(x_1,x_2)$ and the right-hand side of the model (\ref{eq:DiscreteSystem}) in $[-0.91,0]\times[-0.91,0]$. (c),(d) PINN  trained over the entire domain $[-0.91,0]\times[-0.91,0]$. (e),(f) PINN  trained via the greedy-wise approach.}
\label{fig:Modelavailable_TrS_Ex2}
\end{figure}

This particular benchmark problem introduces the challenge of unavailability of the state $x_2$, with only the state $x_1$ being measurable. Consequently, the focus of this problem is on the reconstruction of the second state. As highlighted earlier, the unmeasured state $x_2$ is viewed as either an unidentified system parameter or unknown disturbance, thus adding complexity to the estimation task. 
As a consequence and in alignment with the methodology outlined so far, the task is to address the set of functional equations \eqref{FEs2} utilizing the established PINN scheme. Figure (\ref{fig:Modelavailable_TrS_Ex2}) depicts the numerical approximation accuracy, as the difference between computed and analytical solutions \eqref{Nonlinear_Transformations2}, obtained by various schemes for the training set. This particular benchmark problem is more challenging as the gradients it contains (see Fig. \ref{fig:Analytical_transformation_ex2}) are notably bigger compared to the previous one.\par 
Figures (\ref{fig:Modelavailable_TrS_Ex2})(a),(b) demonstrate the results achieved with a 6th-order truncation in the power series expansion of the nonlinear transformation and the right-hand side of the discrete model. In this particular problem, the steep gradient emanating from the singular point at $x_1=-1$ poses a significant challenge as mentioned before. The numerical approximation provided by the series expansion method is notably poor across the domain, particularly of the order of $10^1$ at the grid's edge, as mentioned previously, near the singular point.\par 
Furthermore, Figures (\ref{fig:Modelavailable_TrS_Ex2})(c),(d) depict the outcomes obtained with the PINN approach to learn the transformation maps \eqref{Nonlinear_Transformations2} across the entire domain without the greedy approach. The approximation accuracy of this scheme is relatively modest, especially in the vicinity of the singularity, approximately of the order of $10^{0}$ for $T_1(x_1,x_2)$ and of the order of $10^{1}$ for $T_2(x_1,x_2)$. Nevertheless, it outperforms the conventional Taylor series expansion method in terms of numerical error accuracy across the entire domain in general.\par
\begin{figure}[htbp]
 \centering
 \begin{subfigure}[h]{0.45\textwidth}
     \centering
     \includegraphics[width=7.5cm]{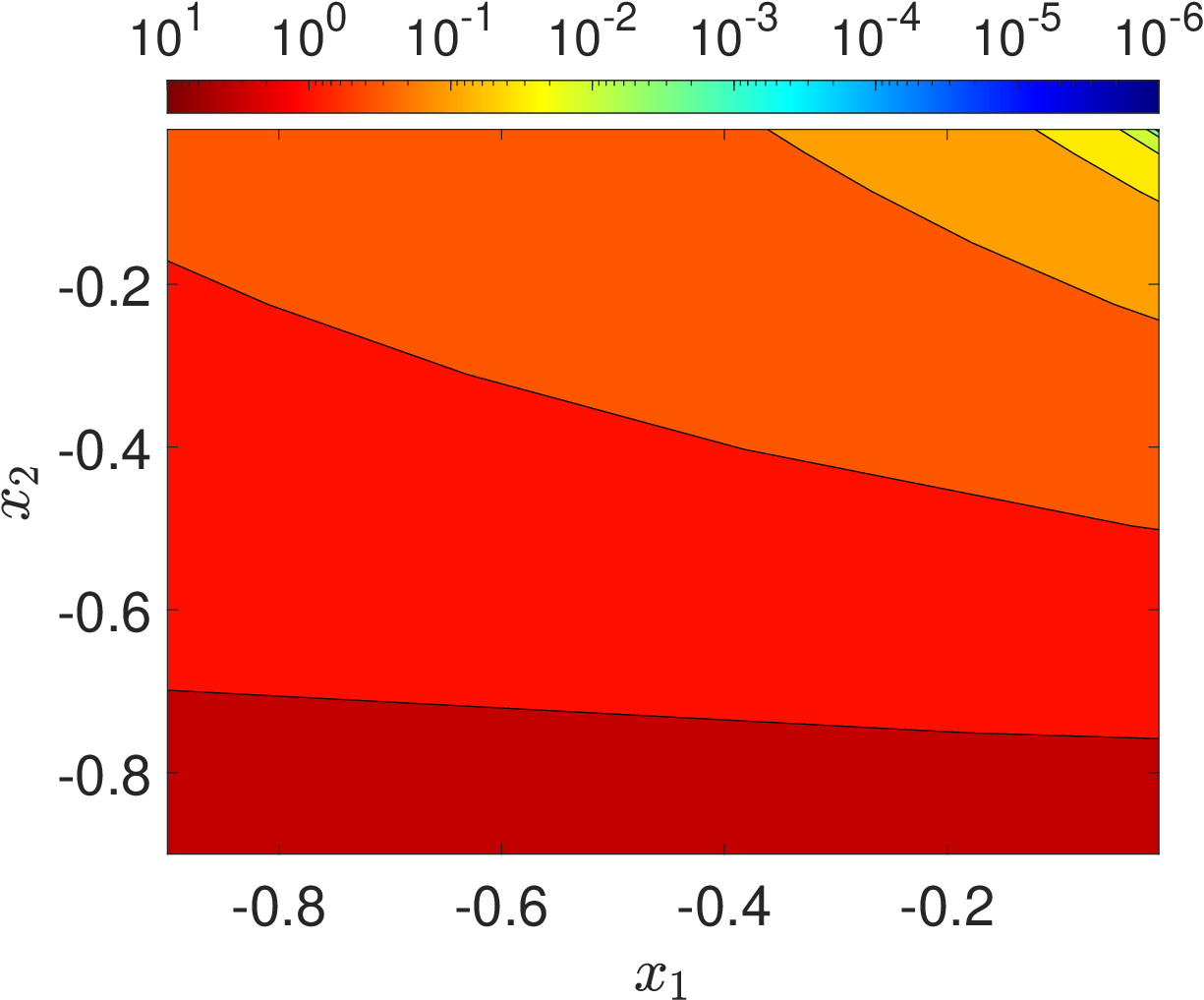}
     \caption{}
 \end{subfigure}
 \hfill
 \begin{subfigure}[h]{0.5\textwidth}
     \centering
     \includegraphics[width=7.5cm]{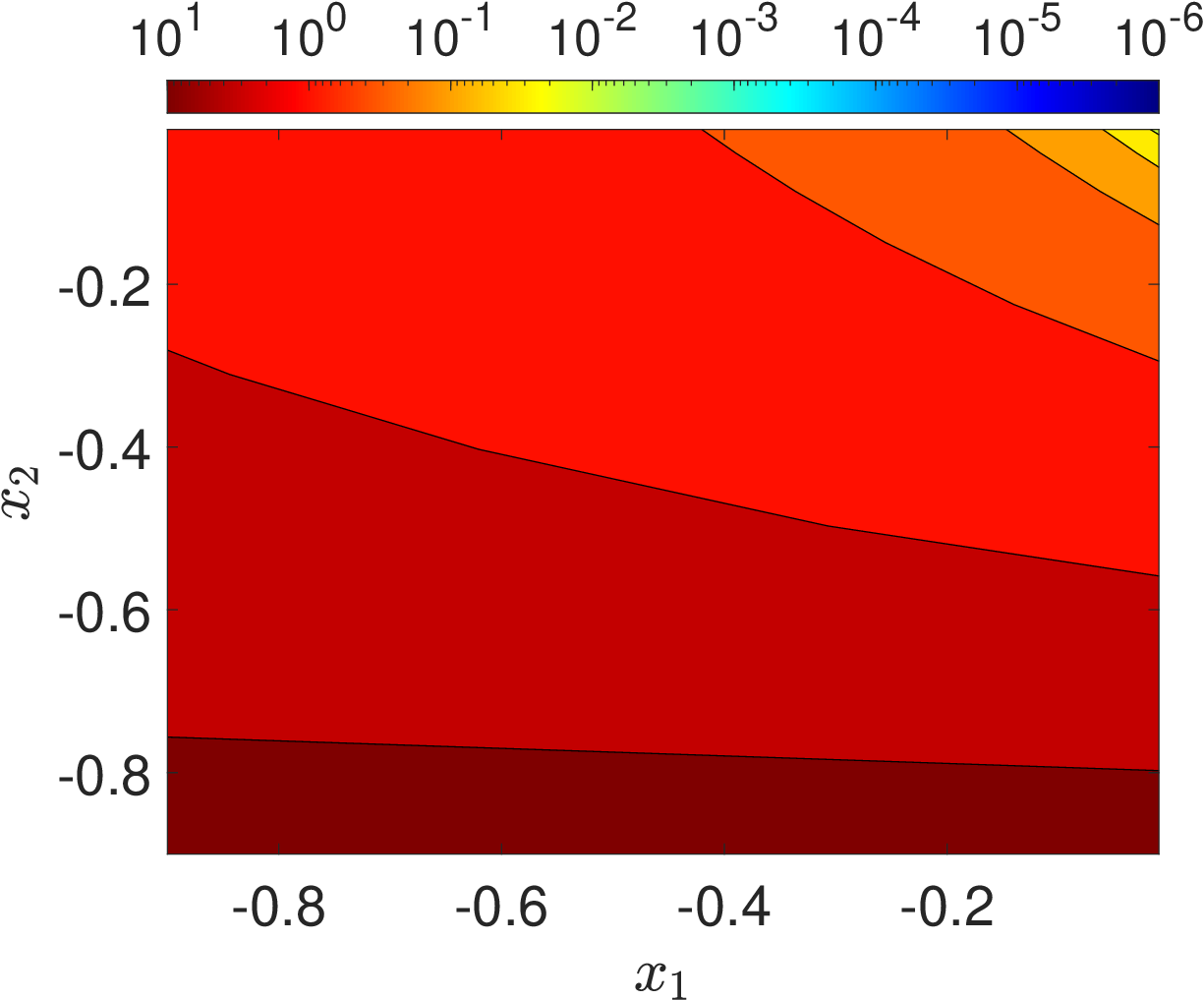}
     \caption{}
 \end{subfigure}
 \begin{subfigure}[h]{0.45\textwidth}
     \centering
     \includegraphics[width=7.2cm]{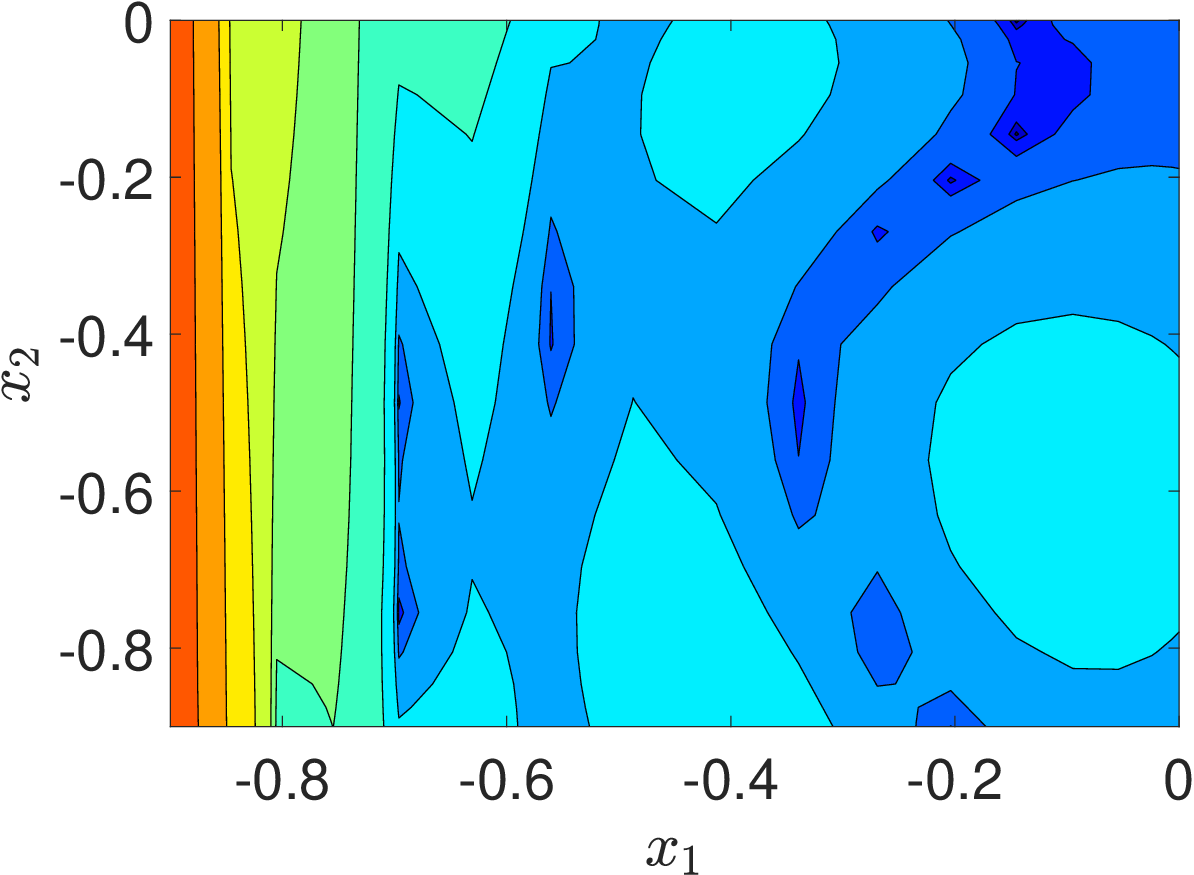}
     \caption{}
 \end{subfigure}
 \hfill
 \begin{subfigure}[h]{0.5\textwidth}
     \centering
     \includegraphics[width=7.2cm]{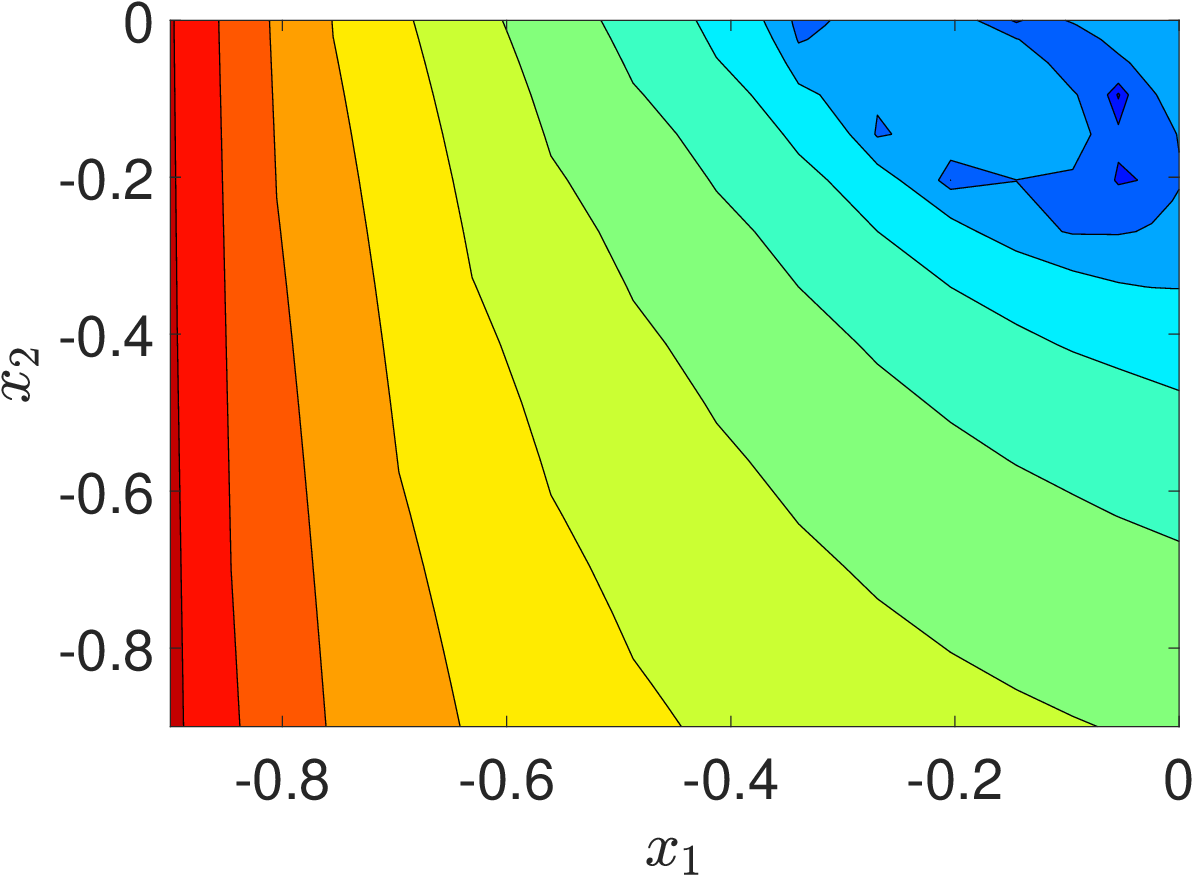}
     \caption{}
 \end{subfigure}
   \centering
 \begin{subfigure}[h]{0.45\textwidth}
     \centering
     \includegraphics[width=7.2cm]{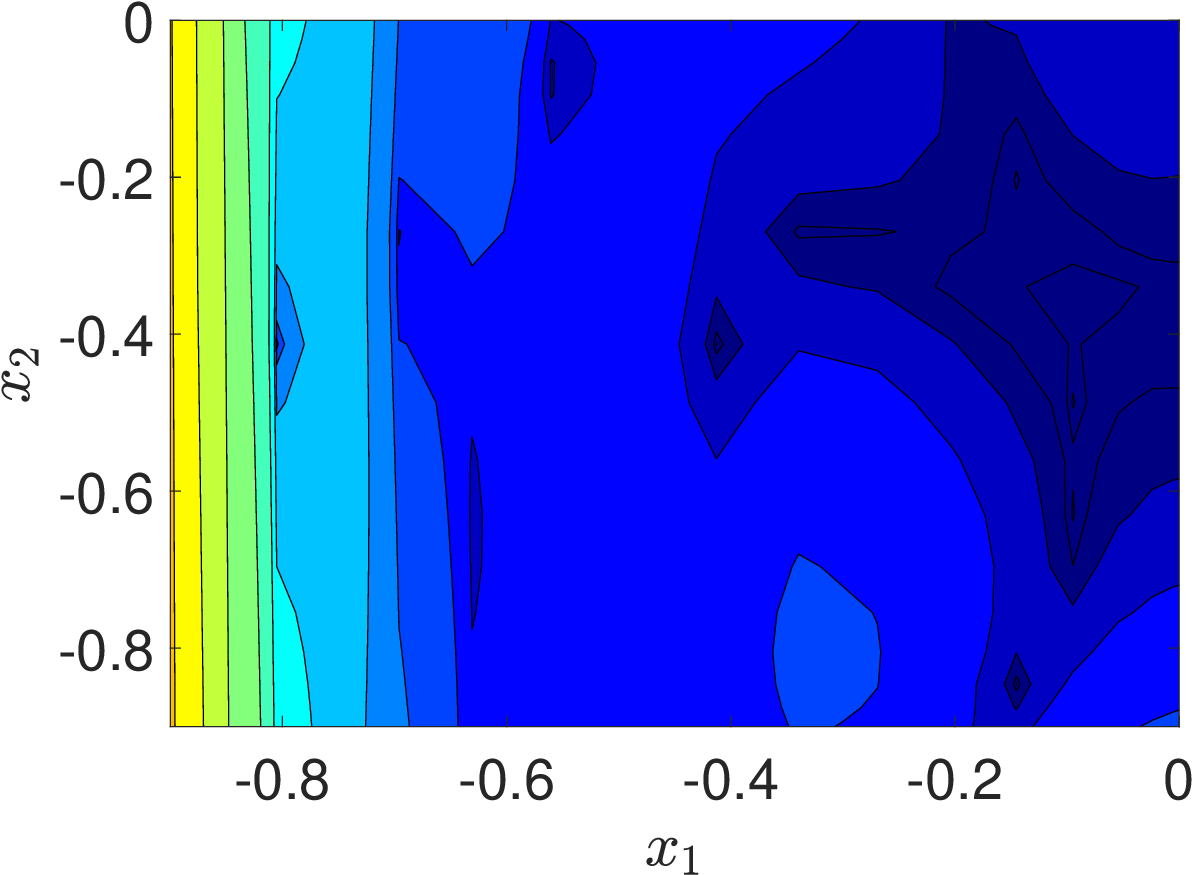}
     \caption{}
 \end{subfigure}
 \hfill
 \begin{subfigure}[h]{0.5\textwidth}
     \centering
     \includegraphics[width=7.2cm]{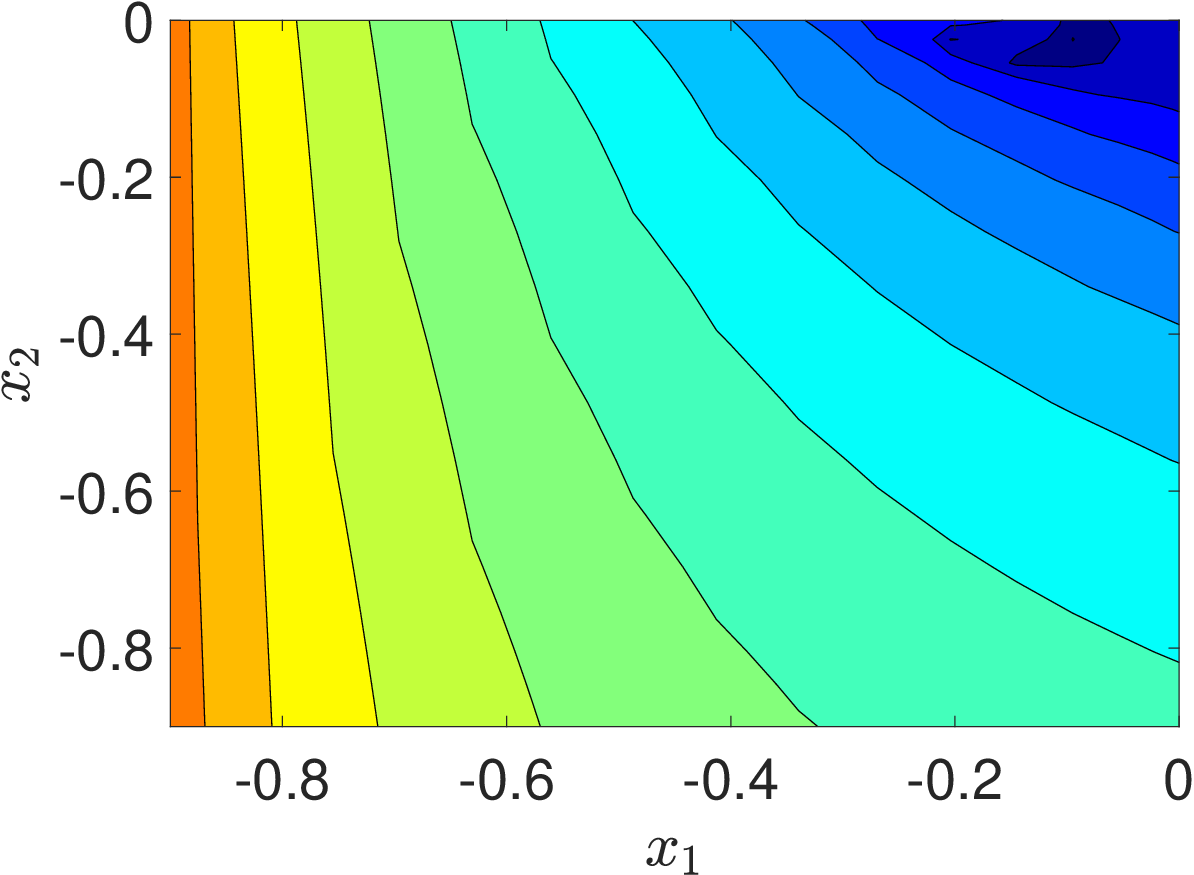}
     \caption{}
 \end{subfigure}
\caption{Test sets (grids of $20 \times 20$ Chebyshev-Lobatto distributed points). Numerical approximation accuracy (difference between the computed and analytical solution) of $T_1(x_1,x_2)$ (left column) and $T_2(x_1,x_2)$ (right column) using the various schemes. (a),(b) $6th$ order power-series expansion of $T_1(x_1,x_2)$ and $T_2(x_1,x_2)$ and the right-hand side of the model (\ref{eq:DiscreteSystem}) in $[-0.91,0]\times[-0.91,0]$. (c),(d) PINN trained over the entire domain $[-0.91,0]\times[-0.91,0]$. (e),(f) PINN trained via the greedy-wise approach.}
\label{fig:Modelavailable_TS_Ex2}
\end{figure}

 
Figures (\ref{fig:Modelavailable_TrS_Ex2})(e),(f) present the results obtained from the PINN implementation with the greedy approach. Notably, the adoption of the greedy-wise training method substantially enhances the numerical approximation accuracy when compared to PINN schemes trained across the entire domain in a single training step. Here, the numerical approximation error is of the order of $10^{-2}$ in the region close to the edge point $(-0.91,-0.91)$. Additionally, Figure (\ref{fig:Modelavailable_TS_Ex2}) portrays the performance of the schemes when tested on a Chebyshev-Lobatto grid, with similar results observed for both the training and test sets.\par
Table (\ref{tab:NormsS2_Ex2}) allows a comparison of the approximation errors between the analytical and numerical solutions. Specifically, the 6th order power series expansions of $T_1(x_1, x_2)$ and $T_2(x_1, x_2)$ with the two PINN solution approximations (with and without the greedy approach) are compared. To test the convergence of the schemes, we used  a Chebyshev-Lobatto grid with $20 \times 20$ points. The reported errors include the median, 5th percentile, and 95th percentile aggregated over 100 runs.\par
Table (\ref{tab:4}),  illustrates the results of uncertainty quantification performed to investigate the convergence of the PINN scheme in approximating the nonlinear transformation maps $T_1(x_1, x_2)$ and $T_2(x_1,x_2)$ on a test set (a grid of $20 \times 20$ Chebyshev-Lobatto distributed points). In our computations, we considered 100, 200, 300, and 400 MC runs. As shown, the medians and variances, are for any practical means of the same order, this indicating convergence of the training process.\par
\begin{figure}[!ht]
    \centering
    \begin{subfigure}{0.33\textwidth}
        \includegraphics[height=4.25cm]{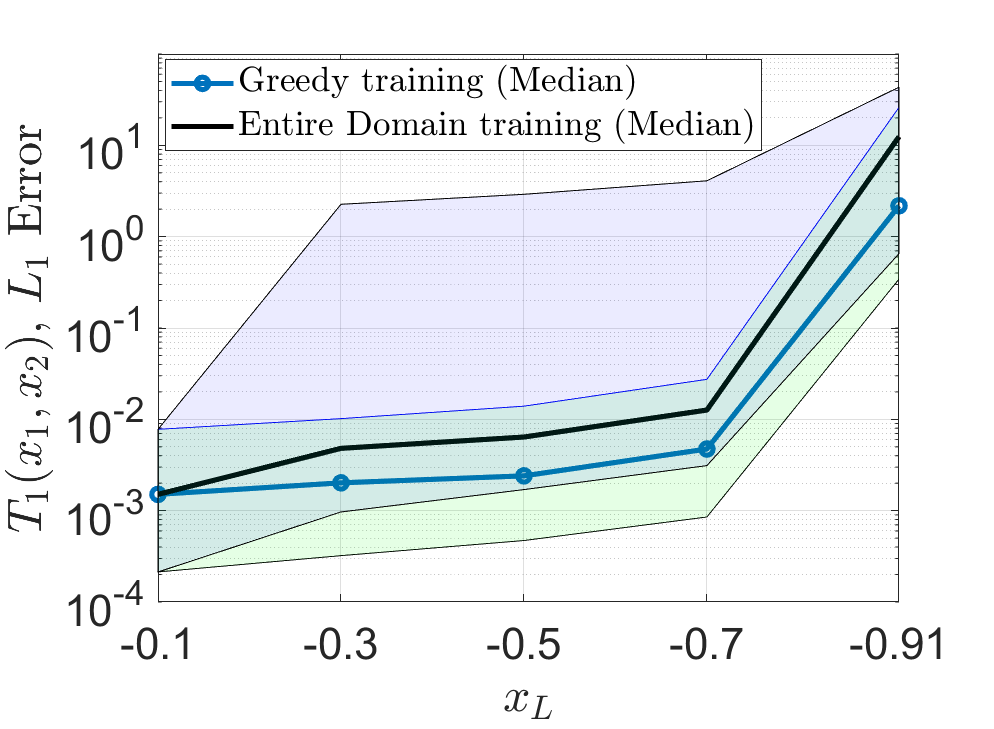}
        \caption{Greedy training, $L_1$ norm }
    \end{subfigure}
    \hfill
    \begin{subfigure}{0.33\textwidth}
        \includegraphics[height=4.25cm]{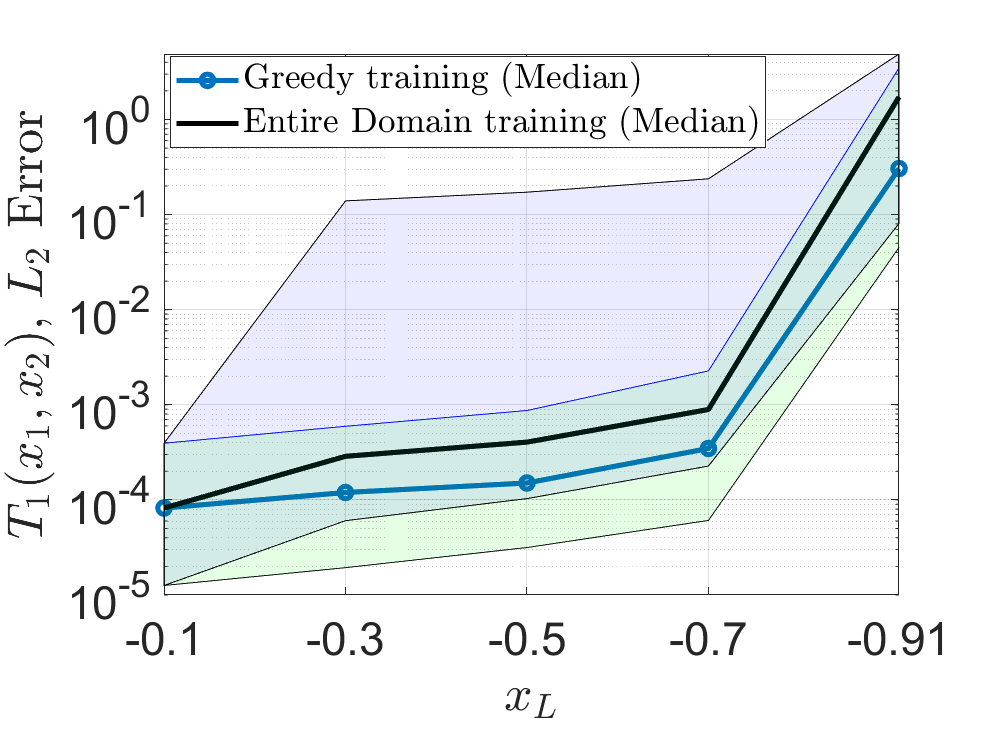}
        \caption{Greedy training, $L_2$ norm}
    \end{subfigure}
    \hfill
    \begin{subfigure}{0.33\textwidth}
        \includegraphics[height=4.25cm]{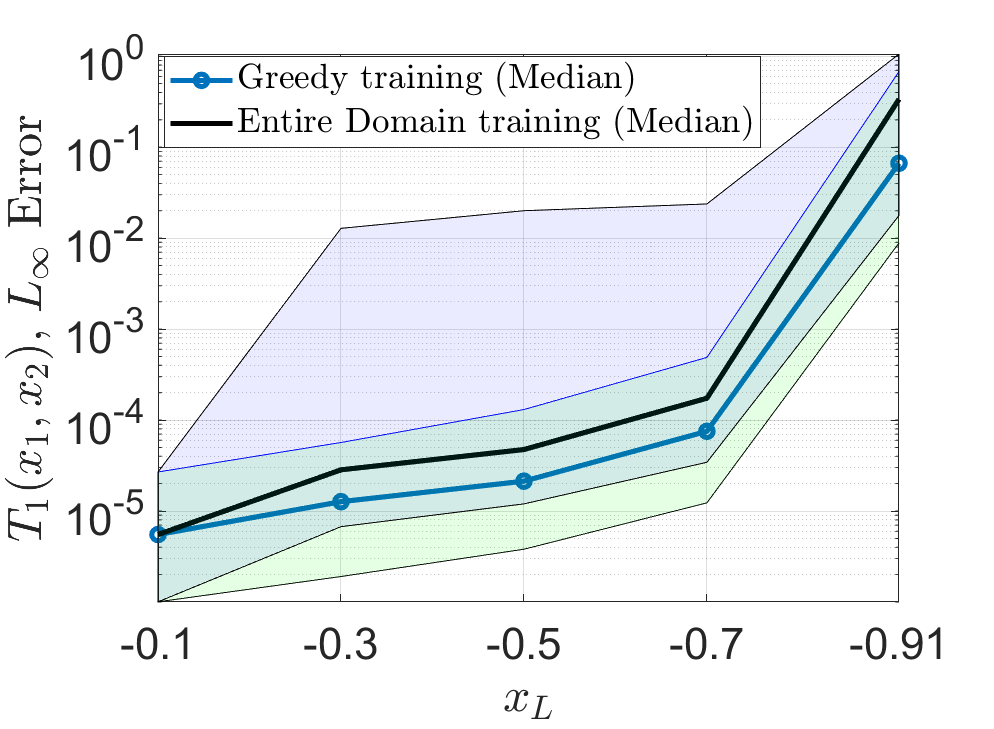}
        \caption{Greedy training, $L_\infty$ norm}
    \end{subfigure}

    \medskip

    \begin{subfigure}{0.33\textwidth}
        \includegraphics[height=4.25cm]{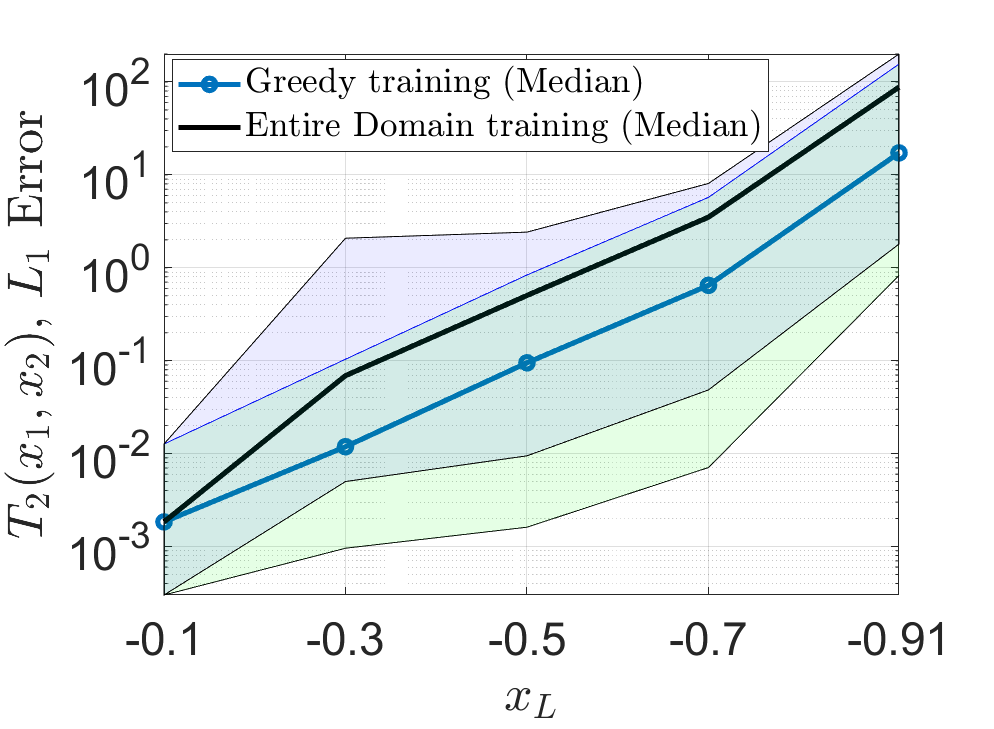}
        \caption{Single training, $L_1$ norm}
    \end{subfigure}
    \hfill
    \begin{subfigure}{0.33\textwidth}
        \includegraphics[height=4.25cm]{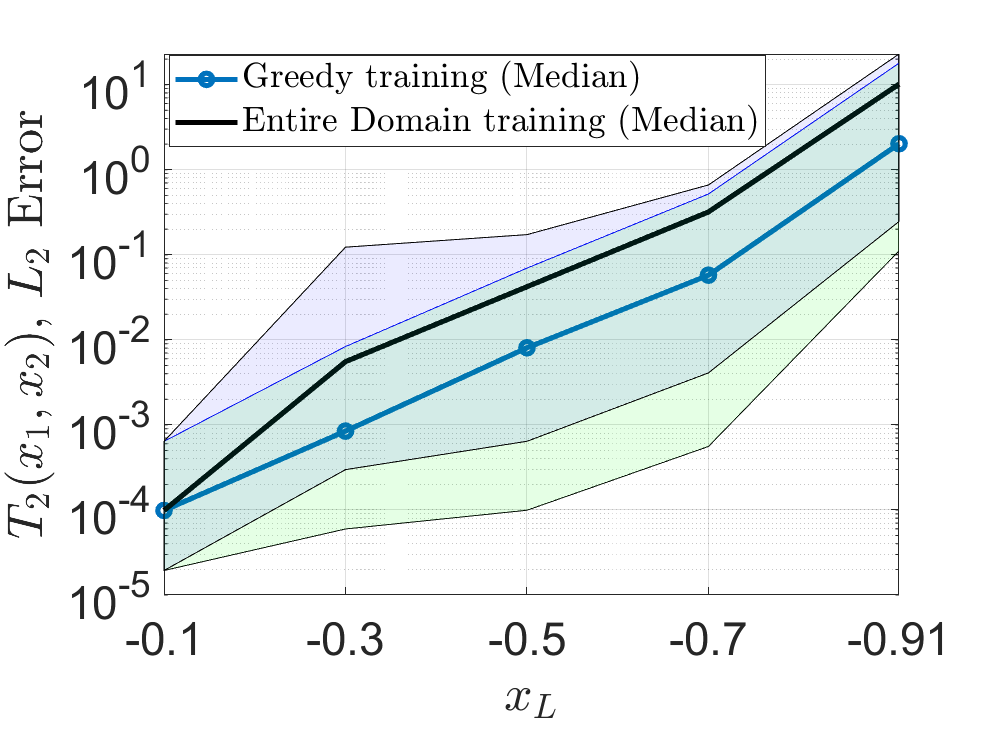}
        \caption{Single training, $L_2$ norm}
    \end{subfigure}
    \hfill
    \begin{subfigure}{0.33\textwidth}
        \includegraphics[height=4.25cm]{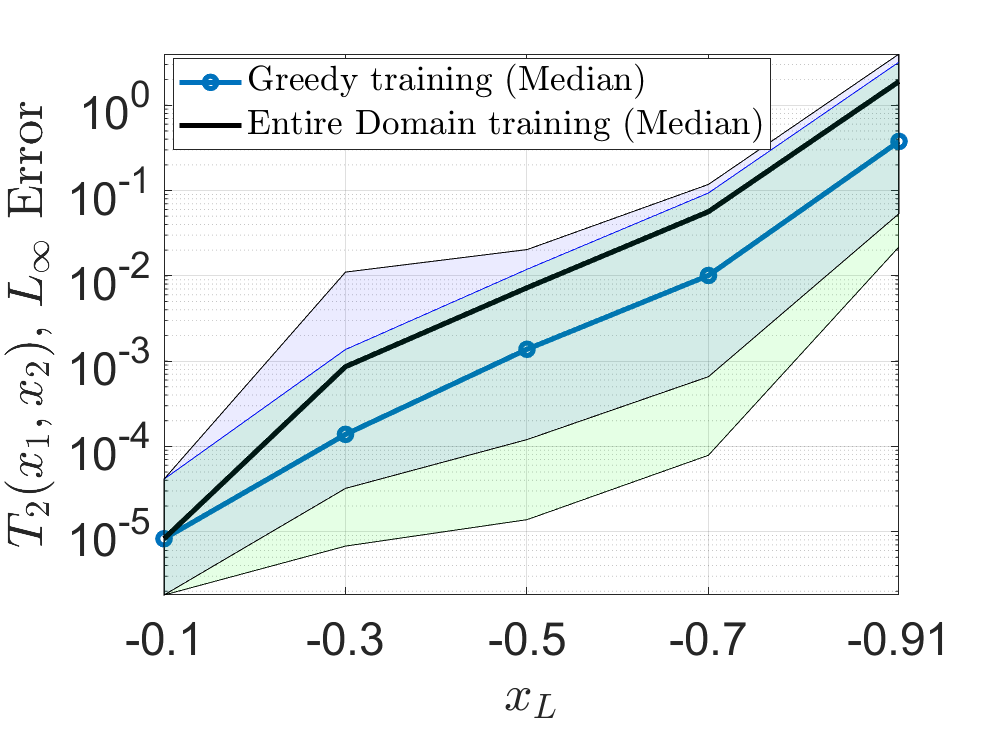}
        \caption{Single training, $L_\infty$ norm}
    \end{subfigure}

    \caption{Benchmark problem 2. (a)-(c) Panels involve a comprehensive comparison between the performance of Physics-Informed Neural Networks (PINN) trained using two distinct approaches: one trained on the entire domain at once (Black line) and another trained using a greedy approach (blue line). The evaluation is conducted on a test set comprising $20\times20$ Chebyshev nodes of the second type. The reported performance metrics include the $L_1$, $L_2$, and $L_\infty$ error norms, presented in terms of the median and the 5th to 95th percentiles. The evaluation is performed over 100 independent runs for the transformation $T_1(x_1, x_2)$. Furthermore Panels (d)-(f) depicts the same error norms with the same configuration for the transformation $T_2(x_1, x_2)$.}
    \label{fig:10}
\end{figure}

 
Figure (\ref{fig:10}) depicts the approximation error behavior between the analytical solution and the one produced by PINN for $T_1(x_1, x_2)$ and  $T_2(x_1, x_2)$. Please notice that $T_2(x_1,x_2)$ is the transformation map of interest in this problem because it is related to the state we seek to observe as mentioned earlier. The approximation accuracy of the PINN trained both across the entire domain and in a greedy-wise manner is assessed. Figure (\ref{fig:10})(a) shows the $L_1$ norm of the transformation map $T_1(x_1, x_2)$. The blue line signifies the median, and the green band represents a confidence interval spanning the 5th to the 95th percentile range for the PINN model trained in a greedy-wise fashion. The black line and the blue band corresponds to the PINN model trained across the entire domain without the greedy approach. Figure (\ref{fig:10})(b)  depicts the corresponding $L_2$ error norms, and Figure (\ref{fig:10})(c) depicts the corresponding $L_\infty$ error norms.  %In these cases, PINN is employed to approximate the transformation $T_1(x_1, x_2)$, The analogous representation across different norm types allows for a comprehensive examination of the model's performance under varying error metrics and as the figure shows the PINN scheme outperforms the performance of the single training in terms of numerical approximation errors.\\
Figure (\ref{fig:10})(d-f) depict similar results for $T_2(x_1, x_2)$.\par
\begin{figure}[!ht]
    \centering    \includegraphics[width=0.550\linewidth]{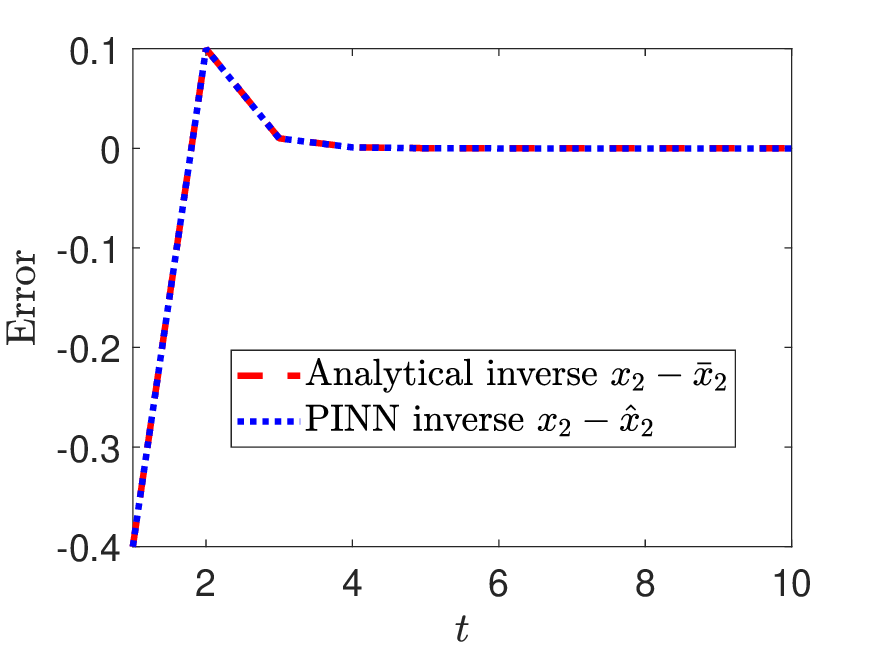}
    \caption{Comparison between observation errors. The red line represents the difference between the analytical inverse $\Bar{x}_2$ of the transformation $T_2(x_1, x_2)$ and the real state $x_2$, whereas the blue line depicts the difference between the true state $x_2$ and its approximation obtained through Newton's method, denoted as $\hat{x}_2$.}    \label{fig:Erroydynamics2}
\end{figure}

It is evident, that the greedy-wise training procedure outperforms the single training one for the approximation of $T_2(x_1, x_2)$.Finally, Figure (\ref{fig:Erroydynamics2}) depicts the error between the trajectories computed by applying the analytically derived inverse $T^{-1}(x_1,x_2)$ \eqref{AinverseEx2} to the observer \eqref{observerEx2} and the numerical approximation of the inverse $\hat{T}^{-1}(x_1,x_2)$ computed by the Newton's method. For this task, we have initialized the $z$-states as: $z_1(1)=0$ and $z_2(1)=0$, and we have used as an initial guess the values $x_0 = [0.1,0.1]$. The analytical inverse expression was initialized using the values: $\hat{x}_1(1) = 1$ and $\hat{x}_2(1) = - 0.4$. Moreover, as mentioned earlier, the tolerance level
for Newton's method was set at $tol<$1E$-$06.
%-----Test set --------------------------------
\begin{table}[ht!]
    \centering
    \rowcolors{1}{mygray}{white}
    %\scriptsize
    \caption{Test sets consist of grids with dimensions of $20 \times 20$ Chebyshev-distributed points. Error norms, including $L_1$, $L_2$, and $L_{\infty}$, quantify the disparities between the analytical and computed solutions of $T_1(x_1,x_2)$ and $T_2(x_1,x_2)$ using various training schemes. These schemes encompass training both with a greedy-wise approach and across the entire domain $[-0.91, 0] \times [-0.91, 0]$. The analysis of error norms includes statistics such as the median and 95th percentiles, aggregated over 400 runs, for each norm corresponding to both $T_1(x_1,x_2)$ and $T_2(x_1,x_2)$.}
    \begin{tabular}{||c|c| c c c ||}
    \toprule
    \toprule
    & Error  & Power-Series & PINN Entire domain & PINN Greedy \\
    $T(x_1,x_2)$& norm &$6th$ order &Median(5th,95th)  &Median(5th,95th) \\
    \midrule
    \midrule 
    &$\lVert \cdot \rVert_1$ & 5.19E$+$02
        &7.03E$+$01,(5.75E$+$01, 8.31E$+$01)  & 2.17E$+$00,(3.43E$-$01,  2.60E$+$01)\\
     $T_1(x_1,x_2)$& $\lVert \cdot \rVert_{2}$ & 4.13E$+$01&2.16E$+$01,(1.73E$+$01,2.62E$+$01)&3.55E$-$01,(4.52E$-$02,3.48E$+$00)\\
   & $\lVert \cdot \rVert_{\infty}$ &1.71E$+$01& 9.39E$+$00,(7.32E$+$00,1.16E$+$01)&3.37E$-$01,(1.78E$-$02,1.058E$+$00)\\
    \midrule
     &  $\lVert \cdot \rVert_1$ &5.87E$+$01
 & 2.14E$+$02,(1.82E$+$02, 2.15E$+$02) & 1.72E$+$01,(8.25E$-$01, 7.50E$+$01)\\
  $T_2(x_1,x_2)$ & $\lVert \cdot \rVert_{2}$
  &5.30E$+$01& 6.92E$+$01,(5.84E$+$01,8.04E$+$01)&2.02E$+$00,(1.11E$-$01,1.92E$+$01)\\
  &   $\lVert \cdot \rVert_{\infty}$ &4.40E$+$01& 3.40E$+$01,(2.89E$+$01,4.02E$+$01) & 1.90E$+$00,(5.36E$-$02, 3.99E$+$00) \\
    \bottomrule
     \bottomrule
    \end{tabular}
    \label{tab:NormsS2_Ex2}
\end{table}
\begin{table}[ht!]
    \scriptsize
    \centering
    \rowcolors{1}{mygray}{white}
    \caption{Benchmark problem 1. Uncertainty quantification systematically conducted on a test set employing a grid of $20 \times 20$ Chebyshev-Lobatto distributed points to investigate the convergence behavior of the Physics-Informed Neural Networks (PINN). The study involves performing 100, 200, 300, and 400 independent runs to thoroughly assess the variability in the results. The analysis focused on capturing the central tendency by calculating the median and examining the dispersion of outcomes through the 5th and 95th percentile confidence intervals.}
    \begin{tabular}{||c|c| c c c c ||}
    \toprule
    \toprule
   & Error  & 100 runs & 200 runs & 300 runs & 400 runs \\
    $T(x_1,x_2)$& norm &Median&Median &Median  &Median \\
    &&(5th,95th)&(5th,95th) &(5th,95th) &(5th,95th)
    \\
    \midrule
    \midrule
    & $\lVert \cdot \rVert_1$ &2.17E$+$00,&2.88E$+$00,&2.83E$+$00,&2.72E$+$00,
    \\
    &&(3.43E$-$01,2.60E$+$01)&(3.97E$-$01,  2.81E$+$01)&(3.86E$-$01, 2.721E$+$01)&(3.81E$-$01, 2.715E$+$01)
    \\
    $T_1(x_1,x_2)$ & $\lVert \cdot \rVert_{2}$ &3.55E$-$01,&4.23E$-$01,&4.01E$-$01,&3.98E$-$01,\\
    &&(4.52E$-$02,3.48E$+$00)&(5.12E$-$02,3.93E$+$00)&(5.076E$-$02,3.66E$+$00)&(4.93E$-$02,3.563E$+$00)
    \\
    & $\lVert \cdot \rVert_{\infty}$ &3.37E$-$01,&3.90E$-$01,&3.84E$-$01,&3.81E$-$01,\\
    &&(1.78E$-$02,1.0581E$+$00)&(1.69E$-$02,1.210E$+$00)&(1.03E$-$02,7.01E$-$01)& (1.00E$-$02,6.91E$-$01)
    \\
    \midrule
    & $\lVert \cdot \rVert_1$ &1.72E$+$01,&1.95E$+$01,&1.49E$+$01,&1.20E$+$01,\\
    &&(8.25E$-$01, 7.50E$+$01)&(8.76E$-$01, 7.85E$+$01)&(8.39E$-$01, 6.55E$+$01)&(8.32E$-$01, 6.13E$+$01)
    \\
   $T_2(x_1,x_2)$ & $\lVert \cdot \rVert_{2}$ &2.02E$+$00,&2.32E$+$00,&1.787E$+$00,&1.59E$+$00,
   \\
   &&(1.11E$-$01,1.92E$+$01)&(1.52E$-$01,2.05E$+$01)&(1.23E$-$01,1.786E$+$01)&(1.09E$-$01,1.778E$+$01)
    \\
   & $\lVert \cdot \rVert_\infty$ &1.90E$+$00,&1.92E$+$00,&1.72E$+$00,&1.49E$+$00, \\ 
   &&(5.36E$-$02, 3.99E$+$00)&(5.37E$-$02, 4.00E$+$00)&(5.64E$-$02, 3.21E$+$00)&(5.35E$-$02, 3.20E$+$00)
   \\
    \bottomrule
    \bottomrule
    \end{tabular}
    \label{tab:4}
\end{table}

\section{Conclusions}
\label{sec:conclusion}
In the present research study, a novel Physics-Informed Machine Learning (PINN) methodological and simulation framework is introduced to solve the nonlinear state observation problem in the discrete-time domain. This procedure stands in contrast to conventional methodologies \cite{Isidori1995},\cite{Luenberger1963}, as it is executed in a single step. The proposed PINN approach facilitates the acquisition of a nonlinear transformation map, effectively reconstructing the dynamics of the nonlinear system through linearization up to an additional output injection term in the transformed system. This transformation map is required to satisfy a system of first-order inhomogeneous functional equations. The principal advantages of this approach lie in its robustness in reconstructing system states, given a reasonably plausible set of assumptions. Regarding practical applications, it is worth noting that such transformation maps may exhibit abrupt gradients due to inherent singularities. Consequently, a greedy-wise training procedure has been implemented to refine our Physics-Informed Machine Learning (PINN) scheme's convergence characteristics. In this context, we employ a rudimentary zeroth-order continuation method to furnish improved initial estimations for the unknown PINN scheme parameters in proximity to singular regions. The incorporation of continuation techniques bears significant potential in simplifying computational experimentation and enhancing the comprehension of these singularities, possibly suggesting strategies for their alleviation. Moreover, we have demonstrated that our PINN scheme exhibits remarkable numerical precision, and thus outperforming established numerical solutions as power-series expansion. % Ultimately, we have performed uncertainty quantification on the approximation of the nonlinear transformation to test the convergence of our Greedy-wise Physics Informed Machine Learning approach to the real (analytical) transformations.

%% maybe this can be consider for the conclusion
\newpage
\appendix

\end{document}